\documentclass[final]{siamltex}
\usepackage{etex}
\reserveinserts{28}
\usepackage{amscd}
\usepackage{amsmath}
\usepackage{amsfonts}
\usepackage{amssymb}
\usepackage{cases}
\usepackage{t1enc}
\usepackage{array}
\usepackage[utf8]{inputenc}
\usepackage{graphicx}
\usepackage[usenames, svgnames]{xcolor}
\usepackage{pgf,tikz}
\usepackage{pgf,tikz}
\usetikzlibrary{matrix,arrows,decorations.pathmorphing,calc}
\usepackage{epstopdf}
\usepackage{showkeys}

\usepackage[ruled]{algorithm}            
\usepackage[matha,mathx,mathb]{mathabx}  
\usepackage{listings}

\usepackage{RR}      

\usepackage{datetime}

\usepackage{xkvltxp}
\usepackage[draft,multiuser,french,inlineface=\bfseries]{fixme}
\fxusetheme{colorsig}
\FXRegisterAuthor{fc}{afc}{FC} 
\FXRegisterAuthor{cj}{acj}{CJ} 
\FXRegisterAuthor{gs}{ags}{GS} 
\fxuselayouts{pdfnote}

\newcommand{\zoom}{0.36}
\newcommand{\zoomdeux}{0.35}
\graphicspath{{./images/}{./}}
\DeclareGraphicsExtensions{.eps}
\newcommand{\imageps}[2]{
  \begin{center}
    \hspace{-0.5cm}
    \includegraphics*[scale={#2}]{#1}
  \end{center}}
\newcommand{\dblimageps}[4]{
  \begin{center}
    \hspace{-0.5cm}
    \includegraphics*[scale={#2}]{#1}
    \includegraphics*[scale={#4}]{#3}
  \end{center}}

\newcommand{\dblimagepss}[4]{
  \begin{center}
    \includegraphics*[scale={#2}]{#1}
     \hspace{-0.8cm}
    \includegraphics*[scale={#4}]{#3}
  \end{center}}

\newtheorem{remark}[theorem]{Remark}

\newcommand{\dsp}{\displaystyle}
\def \R{\mathbb{R}}
\def\GRAD{\mathop{\rm \nabla}\nolimits}
\def\GrandO{\mathop{\mathcal{ O}}\nolimits}
\newcommand{\ENS}[2]{\left\{ {{#1},\hdots,{#2}}\right\}}

\newcommand{\vecb}[1]{\pmb{#1}}
\newcommand{\VEC}[1]{\vecb{#1}}                   
\newcommand{\MAT}[1]{\mathbb{#1}}
\newcommand{\MATT}[1]{\mathbb{#1}^{\text{\texttt{t}}}} 
\def\MatSet{\mathcal{M}}                          
\newcommand{\MS}[2]{\MatSet_{#1}(#2)} 
\newcommand{\DOT}[2]{%
\left\langle #1,#2 \right\rangle}

\newcommand{\DP}[2]{\frac{\partial #1}{\partial #2}}

\newcommand\Local[1]{\tilde{#1}}                  

\def\DOM{\Omega}
\def\DOMH{{\Omega_h}}
\def\nq{{\mathop{\rm n_q}\nolimits}}
\def\nqs{{\mathop{\rm n^2_q}\nolimits}}
\def\nme{{\mathop{\rm n_{me}}\nolimits}}
\def\q{{\rm q}}                                   
\def\qL{\Local{\q}}                               
\def\me{\mathop{\rm me}\nolimits}                 
\def\Taires{\mathop{\rm areas}\nolimits}
\def\areas{{\Taires}}
\def\ndf{{\mathop{\rm n_{df}}\nolimits}}
\def\ndfs{{\mathop{\rm n^2_{df}}\nolimits}}

\def\al{\mathop{\rm \nu}\nolimits}
\def\il{{\rm \alpha}}
\def\jl{{\rm \beta}}

\def\Masse{\mathbb{M}}
\def\MasseElem{\Masse^e}
\newcommand{\MasseF}[1]{\Masse^{[{#1}]}}
\newcommand{\MasseFElem}[1]{\Masse^{e,[{#1}]}}
\newcommand{\MassWElem}[1]{\MasseFElem{#1}}
\def\Stiff{\mathbb{S}}
\def\StiffElem{\mathbb{S}^e}

\def\FoncBase{\varphi}
\def\FoncBaseLocale{\Local{\FoncBase}}
\def\FoncBaseDeuxDLocale{\ASymbol{\Local{\vecb{\psi}}}} 

\def\HUnH{{H_h^1}}
\newcommand{\HUnHD}[1]{{{\HUnH}({#1})}}

\def\BasisFuncTwoD{\vecb{\psi}}                         
\def\BasisFuncTwoDL{\Local{\BasisFuncTwoD}}
\newcommand{\ASymbol}[1]{#1}                            
\def\StiffElas{\ASymbol{\mathbb{K}}}
\def\StiffElasElem{{\ASymbol{\mathbb{K}}}^e}
\def\Odc{\epsilon}
\def\Od{\vecb{\Odc}}                                    
\def\Odv{\underline{\Od}}                              
\def\Ocv{\mathop{\rm \vecb{\underline{\Occ}}}\nolimits}
\def\Occ{\sigma}

\lstnewenvironment{Listing}{}{%
\csname lst@SetFirstLabel\endcsname}
{\csname lst@SaveFirstLabel\endcsname}

\lstdefinelanguage{FreeFEM++}%
   {morekeywords={border,label,savemesh,mesh,buildmesh,plot,%
         cos,sin,exp,real,int,func,string,%
	 wait,ps},%
    sensitive=false,%
    morecomment=[l]//%
   }[keywords, comments]%

\newcommand{\MonFreeFEMNonumber}{
   \lstset{language=FreeFEM++,
	texcl=true,
	breaklines=true,
        basicstyle=\small,
        keywordstyle=\bfseries,
        stringstyle=\ttfamily,
        numberstyle=\tiny,
        stepnumber=1,
        numbersep=5pt,
	firstnumber=1,
        numbers=none,
        frame=bt}
}

\newcommand{\MonMatlabNonumber}{
\lstset{language=MATLAB,
        basicstyle=\small,
        keywordstyle=\bfseries,
        stringstyle=\ttfamily,
        numberstyle=\tiny,
        numbers=none,
       frame=bt
        }
}

\newcommand{\MonMatlab}{
\lstset{language=MATLAB,
        basicstyle=\small,
        keywordstyle=\bfseries,
        stringstyle=\ttfamily,
        numberstyle=\tiny,
        numbers=left,
       frame=bt
        }
}

\def\ListingWidth{0.8\textwidth}

\def\OptvsHanJunvsRahValDATE{2013_01_23_14_11}
\def\OptvsHanJunvsRahValDATEOctave{2013_01_23_17_35}

\colorlet{colorquiver}{DarkBlue}
\colorlet{matelem}{Maroon}
\colorlet{indelem}{Indigo}

\makeatletter
  \renewcommand*{\@fnsymbol}[1]{\ensuremath{\ifcase#1\or *\or
   \mathsection\or \mathparagraph\or \|\or **\or \dagger\dagger
   \or \ddagger\ddagger \else\@ctrerr\fi}}
\makeatother

\usepackage{hyperref,fancybox}
\makeatletter
     \@addtoreset{figure}{section}
     \@addtoreset{table}{section}
\makeatother

\title{An efficient way to perform the assembly of finite element matrices in Matlab and Octave}

\RRNo{8305}
\RRdate{May 2013}  

\RRauthor{
  François Cuvelier
    \thanks[sfn]{Universit\'e Paris 13,  Sorbonne Paris Cité, LAGA, CNRS UMR 7539, 99 Avenue J-B Cl\'ement, 93430 Villetaneuse, France,
            Emails: \texttt{cuvelier@math.univ-paris13.fr, japhet@math.univ-paris13.fr, scarella@math.univ-paris13.fr}}%
\and Caroline Japhet \footnotemark[1] 
    \thanks{INRIA Paris-Rocquencourt, project-team Pomdapi, 78153 Le Chesnay Cedex, France Email:\texttt{Caroline.Japhet@inria.fr}}%
\and Gilles Scarella\footnotemark[1]}%
\authorhead{Cuvelier \& Japhet \& Scarella}

\RRtitle{Optimisation de l'assemblage de matrices éléments finis sous Matlab et Octave}   
\RRetitle{An efficient way to perform the assembly of finite element matrices in Matlab and Octave}  

\RRnote{This work was partially supported by the GNR MoMaS (PACEN/CNRS, ANDRA, BRGM, CEA, EDF, IRSN)}

\RRresume{L'objectif est de décrire différentes techniques d'optimisation, sous Matlab/Octave, 
de routines d'assemblage de matrices éléments finis, en partant de l'approche classique 
jusqu'aux plus récentes vectorisées, sans utiliser de langage de bas niveau.
On aboutit au final à une version vectorisée rivalisant, en terme de performance, avec des
logiciels dédiés tels que FreeFEM++. Les descriptions des différentes méthodes d'assemblage
étant génériques, on les présente pour différentes matrices dans le cadre des éléments finis 
$P_1-$Lagrange en dimension $2$ et en élasticité linéaire. Des résultats numériques sont donnés pour 
illustrer les temps calculs des méthodes proposées.}

\RRabstract{
We describe different optimization techniques to perform the assembly
of finite element matrices in Matlab and Octave, from the standard
approach to recent vectorized ones, without any low level language
used.  We finally obtain a simple and efficient vectorized algorithm
able to compete in performance with dedicated software such as
FreeFEM++. The principle of this assembly algorithm is general, we
present it for different matrices in the $P_1$ finite elements case
and in linear elasticity.  We present numerical results which
illustrate the computational costs of the different approaches.}

\RRmotcle{éléments finis $P_1$, assemblage de matrices, vectorisation, élasticité linéaire, Matlab, Octave}

\RRkeyword{$P_1$ finite elements, matrix assembly, vectorization, linear elasticity, Matlab, Octave}

\RRprojets{Pomdapi}
\URRocq 
\RCParis 

\begin{document}

\makeRR

\renewcommand{\thefootnote}{\fnsymbol{footnote}}

\footnotetext[1]{Universit\'e Paris 13, Sorbonne Paris Cité, LAGA, CNRS UMR 7539,
  99 Avenue J-B Cl\'ement, 93430~Villetaneuse, France (\texttt{cuvelier@math.univ-paris13.fr, japhet@math.univ-paris13.fr},
 \texttt{scarella@math.univ-paris13.fr})}
\footnotetext[2]{INRIA Paris-Rocquencourt, project-team Pomdapi, 78153 Le Chesnay Cedex, France}

\renewcommand{\thefootnote}{\arabic{footnote}}

\pagestyle{myheadings}
\thispagestyle{plain}
\markboth{F. Cuvelier, C. Japhet, G. Scarella}{An efficient way to perform the assembly of
          finite element matrices in Matlab and Octave}

\MonMatlab
\section{Introduction}
Usually, finite elements methods \cite{Ciarlet:TFE:2002,Johnson:NSP:2009} are used to solve partial
differential equations (PDEs) occurring in many applications such as mechanics, fluid dynamics and computational
electromagnetics. These methods are based on a discretization of a weak formulation of the PDEs and
need the assembly of large sparse matrices (e.g. mass or stiffness matrices).
They enable complex geometries and various boundary conditions and they may be coupled with other discretizations,
using a weak coupling between different subdomains with nonconforming meshes~\cite{Bernardi:NNA:1989}.
Solving accurately these problems requires meshes containing a large number of elements and thus the assembly
of large sparse matrices.

Matlab \cite{Matlab:2012} and GNU Octave \cite{Octave:2012} are
efficient numerical computing softwares using matrix-based language
for teaching or industry calculations.  However, the classical
assembly algorithms (see for example~\cite{Cuvelier:MEF:2008,Lucquin:ISC:1998})
basically implemented in Matlab/Octave are much less efficient than
when implemented with other languages.

In~\cite{Davis:DMS:2006} Section~10, 
T. Davis describes different assembly techniques applied to random matrices of finite element type,
while the classical matrices are not treated. 
A first vectorization technique is proposed in~\cite{Davis:DMS:2006}.
Other more efficient algorithms have been proposed recently
in~\cite{Chen:PFE:2011,Chen:iFEM:2013,Hannukainen:IFE:2012,Rahman:FMA:2011}.
More precisely, in~\cite{Hannukainen:IFE:2012}, a vectorization is proposed, based on the permutation
of two local loops with the one through the elements. This more formal technique
allows to easily assemble different matrices, from  a reference element by affine transformation
and by using a numerical integration.
In~\cite{Rahman:FMA:2011}, the implementation is based on extending element operations on arrays
into operations on arrays of matrices, calling it a matrix-array operation, where the array elements are matrices
rather than scalars, and the operations are defined by the rules of linear algebra.
Thanks to these new tools and a quadrature formula, different matrices are computed 
without any loop.
In~\cite{Chen:iFEM:2013}, L.~Chen builds vectorially the nine sparse matrices corresponding to the 
nine elements of the element matrix and adds them to obtain the global matrix.

In this paper we present an optimization approach, in Matlab/Octave, using a vectorization of the algorithm.
This finite element assembly code is entirely vectorized (without loop) and without any quadrature formula.
Our vectorization is close to the one proposed in~\cite{Chen:PFE:2011},
with a full vectorization of the arrays of indices. 

Due to the length of the paper, we restrict ourselves to $P_1$ Lagrange finite elements in 2D
with an extension to linear elasticity.
Our method extends easily to the $P_k$ finite elements case, $k \ge 2$, and in 3D, see~\cite{Cuvelier:AEW:2013}.
We compare the performances of this code with the ones obtained with the standard algorithm and with those
proposed in~\cite{Chen:PFE:2011,Chen:iFEM:2013,Hannukainen:IFE:2012,Rahman:FMA:2011}.
We also show that this implementation is able to compete in performance
with dedicated software such as FreeFEM++~\cite{Hecht:FreeFEM:2012}.
All the computations are done on our reference
computer\footnote{2 x Intel Xeon E5645 (6 cores) at 2.40Ghz, 32Go RAM, supported by GNR MoMaS}
with the releases R2012b for Matlab, 3.6.3 for Octave
and 3.20 for FreeFEM++.
The entire Matlab/Octave code may be found in~\cite{Cuvelier:OptFEM2DP1:2012}.
The Matlab codes are fully compatible with Octave.

The remainder of this paper is organized as follows: in Section~\ref{sssec:PositionduProbleme} we give the notations
associated to the mesh and we define three finite element matrices. Then, in Section~\ref{sssec:AlgoClassique} we
recall the classical algorithm to perform the assembly of these matrices and show its inefficiency compared to
FreeFEM++. This is due to the storage of sparse matrices in Matlab/Octave as explained in
Section~\ref{sssec:StockageMatCreuse}.
In Section~\ref{sec:CodeMassOptV1} we give a method to best use  Matlab/Octave \texttt{sparse} function,
the ``optimized version 1'', suggested in~\cite{Davis:DMS:2006}. Then, in Section~\ref{sssec:CodeMassOptV2}
we present a new vectorization approach, the ``optimized version 2'', and compare its performances
to those obtained with FreeFEM++ and the codes given
in~\cite{Chen:PFE:2011,Chen:iFEM:2013,Hannukainen:IFE:2012,Rahman:FMA:2011}.
Finally, in Section~\ref{sec:Elasticity}, we present an extension to linear elasticity.
The full listings of the routines used in the paper are given in Appendix~\ref{App:codes} (see also \cite{Cuvelier:OptFEM2DP1:2012}).

\section{Notations}
\label{sssec:PositionduProbleme}
\begin{figure}[H]
  \centering
\imageps{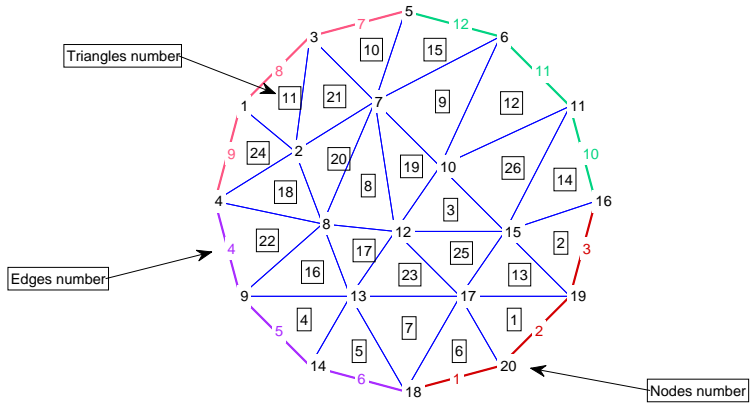}{0.55}
\caption{Description of the mesh.
  \label{mesh}}
\end{figure}

Let $\Omega$ be an open bounded subset of $\mathbb{R}^2$. It is provided with
its mesh ${\cal T}_h$ (classical and locally conforming).
We use a triangulation $\DOMH=\bigcup_{T_k \in {\cal T}_h} T_k$ of $\DOM$ (see Figure~\ref{mesh}) described by :
\vspace{-1mm}
\begin{center}
\begin{tabular}{lccl}
\hline
\textbf{name} & \textbf{type} & \textbf{dimension} & \textbf{description}\\
\hline
$\nq$ & integer & 1 & number of vertices\\
$\nme$ & integer & 1 & number of elements\\
$\q$   & double  &$2\times \nq$ &
\begin{minipage}[t]{7.9cm}
array of vertices coordinates.~${\q}(\al,j)$ is the $\al$-th coordinate of the $j$-th vertex,
$\al\in\{1,2\}$, $j\in\{1,\hdots,\nq\}.$
The $j$-th vertex will be also denoted by $\q^j$ with $\q^j_x=\q(1,j)$ and $\q^j_y=\q(2,j)$
\end{minipage}\\
$\me$  & integer & $3\times \nme$ &  
\begin{minipage}[t]{7.9cm}
connectivity array. $\me(\jl,k)$ is the storage index of the $\jl$-th vertex
of the $k$-th triangle, in the array~$q$, for $\jl\in\{1,2,3\}$ and $k\in\{1,\hdots,\nme\}$
\end{minipage}\\
$\areas$  & double & $1\times \nme$ &
\begin{minipage}[t]{7.9cm}
array of areas. $\areas(k)$ is the $k$-th triangle area, $k\in\{1,\hdots,\nme\}$
\end{minipage}\\
\hline
\end{tabular}
\end{center}
\indent
In this paper we will consider the assembly of the mass, weighted mass and stiffness matrices
denoted by $\Masse,$ $\MasseF{w}$ and $\Stiff$ respectively. These matrices of size $\nq$ are sparse, and
their coefficients are defined by
{\small
$$
\Masse_{i,j}=\int_\DOMH \FoncBase_i(\q) \FoncBase_j(\q) d\q,\ \ \MasseF{w}_{i,j}
=\int_\DOMH w(\q)\FoncBase_i(\q) \FoncBase_j(\q) d\q,\ 
\ 
\Stiff_{i,j}=\int_\DOMH \DOT{\GRAD\FoncBase_i(\q)}{\GRAD\FoncBase_j(\q)}d\q,
$$
}%
 where $\FoncBase_i$ are the usual $P_1$ Lagrange basis functions,
$w$ is a function defined on $\DOM$ and
$\DOT{\cdot}{\cdot}$ is the usual scalar product in $\mathbb{R}^2$. 
More details are given in \cite{Cuvelier:MEF:2008}.
To assemble these matrices, one needs to compute its associated element matrix.
On a triangle $T$ with local vertices $\qL^1,$ $\qL^2,$ $\qL^3$ and area $|T|$,
the element mass matrix is given by 
\begin{equation}\label{MassElem}
\MasseElem(T)=\frac{|T|}{12}\begin{pmatrix} 2 & 1 & 1 \\ 1 & 2 & 1\\ 1 & 1 & 2\end{pmatrix}.
\end{equation}
Let $\Local{w}_\il=w(\qL^{\il}),$ $\forall \il\in\ENS{1}{3}$.
The element weighted mass matrix is approximated by
\begin{equation}\label{MassWElem}
\MassWElem{\Local{w}}(T)=\frac{|T|}{30}
\begin{pmatrix}
   3\Local{w}_1 + \Local{w}_2 +\Local{w}_3  & 
   \Local{w}_1+\Local{w}_2+\frac{\Local{w}_3}{2} & 
   \Local{w}_1+\frac{\Local{w}_2}{2}+\Local{w}_3\\
   \Local{w}_1+\Local{w}_2+\frac{\Local{w}_3}{2}& 
   \Local{w}_1 + 3\Local{w}_2 +\Local{w}_3& 
   \frac{\Local{w}_1}{2}+\Local{w}_2+\Local{w}_3\\
   \Local{w}_1+\frac{\Local{w}_2}{2}+\Local{w}_3& 
   \frac{\Local{w}_1}{2}+\Local{w}_2+\Local{w}_3& 
   \Local{w}_1 + \Local{w}_2 +3\Local{w}_3
\end{pmatrix}.
\end{equation}
Denoting $\vecb{u}=\qL^2-\qL^3,$ $\vecb{v}=\qL^3-\qL^1$ and $\vecb{w}=\qL^1-\qL^2,$ 
the element stiffness matrix is 
\begin{equation}\label{StiffMatElem}
\StiffElem(T) \\=\\ \frac{1}{4|T|}
\displaystyle\begin{pmatrix}
   \DOT{\vecb{u}}{\vecb{u}} &\DOT{\vecb{u}}{\vecb{v}} &\DOT{\vecb{u}}{\vecb{w}}\\
   \DOT{\vecb{v}}{\vecb{u}} &\DOT{\vecb{v}}{\vecb{v}} &\DOT{\vecb{v}}{\vecb{w}}\\
   \DOT{\vecb{w}}{\vecb{u}} &\DOT{\vecb{w}}{\vecb{v}} &\DOT{\vecb{w}}{\vecb{w}}
\end{pmatrix}.
\end{equation}
The listings of the routines to compute the previous element matrices are given in Appendix~\ref{App:ElemMat}
We now give the classical assembly algorithm using these element matrices
with a loop through the triangles.

\section{The classical algorithm}
\label{sssec:AlgoClassique}
\begin{figure}[H]
\dblimageps{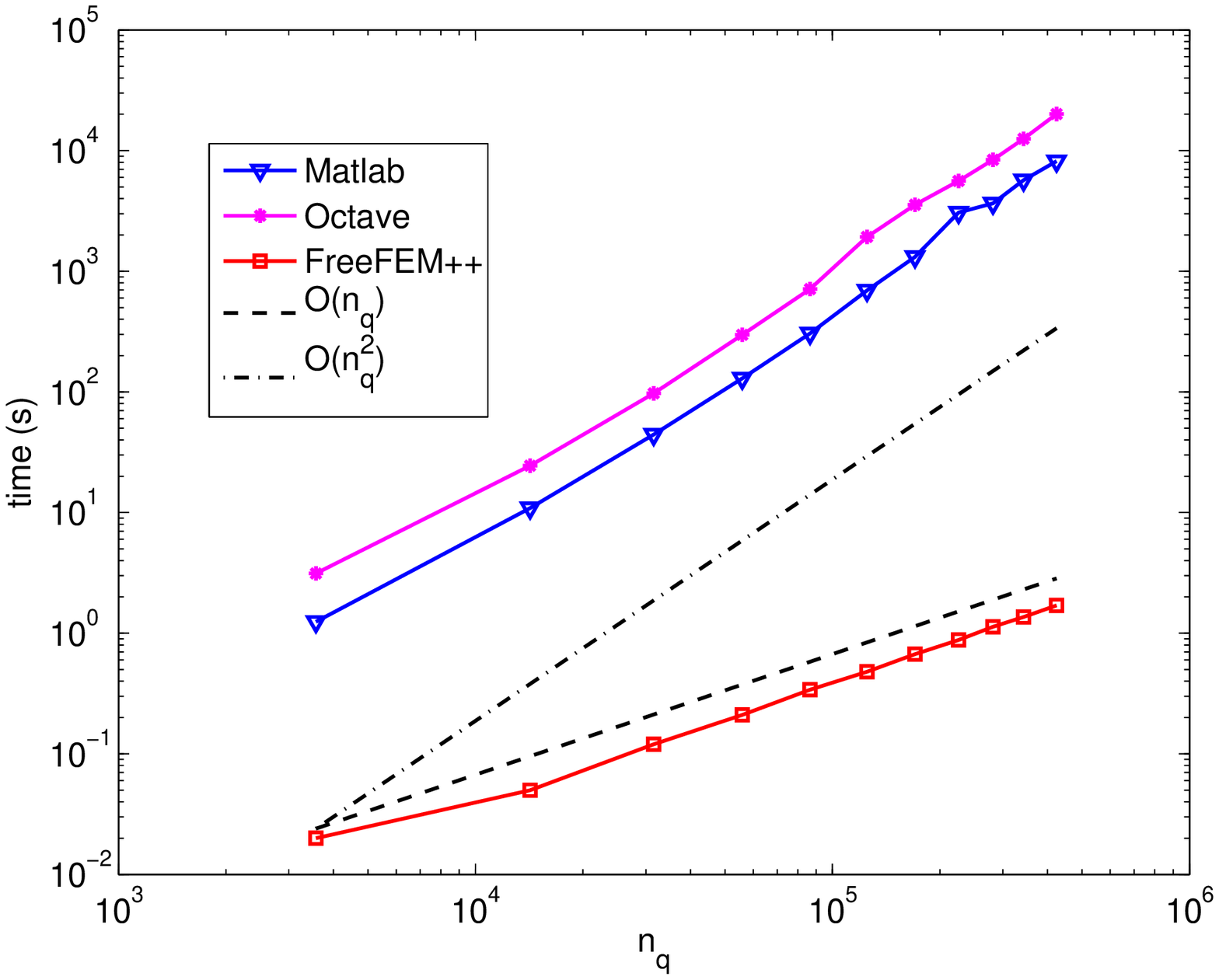}{\zoom}{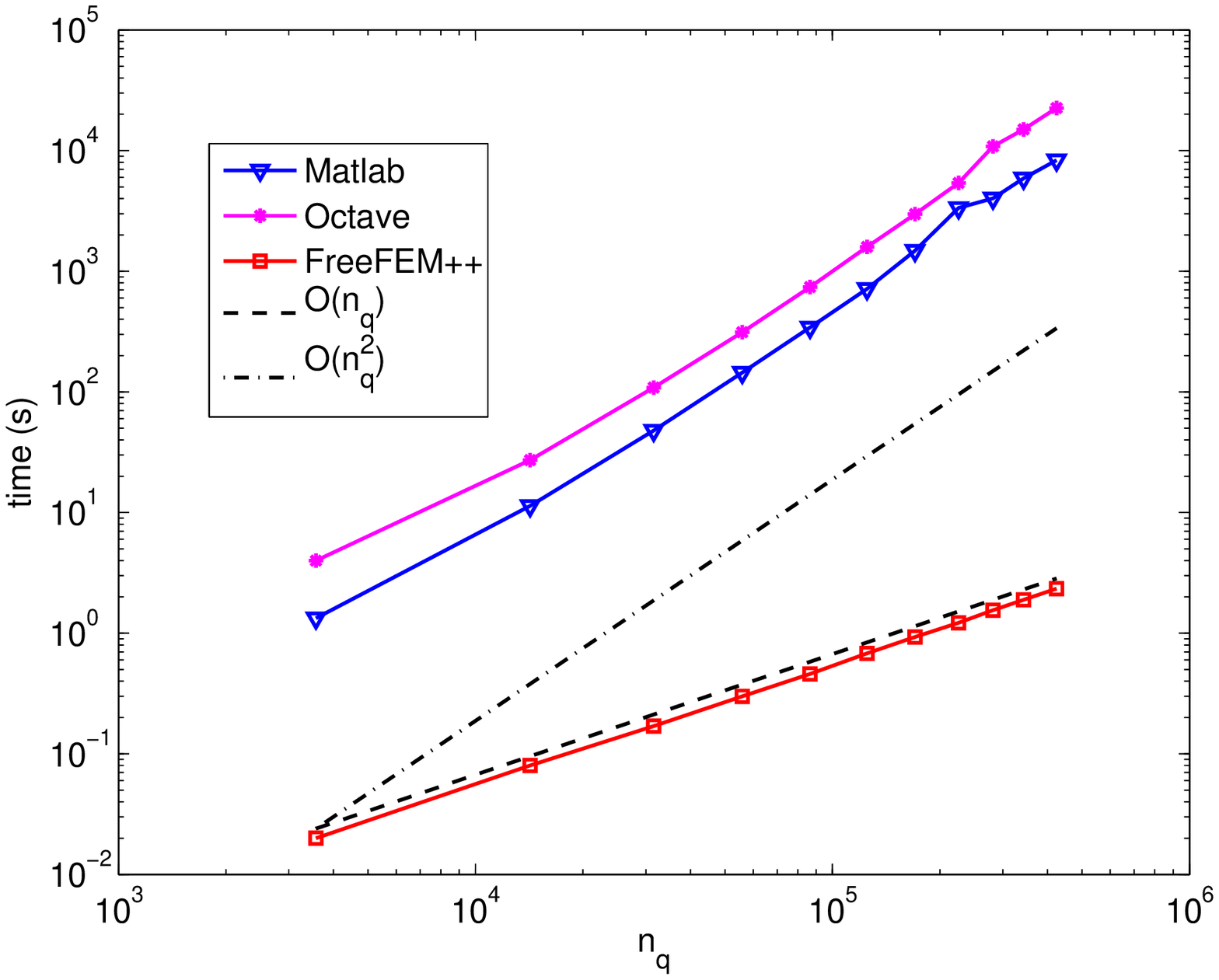}{\zoom}
\imageps{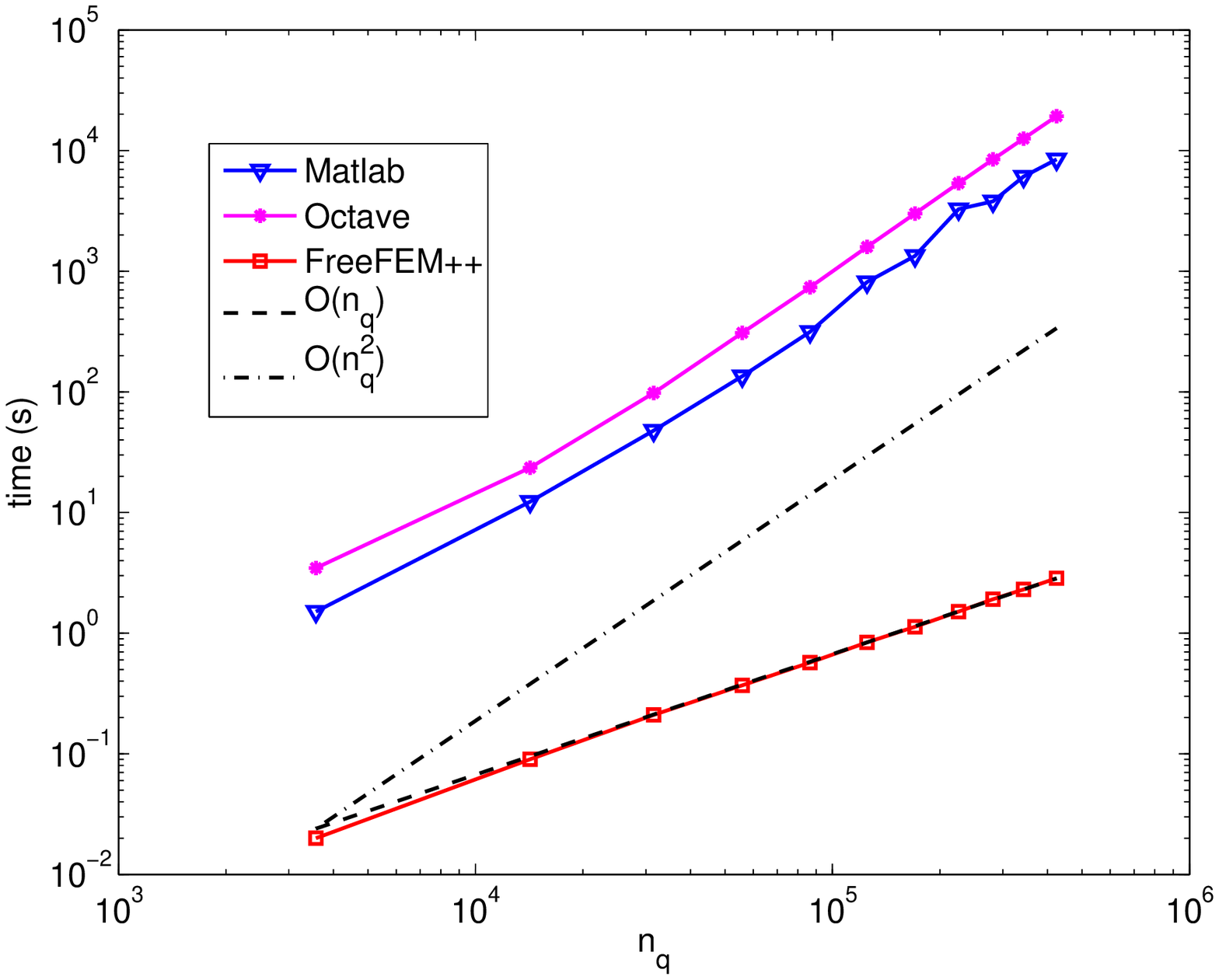}{\zoom}  
\caption{Comparison of the classical matrix assembly code in Matlab/Octave with FreeFEM++,
  for the mass (top left), weighted mass (top right) and stiffness (bottom) matrices. \label{BaseVsFreeFEM}}
\end{figure}

We describe the assembly of a given $\nq\times\nq$ matrix \lstinline{M} from its
associated $3\times 3$ element matrix \lstinline{E}. We denote by
``\lstinline{ElemMat}'' the routine which computes the element matrix
\lstinline{E}.
\vspace{-1.5mm}
\begin{center}
\begin{minipage}[t]{\ListingWidth}
\MonMatlabNonumber
\begin{lstlisting}[caption=Classical matrix assembly code in Matlab/Octave,label=Assemblage:basique]
M=sparse(nq,nq);
for k=1:nme
    E=ElemMat(areas(k),...);
    for il=1:3
        i=me(il,k);
        for jl=1:3
            j=me(jl,k);
            M(i,j)=M(i,j)+E(il,jl);
        end
    end
end
\end{lstlisting}
\end{minipage}
\end{center}
We aim to compare the performances of this code (see Appendix~\ref{App:codebase} for the complete listings)
with those obtained with FreeFEM++~\cite{Hecht:FreeFEM:2012}.
The FreeFEM++ commands to build the mass, weighted mass and stiffness matrices are given
in Listing~\ref{codeAssembeFreeFEM}.
On Figure~\ref{BaseVsFreeFEM}, we show the computation times (in seconds) versus the number of vertices $\nq$
of the mesh (unit disk), for the classical assembly and FreeFEM++ codes.
The values of the computation times  are given in Appendix~\ref{App:basevsFreeFem}.
We observe that the complexity is $\GrandO(\nqs)$ (quadratic) for the Matlab/Octave codes,
while the complexity seems to be $\GrandO(\nq)$ (linear) for FreeFEM++.
\vspace{-1.5mm}
\begin{center}
\begin{minipage}[t]{0.83\textwidth}
\MonFreeFEMNonumber
\begin{lstlisting}[caption=Matrix assembly code in FreeFEM++,label=codeAssembeFreeFEM]
  mesh Th(...);
  fespace Vh(Th,P1); // P1 FE-space
  varf vMass (u,v)= int2d(Th)( u*v);
  varf vMassW (u,v)= int2d(Th)( w*u*v);
  varf vStiff (u,v)= int2d(Th)( dx(u)*dx(v)
                                 + dy(u)*dy(v) );
  
  matrix M= vMass(Vh,Vh);    // Mass matrix assembly
  matrix Mw = vMassW(Vh,Vh); // Weighted mass matrix assembly
  matrix S = vStiff(Vh,Vh);  // Stiffness matrix assembly
\end{lstlisting}
\end{minipage}
\end{center}
We have surprisingly observed that the Matlab performances may be improved using an older Matlab release
(see Appendix~\ref{MatlabTrouble})

Our objective is to propose optimizations of the classical code that lead to more efficient codes
with computational costs comparable to those obtained with FreeFEM++.
A first improvement of the classical algorithm (Listing~\ref{Assemblage:basique})
is to vectorize the two local loops, see Listing~\ref{Assemblage_OptV0}
(the complete listings are given in Appendix~\ref{App:codeV0}).
\begin{center}
\begin{minipage}[t]{\ListingWidth}
\MonMatlabNonumber
\begin{lstlisting}[caption=Optimized matrix assembly code - version 0,label=Assemblage_OptV0]
M=sparse(nq,nq);
for k=1:nme
  I=me(:,k);
  M(I,I)=M(I,I)+ElemMat(areas(k),...);
end
\end{lstlisting}
\end{minipage}
\end{center}
\noindent\\
However the complexity of this algorithm is still quadratic (i.e. $\GrandO(\nqs)$). 

In the next section, we explain the storage of sparse matrices in Matlab/Octave in order to justify this
lack of efficiency.

\section{Sparse matrices storage}
\label{sssec:StockageMatCreuse}
\noindent
In Matlab/Octave, a sparse matrix $\MAT{A}\in\MS{M,N}{\R}$, with $nnz$
non-zeros elements, is stored with CSC (Compressed Sparse Column)
format using the following three arrays:
\begin{eqnarray}
  \begin{array}{rcl}
aa(1:nnz) &:& \mbox{which contains the } nnz \mbox{ non-zeros elements of } \MAT{A} \mbox{ stored column-wise,}\\
ia(1:nnz) &:& \mbox{which contains the row numbers of the elements stored in } aa,\\
ja(1:N+1) &:& \mbox{which allows to find the elements of a column of } \MAT{A}, \mbox{ with the infor-}\\
&&\hspace{-1.5cm}\mbox{mation
that the first non-zero element of the column } k \mbox{ of } \MAT{A} \mbox{ is in the } ja(k)\mbox{-th}\\
&& \hspace{-1.5cm} \mbox{position in the
  array } aa. \mbox{ We have } ja(1)=1 \mbox{ and } ja(N+1)=nnz+1. 
\end{array}\nonumber
\end{eqnarray}
For example, with the matrix 
$$
\MAT{A}=\begin{pmatrix}
1. & 0. & 0. & 6.\\
0. & 5. & 0. & 4.\\
0. & 1. & 2. & 0.
\end{pmatrix},$$
we have $M=3,$ $N=4,$ $nnz=6$ and
$$\begin{array}{ll}
aa & 
\begin{array}{|*{6}{>{$}m{0.5cm}<{$}|}}
\hline
1. &  5. & 1. & 2.& 6. & 4.\\
\hline
\end{array}
\\ \\
ia & \begin{array}{|*{6}{>{$}m{0.5cm}<{$}|}}
\hline
1 &  2 & 3 & 3 & 1 & 2\\
\hline
\end{array}
\\ \\
ja & \begin{array}{|*{5}{>{$}m{0.5cm}<{$}|}}
\hline
1 &  2 &  4 & 5 & 7\\
\hline
\end{array}
\end{array}
$$
The first non-zero element in column $k=3$ of $\MAT{A}$ is $2$, the position of this number
in $aa$ is $4$, thus $ja(3)=4$.

We now describe the operations to be done on the arrays $aa,$ $ia$ and $ja$ if we modify
the matrix $\MAT{A}$ by taking $\MAT{A}(1,2)=8.$
It becomes
$$
\MAT{A}=\begin{pmatrix}
1. & \colorbox{gray!25}{8.} & 0. & 6.\\
0. & 5. & 0. & 4.\\
0. & 1. & 2. & 0.
\end{pmatrix}.
$$
In this case, a zero element of $\MAT{A}$ has been replaced by
the non-zero value $8$ which must be stored in the arrays while no space is provided.
We suppose that the arrays are sufficiently large (to avoid memory space problems),
we must then shift one cell all the values in the arrays $aa$ and $ia$ from the third position
and then copy the value~$8$ in $aa(3)$ and the value $1$ (row number) in $ia(3)$ :
$$
\begin{array}{ll}
aa & 
\begin{array}{|*{7}{>{$}m{0.5cm}<{$}|}}
\hline
1. &  \colorbox{gray!25}{8.} & 5. & 1. & 2. & 6. & 4.\\
\hline
\end{array}
\\ \\
ia & \begin{array}{|*{7}{>{$}m{0.5cm}<{$}|}}
\hline
1 &  \colorbox{gray!25}{1} & 2 & 3 & 3& 1 & 2\\
\hline
\end{array}
\end{array}
$$
For the array $ja,$ we increment of $1$ the values after the position $2$ :
$$
\begin{array}{ll}
ja & \begin{array}{|*{5}{>{$}m{0.5cm}<{$}|}}
\hline
1 & 2 &  \colorbox{gray!25}{5} & \colorbox{gray!25}{6} & \colorbox{gray!25}{8}\\
\hline
\end{array}
\end{array}
$$
The repetition of these operations is expensive upon assembly of the matrix in the previous codes.
Moreover, we haven't considered dynamic reallocation problems that may also occur.

We now present the optimized version 1 of the code that will allow to improve the performance of the
classical code.

\section{Optimized matrix assembly - version 1 (OptV1)}
\label{sec:CodeMassOptV1}
We will use the following call of the \texttt{sparse} Matlab function:\\
\centerline{\lstinline{M = sparse(I,J,K,m,n);}}
This command returns a sparse matrix $M$ of size \lstinline{m}
$\times$ \lstinline{n} such that \lstinline{M(I(k),J(k)) = K(k)}.  The
vectors \lstinline{I}, \lstinline{J} and \lstinline{K} have the same
length.  The zero elements of \lstinline{K} are not taken into account
and the elements of \lstinline{K} having the same indices in
\lstinline{I} and \lstinline{J} are summed.

The idea is to create three global 1d-arrays $\vecb{I}_g$, $\vecb{J}_g$ and $\vecb{K}_g$ 
allowing the storage of the element matrices as well as the position of their elements in the global matrix.
The length of each array is $9 \nme.$ 
Once these arrays are created, the matrix assembly is obtained with the command\\
\centerline{\lstinline{M = sparse(Ig,Jg,Kg,nq,nq);}}

To create these three arrays, we first define three local arrays $\vecb{K}^e_k,$ $\vecb{I}^e_k$ and 
$\vecb{J}^e_k$ of nine elements obtained from a generic element matrix $\MAT{E}(T_k)$ of dimension $3$ :\\
\begin{center}
\begin{tabular}{lcl}
$\vecb{K}^e_k$ &:& elements of the matrix $\MAT{E}(T_k)$ stored column-wise,\\
$\vecb{I}^e_k$ &:& global row indices associated to the elements stored in $\vecb{K}^e_k$,\\
$\vecb{J}^e_k$ &:& global column indices associated to the elements stored in $\vecb{K}^e_k.$
\end{tabular}
\end{center}
We have chosen a column-wise numbering for 1d-arrays in Matlab/Octave implementation, 
but for representation convenience we draw them in line format,
{\scriptsize
$$
\MAT{E}(T_k)=\begin{pmatrix}
\color{matelem} e_{1,1}^k & \color{matelem} e_{1,2}^k & \color{matelem} e_{1,3}^k\\
\color{matelem} e_{2,1}^k & \color{matelem} e_{2,2}^k & \color{matelem} e_{2,3}^k\\
\color{matelem} e_{3,1}^k & \color{matelem} e_{3,2}^k & \color{matelem} e_{3,3}^k
\end{pmatrix}\ \Longrightarrow
\begin{array}{rl}
\vecb{K}^e_k  :
&
\hspace{-1.5mm}
\begin{array}{|*{9}{>{$}m{0.5cm}<{$}|}}
\hline
\color{matelem} e_{1,1}^k & \color{matelem} e_{2,1}^k & \color{matelem} e_{3,1}^k &
\color{matelem} e_{1,2}^k & \color{matelem} e_{2,2}^k & \color{matelem} e_{3,2}^k &
\color{matelem} e_{1,3}^k & \color{matelem} e_{2,3}^k & \color{matelem} e_{3,3}^k \\
\hline
\end{array}
\\ \\
\vecb{I}^e_k:
&
\hspace{-1.5mm}
\begin{array}{|*{9}{>{$}m{0.5cm}<{$}|}}
\hline
\color{indelem} i_1^k & \color{indelem} i_2^k & \color{indelem} i_3^k &
\color{indelem} i_1^k & \color{indelem} i_2^k & \color{indelem} i_3^k &
\color{indelem} i_1^k & \color{indelem} i_2^k & \color{indelem} i_3^k \\
\hline
\end{array}
\\ \\
\vecb{J}^e_k:
&
\hspace{-1.5mm}
\begin{array}{|*{9}{>{$}m{0.5cm}<{$}|}}
\hline
\color{indelem} i_1^k  & \color{indelem} i_1^k & \color{indelem} i_1^k &
\color{indelem} i_2^k & \color{indelem} i_2^k & \color{indelem} i_2^k &
\color{indelem} i_3^k & \color{indelem} i_3^k & \color{indelem} i_3^k \\
\hline
\end{array}
\\
\end{array}
$$
}
with $i_1^k=\me(1,k),$  $i_2^k=\me(2,k),$ $i_3^k=\me(3,k).$

To create the three arrays $\vecb{K}^e_k,$ $\vecb{I}^e_k$ and $\vecb{J}^e_k,$ in Matlab/Octave,
one can use the following commands :
\begin{center}
\begin{minipage}[t]{\ListingWidth}
\MonMatlabNonumber
\begin{lstlisting}
E  = ElemMat(areas(k), ...);     % E  : 3-by-3 matrix 
Ke = E(:);                       % Ke : 9-by-1 matrix 
Ie = me([1 2 3 1 2 3 1 2 3],k);  % Ie : 9-by-1 matrix 
Je = me([1 1 1 2 2 2 3 3 3],k);  % Je : 9-by-1 matrix 
\end{lstlisting}
\end{minipage}
\end{center}

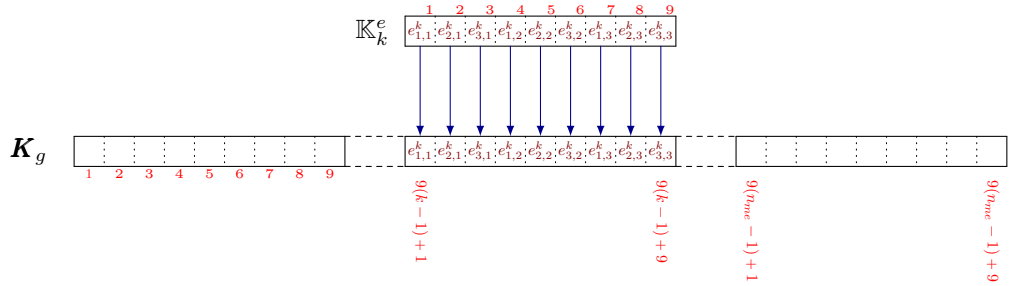
\begin{figure}[H]
\definecolor{ffqqff}{rgb}{1,0,1}
\definecolor{ffqqqq}{rgb}{1,0,0}
\definecolor{qqqqff}{rgb}{0,0,1}
\definecolor{cqcqcq}{rgb}{0.75,0.75,0.75}
\colorlet{colorquiver}{DarkBlue}
\colorlet{matelem}{Maroon}
\colorlet{indelem}{Indigo}

\hspace{-.2cm}
\begin{tikzpicture}[line cap=round,line join=round,>=triangle 45,x=1.0cm,y=1.0cm,>=latex,scale=0.4]
\draw (-4,0)-- (5,0);
\draw (-4,-1)-- (5,-1);
\draw (5,0)-- (5,-1);
\draw (-4,0)-- (-4,-1);
\draw [dotted] (-3,0)-- (-3.,-1);
\draw [dotted] (-2,0)-- (-2.,-1);
\draw [dotted] (-1,0)-- (-1.,-1);
\draw [dotted] (0,0)-- (0,-1);
\draw [dotted] (1,0)-- (1,-1);
\draw [dotted] (2,0)-- (2,-1);
\draw [dotted] (3,0)-- (3,-1);
\draw [dotted] (4,0)-- (4,-1);
\draw (7,4)-- (16,4);
\draw (7,3)-- (16,3);
\draw (16,4)-- (16,3);
\draw (7,4)-- (7,3);
\draw [dotted] (8,4)-- (8,3);
\draw [dotted] (9,4)-- (9,3);
\draw [dotted] (10,4)-- (10,3);
\draw [dotted] (11,4)-- (11,3);
\draw [dotted] (12,4)-- (12,3);
\draw [dotted] (13,4)-- (13,3);
\draw [dotted] (14,4)-- (14,3);
\draw [dotted] (15,4)-- (15,3);
\draw (7,0)-- (16,0);
\draw (7,-1)-- (16,-1);
\draw (16,0)-- (16,-1);
\draw (7,0)-- (7,-1);
\draw [dotted] (8,0)-- (8,-1);
\draw [dotted] (9,0)-- (9,-1);
\draw [dotted] (10,0)-- (10,-1);
\draw [dotted] (11,0)-- (11,-1);
\draw [dotted] (12,0)-- (12,-1);
\draw [dotted] (13,0)-- (13,-1);
\draw [dotted] (14,0)-- (14,-1);
\draw [dotted] (15,0)-- (15,-1);
\draw (18,0)-- (27,0);
\draw (18,-1)-- (27,-1);
\draw (27,0)-- (27,-1);
\draw (18,0)-- (18,-1);
\draw [dotted] (19,0)-- (19,-1);
\draw [dotted] (20,0)-- (20,-1);
\draw [dotted] (21,0)-- (21,-1);
\draw [dotted] (22,0)-- (22,-1);
\draw [dotted] (23,0)-- (23,-1);
\draw [dotted] (24,0)-- (24,-1);
\draw [dotted] (25,0)-- (25,-1);
\draw [dotted] (26,0)-- (26,-1);
\draw [dash pattern=on 2pt off 2pt] (5,0)-- (7,0);
\draw [dash pattern=on 2pt off 2pt] (5,-1)-- (7,-1);
\draw [dash pattern=on 2pt off 2pt] (16,0)-- (18,0);
\draw [dash pattern=on 2pt off 2pt] (16,-1)-- (18,-1);
\draw [->,color=colorquiver] (7.5,3) -- (7.5,0);
\draw [->,color=colorquiver] (8.5,3) -- (8.5,0);
\draw [->,color=colorquiver] (9.5,3) -- (9.5,0);
\draw [->,color=colorquiver] (10.5,3) -- (10.5,0);
\draw [->,color=colorquiver] (11.5,3) -- (11.5,0);
\draw [->,color=colorquiver] (12.5,3) -- (12.5,0);
\draw [->,color=colorquiver] (13.5,3) -- (13.5,0);
\draw [->,color=colorquiver] (14.5,3) -- (14.5,0);
\draw [->,color=colorquiver] (15.5,3) -- (15.5,0);
\draw (5.04,4.2) node[anchor=north west] {$ \mathbb{K}^e_k$ };
\draw (-5.5,-0.5) node[scale=1] {$\vecb{K}_g$};
\draw [color=ffqqqq](-3.5,-1.25) node[anchor=center,scale=0.5] {$ 1 $};
\draw [color=ffqqqq](-2.5,-1.25) node[anchor=center] {\tiny{$ 2 $}};
\draw [color=ffqqqq](-1.5,-1.25) node[anchor=center] {\tiny{$ 3 $}};
\draw [color=ffqqqq](-0.5,-1.25) node[anchor=center] {\tiny{$ 4 $}};
\draw [color=ffqqqq](0.5,-1.25) node[anchor=center] {\tiny{$ 5 $}};
\draw [color=ffqqqq](1.5,-1.25) node[anchor=center] {\tiny{$ 6 $}};
\draw [color=ffqqqq](2.5,-1.25) node[anchor=center] {\tiny{$ 7 $}};
\draw [color=ffqqqq](3.5,-1.25) node[anchor=center] {\tiny{$ 8 $}};
\draw [color=ffqqqq](4.5,-1.25) node[anchor=center] {\tiny{$ 9 $}};
\draw [color=ffqqqq](7.5,-1.25) node[right,rotate=-90,scale=0.6] {$ 9(k-1)+1 $};
\draw [color=ffqqqq](15.5,-1.25) node[right,rotate=-90,scale=0.6] {$ 9(k-1)+9$};
\draw [color=ffqqqq](26.5,-1.25) node[right,rotate=-90,scale=0.6] {$ 9(n_{me}-1)+9$};
\draw [color=ffqqqq](18.5,-1.25) node[right,rotate=-90,scale=0.6] {$ 9(n_{me}-1)+1$};
\draw [color=ffqqqq](7.36,4.68) node[anchor=north west] {\tiny{$ 1 $}};
\draw [color=ffqqqq](8.36,4.68) node[anchor=north west] {\tiny{$ 2 $}};
\draw [color=ffqqqq](9.36,4.68) node[anchor=north west] {\tiny{$ 3 $}};
\draw [color=ffqqqq](10.38,4.68) node[anchor=north west] {\tiny{$ 4 $}};
\draw [color=ffqqqq](11.4,4.68) node[anchor=north west] {\tiny{$ 5 $}};
\draw [color=ffqqqq](12.36,4.68) node[anchor=north west] {\tiny{$ 6 $}};
\draw [color=ffqqqq](13.38,4.68) node[anchor=north west] {\tiny{$ 7 $}};
\draw [color=ffqqqq](14.34,4.68) node[anchor=north west] {\tiny{$ 8 $}};
\draw [color=ffqqqq](15.34,4.68) node[anchor=north west] {\tiny{$ 9 $}};
\draw [color=matelem](7.5,3.5) node[anchor=center,scale=0.6] {$ e_{1,1}^k $};
\draw [color=matelem](8.5,3.5) node[anchor=center,scale=0.6] {$ e_{2,1}^k $};
\draw [color=matelem](9.5,3.5) node[anchor=center,scale=0.6] {$ e_{3,1}^k $};
\draw [color=matelem](10.5,3.5) node[anchor=center,scale=0.6] {$ e_{1,2}^k $};
\draw [color=matelem](11.5,3.5) node[anchor=center,scale=0.6] {$ e_{2,2}^k $};
\draw [color=matelem](12.5,3.5) node[anchor=center,scale=0.6] {$ e_{3,2}^k $};
\draw [color=matelem](13.5,3.5) node[anchor=center,scale=0.6] {$ e_{1,3}^k $};
\draw [color=matelem](14.5,3.5) node[anchor=center,scale=0.6] {$ e_{2,3}^k $};
\draw [color=matelem](15.5,3.5) node[anchor=center,scale=0.6] {$ e_{3,3}^k $};
\draw [color=matelem](7.5,-0.5) node[anchor=center,scale=0.6] {$ e_{1,1}^k $};
\draw [color=matelem](8.5,-0.5) node[anchor=center,scale=0.6] {$ e_{2,1}^k $};
\draw [color=matelem](9.5,-0.5) node[anchor=center,scale=0.6] {$ e_{3,1}^k $};
\draw [color=matelem](10.5,-0.5) node[anchor=center,scale=0.6] {$ e_{1,2}^k $};
\draw [color=matelem](11.5,-0.5) node[anchor=center,scale=0.6] {$ e_{2,2}^k $};
\draw [color=matelem](12.5,-0.5) node[anchor=center,scale=0.6] {$ e_{3,2}^k $};
\draw [color=matelem](13.5,-0.5) node[anchor=center,scale=0.6] {$ e_{1,3}^k $};
\draw [color=matelem](14.5,-0.5) node[anchor=center,scale=0.6] {$ e_{2,3}^k $};
\draw [color=matelem](15.5,-0.5) node[anchor=center,scale=0.6] {$ e_{3,3}^k $};

\end{tikzpicture}  
\caption{Insertion of an element matrix\label{AssMatV1} in the global array - Version 1}
\end{figure}

\begin{figure}[H]
  \centering
\dblimageps{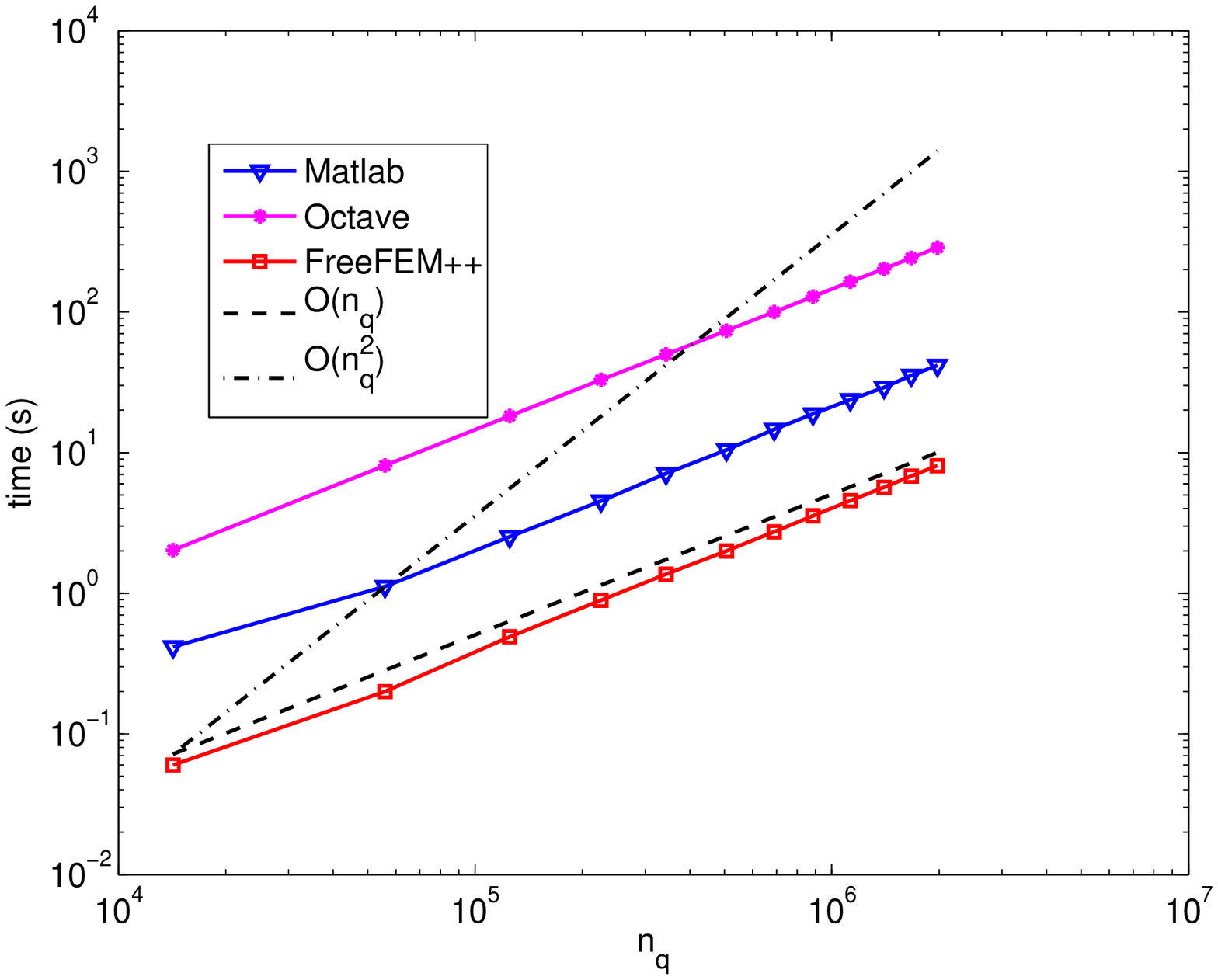}{\zoom}
{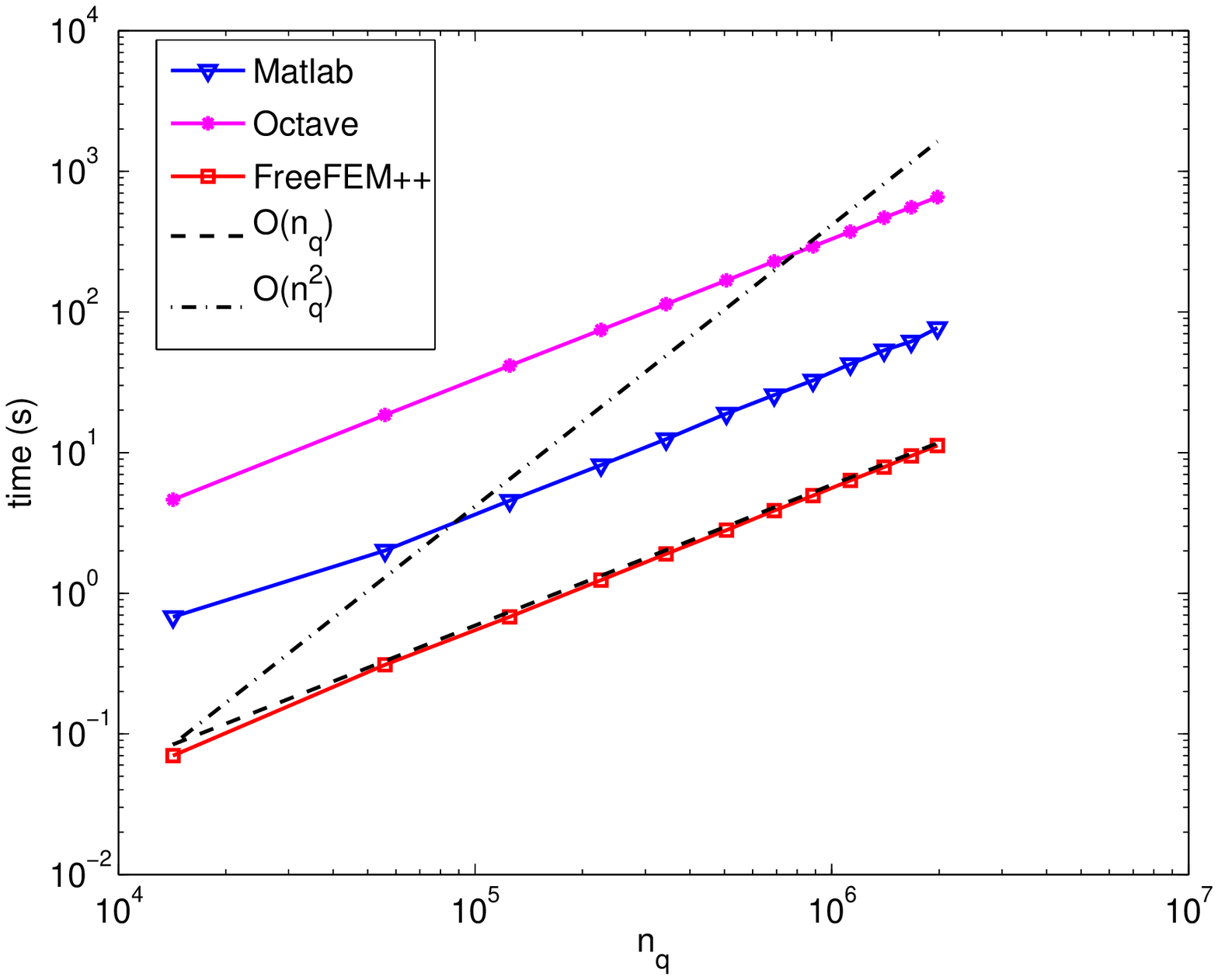}{\zoom}
\imageps{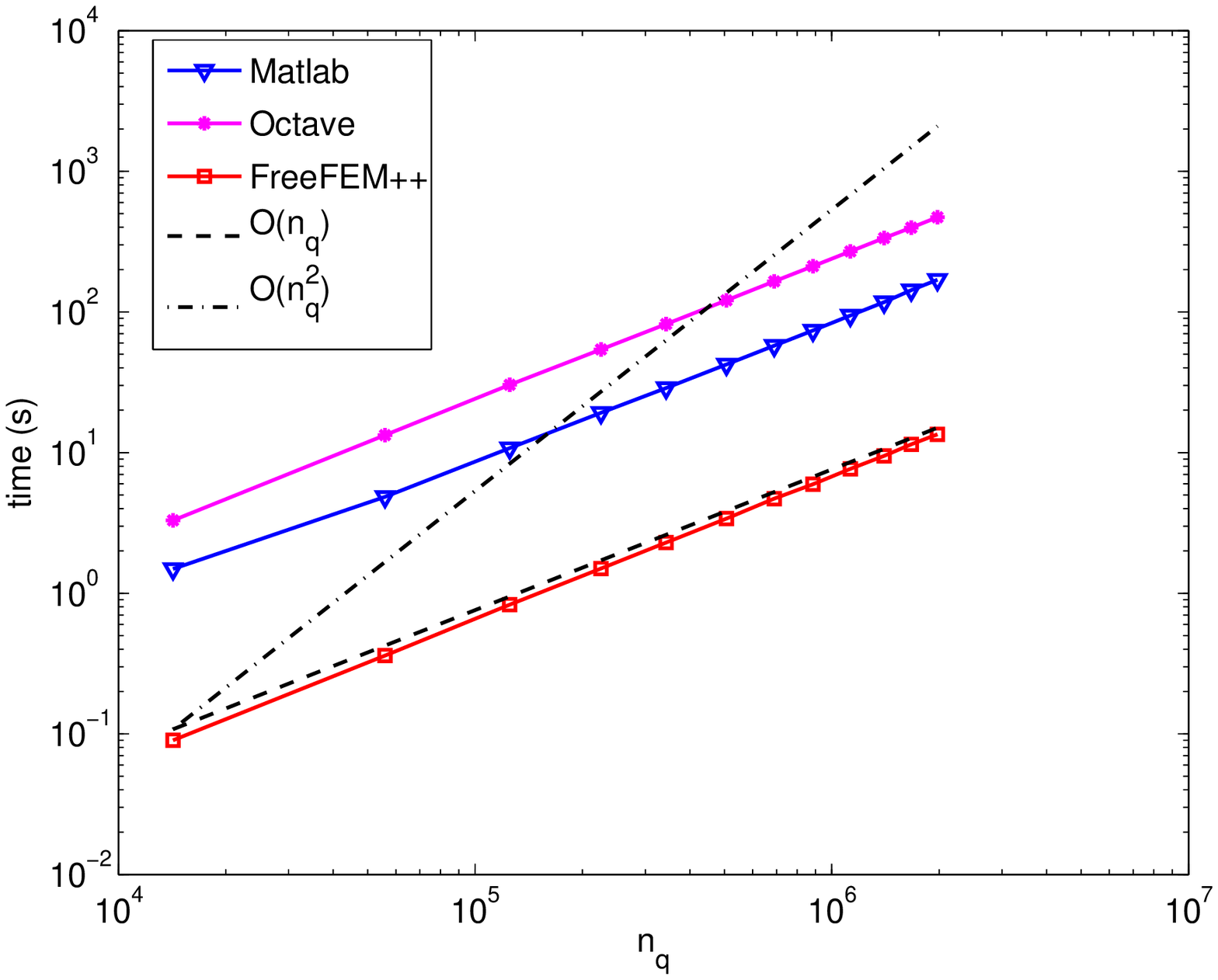}{\zoom}
\caption{Comparison of the matrix assembly codes : \texttt{OptV1} in Matlab/Octave and FreeFEM++,
  for the mass (top left), weighted mass (top right) and stiffness (bottom) matrices.
  \label{OptV1VsFreeFEM}}
\end{figure}

From these arrays, it is then possible to build the three global arrays
$\vecb{I}_g,$ $\vecb{J}_g$ and $\vecb{K}_g,$ of size $9\nme \times 1$
defined by : $\forall k \in \ENS{1}{\nme},$ $\forall il \in\ENS{1}{9}$,
\begin{eqnarray*}
\vecb{K}_g(9(k-1)+il)&=&\vecb{K}^e_k(il),\\
\vecb{I}_g(9(k-1)+il)&=&\vecb{I}^e_k(il),\\
\vecb{J}_g(9(k-1)+il)&=&\vecb{J}^e_k(il).
\end{eqnarray*}
On Figure~\ref{AssMatV1}, we show the insertion of the local array $\vecb{K}^e_k$
into the global 1d-array $\vecb{K}_g,$ and, for representation convenience, we draw them in line format.
We make the same operation for the two other arrays.
\noindent
We give in Listing~\ref{Assemblage_OptV1} the Matlab/Octave associated code where 
the global vectors $\vecb{I}_g,$ $\vecb{J}_g$ and $\vecb{K}_g$ are stored column-wise.
The complete listings and the values of the computation times are given in Appendices \ref{App:codeV1} and
\ref{App:OptV1vsFreeFem} respectively.
On Figure~\ref{OptV1VsFreeFEM}, we show the 
computation times of the Matlab, Octave and FreeFEM++ codes
versus the number of vertices of the mesh (unit disk).
The complexity of the Matlab/Octave codes seems now linear (i.e. $\GrandO(\nq)$) as for FreeFEM++.
However, FreeFEM++ is still much more faster than Matlab/Octave (about a factor 5 for the mass matrix,
6.5 for the weighted mass matrix and 12.5 for the stiffness matrix, in Matlab).
\begin{center}
\begin{minipage}[t]{\ListingWidth}
\MonMatlabNonumber
  \begin{lstlisting}[caption=Optimized matrix assembly code - version 1,label=Assemblage_OptV1]
Ig=zeros(9*nme,1);Jg=zeros(9*nme,1);Kg=zeros(9*nme,1);
ii=[1 2 3 1 2 3 1 2 3]; jj=[1 1 1 2 2 2 3 3 3];
kk=1:9;
for k=1:nme
  E=ElemMat(areas(k), ...);
  Ig(kk)=me(ii,k); 
  Jg(kk)=me(jj,k);
  Kg(kk)=E(:);
  kk=kk+9;
end
M=sparse(Ig,Jg,Kg,nq,nq);
\end{lstlisting}
\end{minipage}
\end{center}
To further improve the efficiency of the codes, we introduce now a second optimized version
of the assembly algorithm.

\section{Optimized matrix assembly - version 2 (OptV2)}
\label{sssec:CodeMassOptV2}
We present the optimized version 2 of the algorithm where no loop is used.

We define three 2d-arrays that allow to store all the element matrices 
as well as their positions in the global matrix. We denote by $\MAT{K}_g,$ $\MAT{I}_g$
and $\MAT{J}_g$ these $9$-by-$\nme$ arrays, defined
$\forall k \in \ENS{1}{\nme},$ $\forall il \in\ENS{1}{9}$ by
\begin{eqnarray*}
\MAT{K}_g(il,k)=\vecb{K}^e_k(il),\quad\quad
\MAT{I}_g(il,k)=\vecb{I}^e_k(il),\quad\quad
\MAT{J}_g(il,k)=\vecb{J}^e_k(il).
\end{eqnarray*}
The three local arrays $\vecb{K}^e_k,$ $\vecb{I}^e_k$ and $\vecb{J}^e_k$ are thus stored
in the $k$-th column of the global arrays $\MAT{K}_g,$ $\MAT{I}_g$ and $\MAT{J}_g$ respectively.

A natural way to build these three arrays consists in using a loop through the triangles
$T_k$ in which we insert the local arrays column-wise, see Figure~\ref{AssMatV2}.
Once these arrays are determined, the matrix assembly is obtained with the Matlab/Octave command
\vspace{-1mm}
\begin{center}
  \lstinline{M = sparse(Ig(:),Jg(:),Kg(:),nq,nq);}
\end{center}
We remark that the matrices containing global indices $\MAT{I}_g$ and $\MAT{J}_g$ may be computed, in Matlab/Octave,
without any loop. 
For the computation of these two matrices, on the left we give the usual code and on the right
the vectorized code :
\begin{center}
\begin{minipage}[t]{0.54\textwidth}
\MonMatlabNonumber
\begin{lstlisting}
Ig=zeros(9,nme);Jg=zeros(9,nme);
for k=1:nme
  Ig(:,k)=me([1 2 3 1 2 3 1 2 3],k);
  Jg(:,k)=me([1 1 1 2 2 2 3 3 3],k);
end
\end{lstlisting}
\end{minipage}\hspace{0.01\textwidth}
\begin{minipage}[t]{0.43\textwidth}
\MonMatlabNonumber
\begin{lstlisting}
Ig=me([1 2 3 1 2 3 1 2 3],:);
Jg=me([1 1 1 2 2 2 3 3 3],:);
\end{lstlisting}
\end{minipage}
\end{center}
Another way to present this computation, used and adapted in
Section~\ref{sec:Elasticity}, is given by
\begin{remark}\label{remq:IgJg}
Denoting $\mathcal{I}_k=[\me(1,k), \, \me(2,k), \, \me(3,k)]$ and
\begin{equation*}\label{def:T1}
\MAT{T}=\begin{pmatrix}
\mathcal{I}_1(1) & \hdots & \mathcal{I}_k(1) & \hdots & \mathcal{I}_\nme(1)\\ 
\mathcal{I}_1(2) & \hdots & \mathcal{I}_k(2) & \hdots &\mathcal{I}_\nme(2)\\ 
\mathcal{I}_1(3) & \hdots & \mathcal{I}_k(3) & \hdots &\mathcal{I}_\nme(3)
\end{pmatrix},
\end{equation*}
then, in that case $\MAT{T}=\me$, and $\MAT{I}_g$ and $\MAT{J}_g$ may
be computed from $\MAT{T}$ as follows:
\begin{center}
\begin{minipage}[t]{0.9\textwidth}
\MonMatlabNonumber
\begin{lstlisting}
ii=[1 1 1; 2 2 2; 3 3 3]; jj=ii';
Ig=T(ii(:),:); Jg=T(jj(:),:);
\end{lstlisting}
\end{minipage}
\end{center}
\end{remark}
\begin{figure}[H]
\footnotesize{
  \colorlet{colorquiver}{DarkBlue}
\colorlet{matelem}{Maroon}
\colorlet{indelem}{Indigo}
\definecolor{ffqqqq}{rgb}{1,0,0}
\definecolor{cqcqcq}{rgb}{0.75,0.75,0.75}
\definecolor{color1}{rgb}{0.75,0.75,0.}
\newcommand{\myunit}{0.25 cm}
\tikzset{
   node style sp/.style={draw,rectangle,rounded corners=3pt},
    node style ge/.style={circle,minimum size=\myunit},
    arrow style mul/.style={draw,sloped,midway,fill=white},
    arrow style plus/.style={midway,sloped,fill=white},
}


\newcommand{\DrawMatrixMee}[3]{%
\matrix (M) [matrix of math nodes,%
             nodes = {node style ge},%
             left delimiter  = (,%
             right delimiter = )] at (#1,#2)
{%
  \node[node style sp,color=blue,text=black] {#3_{1,1}};%
  & \node {};&\node[node style sp,color=magenta,text=black] {#3_{1,2}};%
  & \node {};&\node[node style sp,color=green,text=black] {#3_{1,3}};  \\ \node {};& & \\
    \node[node style sp,color=blue,text=black] {#3_{2,1}};%
  &  \node {};&\node[node style sp,color=magenta,text=black] {#3_{2,2}};%
  &  \node {};&\node[node style sp,color=green,text=black] {#3_{2,3}}; \\ \node {};& & \\
   \node[node style sp,color=blue,text=black] {#3_{3,1}};%
  & \node {};& \node[node style sp,color=magenta,text=black] {#3_{3,2} };%
  &  \node {};&\node[node style sp,color=green,text=black] {#3_{3,3}}; \\ 
};
\node [draw,left=20pt] at (M.west)     { $\mathbb{M}^e(T_k)$ :};
}

\newcommand{\DrawMatrixMe}[3]{%
\matrix (M) [matrix of math nodes,%
             nodes = {node style ge},%
             left delimiter  = (,%
             right delimiter = )] at (0,0)
{%
  \node[node style sp,color=blue,text=black] {m_{1,1}};%
  & \node {};&\node[node style sp,color=magenta,text=black] {m_{1,2}};%
  & \node {};&\node[node style sp,color=green,text=black] {m_{1,3}};  \\ \node {};& & \\%
    \node[node style sp,color=blue,text=black] {m_{2,1}};%
  &  \node {};&\node[node style sp,color=magenta,text=black] {m_{2,2}};%
  &  \node {};&\node[node style sp,color=green,text=black] {m_{2,3}}; \\ \node {};& & \\%
   \node[node style sp,color=blue,text=black] {m_{3,1}};%
  & \node {};& \node[node style sp,color=magenta,text=black] {m_{3,2} };%
  &  \node {};&\node[node style sp,color=green,text=black] {m_{3,3}}; \\%
};
}

\newcommand{\DrawTab}[4]{%
  \draw (#1,#2)-- ({#1},{#2 +#4})-- ({#1 + #3},{#2 +#4}) -- ({#1+#3}, {#2}) -- (#1,#2);

%
 \foreach \k in {1,...,#4}
   {\draw [dotted] (#1,{#2+\k}) -- ({#1+#3},{#2+\k});}
}

\newcommand{\DrawVertTabM}[4]{
  \DrawTab{#1}{#2}{1}{9}
  \foreach \j in {1,...,3}{
    \foreach \i in {1,...,3}{
      \draw [color=matelem](#1+0.5,#2+9.5-\i-3*\j+3) node[anchor=center,scale=0.6] {$ #3_{\i,\j}$};
    }
  }
  \draw (#1+0.5,#2+10) node[anchor=center] {$ #4$ };
}
\newcommand{\DrawVertTabI}[4]{
  \DrawTab{#1}{#2}{1}{9}
  \foreach \j in {1,...,3}{
    \foreach \i in {1,...,3}{
      \draw [color=indelem](#1+0.5,#2+9.5-\i-3*\j+3) node[anchor=center,scale=0.6] {$ #3_{\i}$};
    }
  }
  \draw (#1+0.5,#2+10) node[anchor=center] {$ #4$ };
}

\newcommand{\DrawVertTabJ}[4]{
  \DrawTab{#1}{#2}{1}{9}
  \foreach \j in {1,...,3}{
    \foreach \i in {1,...,3}{
      \draw [color=indelem](#1+0.5,#2+9.5-\i-3*\j+3) node[anchor=center,scale=0.6] {$ #3_{\j}$};
    }
  }
  \draw (#1+0.5,#2+10) node[anchor=center] {$ #4$ };
}

\newcommand{\DrawKg}[2]{
\draw (#1+5-5,#2-4.5+4.5) node[anchor=center] {$ \mathbb{K}_g$ };
\DrawTab{#1+1.5-5}{#2-3.5+4.5}{1}{9}
\DrawTab{#1+2.5-5}{#2-3.5+4.5}{1}{9}
\DrawVertTabM{#1+4.5-5}{#2-3.5+4.5}{e^k}{}
\DrawTab{#1+6.5-5}{#2-3.5+4.5}{1}{9}
\DrawTab{#1+7.5-5}{#2-3.5+4.5}{1}{9}
\draw [dashed] (#1+3.5-5,#2-3.5+4.5) -- (#1+4.5-5,#2-3.5+4.5);
\draw [dashed] (#1+5.5-5,#2-3.5+4.5) -- (#1+6.5-5,#2-3.5+4.5);
\draw [dashed] (#1+3.5-5,#2+5.5+4.5) -- (#1+4.5-5,#2+5.5+4.5);
\draw [dashed] (#1+5.5-5,#2+5.5+4.5) -- (#1+6.5-5,#2+5.5+4.5);

 \foreach \i in {1,...,9}
   \draw [color=red](#1+1-5,#2+6-\i+4.5) node[anchor=center,scale=0.6] {$ \i $};
 \draw [color=red](#1+2-5,#2+6+4.5) node[anchor=center,scale=0.6] {$ 1 $};
 \draw [color=red](#1+3-5,#2+6+4.5) node[anchor=center,scale=0.6] {$ 2 $};
 \draw [color=red](#1+4-5,#2+6+4.5) node[anchor=center,scale=0.6] {$ \hdots $};
 \draw [color=red](#1+5-5,#2+6+4.5) node[anchor=center,scale=0.6] {$ k $};
 \draw [color=red](#1+6-5,#2+6+4.5) node[anchor=center,scale=0.6] {$ \hdots $};
 \draw [color=red](#1+8-5,#2+6+4.5) node[anchor=center,scale=0.6] {$ \nme $};
}

\newcommand{\DrawIg}[2]{
\draw (#1+5-5,#2-4.5+4.5) node[anchor=center] {$ \mathbb{I}_g$ };
\DrawTab{#1+1.5-5}{#2-3.5+4.5}{1}{9}
\DrawTab{#1+2.5-5}{#2-3.5+4.5}{1}{9}
\DrawVertTabI{#1+4.5-5}{#2-3.5+4.5}{i^k}{}
\DrawTab{#1+6.5-5}{#2-3.5+4.5}{1}{9}
\DrawTab{#1+7.5-5}{#2-3.5+4.5}{1}{9}
\draw [dashed] (#1+3.5-5,#2-3.5+4.5) -- (#1+4.5-5,#2-3.5+4.5);
\draw [dashed] (#1+5.5-5,#2-3.5+4.5) -- (#1+6.5-5,#2-3.5+4.5);
\draw [dashed] (#1+3.5-5,#2+5.5+4.5) -- (#1+4.5-5,#2+5.5+4.5);
\draw [dashed] (#1+5.5-5,#2+5.5+4.5) -- (#1+6.5-5,#2+5.5+4.5);

 \foreach \i in {1,...,9}
   \draw [color=red](#1+1-5,#2+6-\i+4.5) node[anchor=center,scale=0.6] {$ \i $};
 \draw [color=red](#1+2-5,#2+6+4.5) node[anchor=center,scale=0.6] {$ 1 $};
 \draw [color=red](#1+3-5,#2+6+4.5) node[anchor=center,scale=0.6] {$ 2 $};
 \draw [color=red](#1+4-5,#2+6+4.5) node[anchor=center,scale=0.6] {$ \hdots $};
 \draw [color=red](#1+5-5,#2+6+4.5) node[anchor=center,scale=0.6] {$ k $};
 \draw [color=red](#1+6-5,#2+6+4.5) node[anchor=center,scale=0.6] {$ \hdots $};
 \draw [color=red](#1+8-5,#2+6+4.5) node[anchor=center,scale=0.6] {$ \nme $};
}

\newcommand{\DrawJg}[2]{
\draw (#1+5-5,#2-4.5+4.5) node[anchor=center] {$ \mathbb{J}_g$ };
\DrawTab{#1+1.5-5}{#2-3.5+4.5}{1}{9}
\DrawTab{#1+2.5-5}{#2-3.5+4.5}{1}{9}
\DrawVertTabJ{#1+4.5-5}{#2-3.5+4.5}{i^k}{}
\DrawTab{#1+6.5-5}{#2-3.5+4.5}{1}{9}
\DrawTab{#1+7.5-5}{#2-3.5+4.5}{1}{9}
\draw [dashed] (#1+3.5-5,#2-3.5+4.5) -- (#1+4.5-5,#2-3.5+4.5);
\draw [dashed] (#1+5.5-5,#2-3.5+4.5) -- (#1+6.5-5,#2-3.5+4.5);
\draw [dashed] (#1+3.5-5,#2+5.5+4.5) -- (#1+4.5-5,#2+5.5+4.5);
\draw [dashed] (#1+5.5-5,#2+5.5+4.5) -- (#1+6.5-5,#2+5.5+4.5);

 \foreach \i in {1,...,9}
   \draw [color=red](#1+1-5,#2+6-\i+4.5) node[anchor=center,scale=0.6] {$ \i $};
 \draw [color=red](#1+2-5,#2+6+4.5) node[anchor=center,scale=0.6] {$ 1 $};
 \draw [color=red](#1+3-5,#2+6+4.5) node[anchor=center,scale=0.6] {$ 2 $};
 \draw [color=red](#1+4-5,#2+6+4.5) node[anchor=center,scale=0.6] {$ \hdots $};
 \draw [color=red](#1+5-5,#2+6+4.5) node[anchor=center,scale=0.6] {$ k $};
 \draw [color=red](#1+6-5,#2+6+4.5) node[anchor=center,scale=0.6] {$ \hdots $};
 \draw [color=red](#1+8-5,#2+6+4.5) node[anchor=center,scale=0.6] {$ \nme $};
}

\begin{center}
\begin{tikzpicture}[>=latex,scale=0.45]
\tikzstyle{operation}=[->,>=latex]
\tikzstyle{etiquette}=[midway,fill=black!20]
\tikzset{LabelStyle/.style =   {draw,
                                  fill           = yellow,
                                  text           = red}}
\matrix (M) [matrix of math nodes,%
             nodes = {node style ge},%
             left delimiter  = (,%
             right delimiter = )] at (-7,8)
{%
  \node[text=matelem] {e_{1,1}^k};%
  & \node {};&\node[text=matelem] {e_{1,2}^k};%
  & \node {};&\node[text=matelem] {e_{1,3}^k};  \\
    \node[text=matelem] {e_{2,1}^k};%
  &  \node {};&\node[text=matelem] {e_{2,2}^k};%
  &  \node {};&\node[text=matelem] {e_{2,3}^k};  \\
   \node[text=matelem] {e_{3,1}^k};%
  & \node {};& \node[text=matelem] {e_{3,2}^k };%
  &  \node {};&\node[text=matelem] {e_{3,3}^k}; \\ 
};
\node [above=10pt] at (M.north) { $\mathbb{E}(T_k)$};

\draw [->,color=colorquiver] (-1,8) -- (2.5,8);

\DrawVertTabM{4}{3}{e^k}{\vecb{K}^e_k}
\DrawVertTabI{7}{3}{i^k}{\vecb{I}^e_k}
\DrawVertTabJ{10}{3}{i^k}{\vecb{J}^e_k}

\DrawKg{-9}{-13}
\DrawIg{0.5}{-13}
\DrawJg{10}{-13}

\draw [->,color=colorquiver] (4.5,3) to[in=90,out=270] (-9,-2);
\draw [->,color=colorquiver] (7.5,3) to[in=90,out=270] (0.5,-2);
\draw [->,color=colorquiver] (10.5,3) to[in=90,out=270] (10,-2);

\end{tikzpicture}
\end{center}

}
  \caption{Insertion of an element matrix\label{AssMatV2} in the global array - Version 2}
\end{figure}
It remains to vectorize the computation of the 2d-array $\MAT{K}_g$.
The usual code, corresponding to a column-wise computation, is :
\vspace{-2mm}
\begin{center}
\begin{minipage}[t]{\ListingWidth}
\MonMatlabNonumber
  \begin{lstlisting}[caption=Usual assembly (column-wise computation),label=UsualAssemblycwc]
Kg=zeros(9,nme);
for k=1:nme
  E=ElemMat(areas(k), ...);
  Kg(:,k)=E(:);
end
\end{lstlisting}
\end{minipage}
\end{center}
The vectorization of this code is done by the computation of the array  $\MAT{K}_g$ row-wise,
for each matrix assembly. This corresponds to the permutation of the loop through the elements
with the local loops, in the classical matrix assembly code (see Listing~\ref{Assemblage:basique}).
This vectorization differs from the one proposed in~\cite{Hannukainen:IFE:2012} as it doesn't use
any quadrature formula and from the one in \cite{Chen:PFE:2011} by the full vectorization
of the arrays $\MAT{I}_g$ and $\MAT{J}_g.$

We describe below this method for each matrix defined in Section~\ref{sssec:PositionduProbleme}.

\subsection{Mass matrix assembly}
\begin{figure}[H]
  \definecolor{ffqqqq}{rgb}{1,0,0}
\definecolor{cqcqcq}{rgb}{0.75,0.75,0.75}
\definecolor{color1}{rgb}{0.75,0.75,0.}
\definecolor{IndiceColor}{rgb}{0,0,0}
\newcommand{\myunit}{0.25 cm}
\tikzset{
   node style sp/.style={draw,rectangle,rounded corners=3pt},
    node style ge/.style={circle,minimum size=\myunit},
    arrow style mul/.style={draw,sloped,midway,fill=white},
    arrow style plus/.style={midway,sloped,fill=white},
}


\newcommand{\DrawMatrixMee}[3]{%
\matrix (M) [matrix of math nodes,%
             nodes = {node style ge},%
             left delimiter  = (,%
             right delimiter = )] at (#1,#2)
{%
  \node[node style sp,color=blue,text=black] {#3_{1,1}};%
  & \node {};&\node[node style sp,color=magenta,text=black] {#3_{1,2}};%
  & \node {};&\node[node style sp,color=green,text=black] {#3_{1,3}};  \\ \node {};& & \\
    \node[node style sp,color=blue,text=black] {#3_{2,1}};%
  &  \node {};&\node[node style sp,color=magenta,text=black] {#3_{2,2}};%
  &  \node {};&\node[node style sp,color=green,text=black] {#3_{2,3}}; \\ \node {};& & \\
   \node[node style sp,color=blue,text=black] {#3_{3,1}};%
  & \node {};& \node[node style sp,color=magenta,text=black] {#3_{3,2} };%
  &  \node {};&\node[node style sp,color=green,text=black] {#3_{3,3}}; \\ 
};
\node [draw,left=20pt] at (M.west)     { $\mathbb{M}^e(T_k)$ :};
}

\newcommand{\DrawMatrixMe}[3]{%
\matrix (M) [matrix of math nodes,%
             nodes = {node style ge},%
             left delimiter  = (,%
             right delimiter = )] at (0,0)
{%
  \node[node style sp,color=blue,text=black] {m_{1,1}};%
  & \node {};&\node[node style sp,color=magenta,text=black] {m_{1,2}};%
  & \node {};&\node[node style sp,color=green,text=black] {m_{1,3}};  \\ \node {};& & \\%
    \node[node style sp,color=blue,text=black] {m_{2,1}};%
  &  \node {};&\node[node style sp,color=magenta,text=black] {m_{2,2}};%
  &  \node {};&\node[node style sp,color=green,text=black] {m_{2,3}}; \\ \node {};& & \\%
   \node[node style sp,color=blue,text=black] {m_{3,1}};%
  & \node {};& \node[node style sp,color=magenta,text=black] {m_{3,2} };%
  &  \node {};&\node[node style sp,color=green,text=black] {m_{3,3}}; \\%
};
}

\newcommand{\DrawTab}[4]{%
  \draw (#1,#2)-- ({#1},{#2 +#4})-- ({#1 + #3},{#2 +#4}) -- ({#1+#3}, {#2}) -- (#1,#2);

%
 \foreach \k in {1,...,#4}
   {\draw [dotted] (#1,{#2+\k}) -- ({#1+#3},{#2+\k});}
}

\newcommand{\DrawVertTabM}[4]{
  \DrawTab{#1}{#2}{1}{9}
  \foreach \j in {1,...,3}{
    \foreach \i in {1,...,3}{
      \draw [color=magenta](#1+0.5,#2+9.5-\i-3*\j+3) node[anchor=center,scale=0.6] {$ #3_{\i,\j}$};
    }
  }
  \draw (#1+0.5,#2+10) node[anchor=center] {$ #4$ };
}
\newcommand{\DrawVertTabI}[4]{
  \DrawTab{#1}{#2}{1}{9}
  \foreach \j in {1,...,3}{
    \foreach \i in {1,...,3}{
      \draw [color=magenta](#1+0.5,#2+9.5-\i-3*\j+3) node[anchor=center,scale=0.6] {$ #3_{\i}$};
    }
  }
  \draw (#1+0.5,#2+10) node[anchor=center] {$ #4$ };
}

\newcommand{\DrawVertTabJ}[4]{
  \DrawTab{#1}{#2}{1}{9}
  \foreach \j in {1,...,3}{
    \foreach \i in {1,...,3}{
      \draw [color=magenta](#1+0.5,#2+9.5-\i-3*\j+3) node[anchor=center,scale=0.6] {$ #3_{\j}$};
    }
  }
  \draw (#1+0.5,#2+10) node[anchor=center] {$ #4$ };
}

\newcommand{\DrawKg}[2]{
\draw (#1+5-5,#2-4.5+4.5) node[anchor=center] {$ \mathbb{K}_g$ };
\DrawTab{#1+1.5-5}{#2-3.5+4.5}{1}{9}
\DrawTab{#1+2.5-5}{#2-3.5+4.5}{1}{9}
\DrawVertTabM{#1+4.5-5}{#2-3.5+4.5}{m^k}{}
\DrawTab{#1+6.5-5}{#2-3.5+4.5}{1}{9}
\DrawTab{#1+7.5-5}{#2-3.5+4.5}{1}{9}
\draw [dashed] (#1+3.5-5,#2-3.5+4.5) -- (#1+4.5-5,#2-3.5+4.5);
\draw [dashed] (#1+5.5-5,#2-3.5+4.5) -- (#1+6.5-5,#2-3.5+4.5);
\draw [dashed] (#1+3.5-5,#2+5.5+4.5) -- (#1+4.5-5,#2+5.5+4.5);
\draw [dashed] (#1+5.5-5,#2+5.5+4.5) -- (#1+6.5-5,#2+5.5+4.5);

 \foreach \i in {1,...,9}
   \draw [color=\IndiceColor](#1+1-5,#2+6-\i+4.5) node[anchor=center,scale=0.6] {$ \i $};
 \draw [color=\IndiceColor](#1+2-5,#2+6+4.5) node[anchor=center,scale=0.6] {$ 1 $};
 \draw [color=\IndiceColor](#1+3-5,#2+6+4.5) node[anchor=center,scale=0.6] {$ 2 $};
 \draw [color=\IndiceColor](#1+4-5,#2+6+4.5) node[anchor=center,scale=0.6] {$ \hdots $};
 \draw [color=\IndiceColor](#1+5-5,#2+6+4.5) node[anchor=center,scale=0.6] {$ k $};
 \draw [color=\IndiceColor](#1+6-5,#2+6+4.5) node[anchor=center,scale=0.6] {$ \hdots $};
 \draw [color=\IndiceColor](#1+8-5,#2+6+4.5) node[anchor=center,scale=0.6] {$ \nme $};
}

\newcommand{\DrawIg}[2]{
\draw (#1+5-5,#2-4.5+4.5) node[anchor=center] {$ \mathbb{I}_g$ };
\DrawTab{#1+1.5-5}{#2-3.5+4.5}{1}{9}
\DrawTab{#1+2.5-5}{#2-3.5+4.5}{1}{9}
\DrawVertTabI{#1+4.5-5}{#2-3.5+4.5}{i^k}{}
\DrawTab{#1+6.5-5}{#2-3.5+4.5}{1}{9}
\DrawTab{#1+7.5-5}{#2-3.5+4.5}{1}{9}
\draw [dashed] (#1+3.5-5,#2-3.5+4.5) -- (#1+4.5-5,#2-3.5+4.5);
\draw [dashed] (#1+5.5-5,#2-3.5+4.5) -- (#1+6.5-5,#2-3.5+4.5);
\draw [dashed] (#1+3.5-5,#2+5.5+4.5) -- (#1+4.5-5,#2+5.5+4.5);
\draw [dashed] (#1+5.5-5,#2+5.5+4.5) -- (#1+6.5-5,#2+5.5+4.5);

 \foreach \i in {1,...,9}
   \draw [color=red](#1+1-5,#2+6-\i+4.5) node[anchor=center,scale=0.6] {$ \i $};
 \draw [color=red](#1+2-5,#2+6+4.5) node[anchor=center,scale=0.6] {$ 1 $};
 \draw [color=red](#1+3-5,#2+6+4.5) node[anchor=center,scale=0.6] {$ 2 $};
 \draw [color=red](#1+4-5,#2+6+4.5) node[anchor=center,scale=0.6] {$ \hdots $};
 \draw [color=red](#1+5-5,#2+6+4.5) node[anchor=center,scale=0.6] {$ k $};
 \draw [color=red](#1+6-5,#2+6+4.5) node[anchor=center,scale=0.6] {$ \hdots $};
 \draw [color=red](#1+8-5,#2+6+4.5) node[anchor=center,scale=0.6] {$ \nme $};
}

\newcommand{\DrawJg}[2]{
\draw (#1+5-5,#2-4.5+4.5) node[anchor=center] {$ \mathbb{J}_g$ };
\DrawTab{#1+1.5-5}{#2-3.5+4.5}{1}{9}
\DrawTab{#1+2.5-5}{#2-3.5+4.5}{1}{9}
\DrawVertTabJ{#1+4.5-5}{#2-3.5+4.5}{i^k}{}
\DrawTab{#1+6.5-5}{#2-3.5+4.5}{1}{9}
\DrawTab{#1+7.5-5}{#2-3.5+4.5}{1}{9}
\draw [dashed] (#1+3.5-5,#2-3.5+4.5) -- (#1+4.5-5,#2-3.5+4.5);
\draw [dashed] (#1+5.5-5,#2-3.5+4.5) -- (#1+6.5-5,#2-3.5+4.5);
\draw [dashed] (#1+3.5-5,#2+5.5+4.5) -- (#1+4.5-5,#2+5.5+4.5);
\draw [dashed] (#1+5.5-5,#2+5.5+4.5) -- (#1+6.5-5,#2+5.5+4.5);

 \foreach \i in {1,...,9}
   \draw [color=red](#1+1-5,#2+6-\i+4.5) node[anchor=center,scale=0.6] {$ \i $};
 \draw [color=red](#1+2-5,#2+6+4.5) node[anchor=center,scale=0.6] {$ 1 $};
 \draw [color=red](#1+3-5,#2+6+4.5) node[anchor=center,scale=0.6] {$ 2 $};
 \draw [color=red](#1+4-5,#2+6+4.5) node[anchor=center,scale=0.6] {$ \hdots $};
 \draw [color=red](#1+5-5,#2+6+4.5) node[anchor=center,scale=0.6] {$ k $};
 \draw [color=red](#1+6-5,#2+6+4.5) node[anchor=center,scale=0.6] {$ \hdots $};
 \draw [color=red](#1+8-5,#2+6+4.5) node[anchor=center,scale=0.6] {$ \nme $};
}

\newcommand{\DrawHorTab}[4]{%
  \draw [#4] (#1,#2)-- ({#1+#3},{#2})-- ({#1 + #3},{#2 +1}) -- ({#1}, {#2+1}) -- (#1,#2);
  \draw [dotted] ({#1+1},{#2}) -- ({#1+1},{#2+1});
  \draw [dotted] ({#1+2},{#2}) -- ({#1+2},{#2+1});
  \draw [loosely dotted] (#1+2.5,#2+0.5) -- ({#1+#3-2.5},{#2+0.5});
  \draw [dotted] ({#1+#3-2},{#2}) -- ({#1+#3-2},{#2+1});
  \draw [dotted] ({#1+#3-1},{#2}) -- ({#1+#3-1},{#2+1});
}
\newcommand{\DrawHorInd}[3]{%
\draw [color=IndiceColor](#1+0.5,#2+1) node[anchor=center,scale=0.6] {$ 1 $};
 \draw [color=IndiceColor](#1+1.5,#2+1) node[anchor=center,scale=0.6] {$ 2 $};
 \draw [color=IndiceColor](#1+2.5,#2+1) node[anchor=center,scale=0.6] {$ \hdots $};
 \draw [color=IndiceColor](#1+#3-1.5,#2+1) node[anchor=center,scale=0.6] {$ \hdots $};
 \draw [color=IndiceColor](#1+#3-0.5,#2+1) node[anchor=center,scale=0.6] {$ \nme $};
}

\newcommand{\DrawVerInd}[2]{
 \foreach \i in {1,...,9}
   \draw [color=IndiceColor](#1,#2+9-\i+0.5) node[anchor=center,scale=0.6] {$ \i $};
}
\newcommand{\DrawKgBis}[2]{%
\foreach \i in {1,5,9}
  \DrawHorTab{#1}{#2+\i-1}{9}{fill=blue!10};
\foreach \i in {2,3,4,6,7,8}
  \DrawHorTab{#1}{#2+\i-1}{9}{fill=blue!30};
}

\begin{center}
\begin{tikzpicture}[>=latex,scale=0.5]
\tikzset{LabelStyle/.style =   {draw,
                                  fill           = yellow,
                                  text           = red}}
\tikzstyle{operation}=[->,>=latex]
\draw (-0.5,0.5) node[anchor=east] { \texttt{areas} };
\DrawHorTab{0}{0.}{9}{}
\DrawHorInd{0}{0.5}{9}

\draw (2.6,-3.5) node[anchor=east] { $A_6$ };
\DrawHorTab{-8}{-4}{9}{fill=blue!10}
\DrawHorInd{-8}{-3.5}{9}
\draw [->,color=blue] (4.5,0) to[in=90,out=270] node[LabelStyle,scale=0.6]{$/6$} (-3.5,-3.0);

\draw (6.1,-3.5) node[anchor=west] { $A_{12}$ };
\DrawHorTab{8}{-4}{9}{fill=blue!30}
\DrawHorInd{8}{-3.5}{9}
\draw [->,color=blue] (4.5,0) to[in=90,out=270] node[LabelStyle,scale=0.6]{$/12$} (12.5,-3.0);

\draw (4.5,-6.0) node[anchor=center] { $\MAT{K}_g$ };
\DrawKgBis{0}{-16}
\DrawHorInd{0}{-7.5}{9}
\DrawVerInd{-0.5}{-16}
\draw [->,color=blue] (-3.5,-4.0) to[in=180,out=270] (-0.7,-7.5);
\draw [->,color=blue] (-3.5,-4.0) to[in=180,out=270] (-0.7,-11.5);
\draw [->,color=blue] (-3.5,-4.0) to[in=180,out=270] (-0.7,-15.5);
\DrawVerInd{9.5}{-16}
\draw [->,color=blue] (12.5,-4.0) to[in=0,out=270] (9.7,-8.5);
\draw [->,color=blue] (12.5,-4.0) to[in=0,out=270] (9.7,-9.5);
\draw [->,color=blue] (12.5,-4.0) to[in=0,out=270] (9.7,-10.5);
\draw [->,color=blue] (12.5,-4.0) to[in=0,out=270] (9.7,-12.5);
\draw [->,color=blue] (12.5,-4.0) to[in=0,out=270] (9.7,-13.5);
\draw [->,color=blue] (12.5,-4.0) to[in=0,out=270] (9.7,-14.5);

\end{tikzpicture}
\end{center}

\caption{Mass matrix assembly\label{AssMatV2c} - Version 2}
\end{figure}
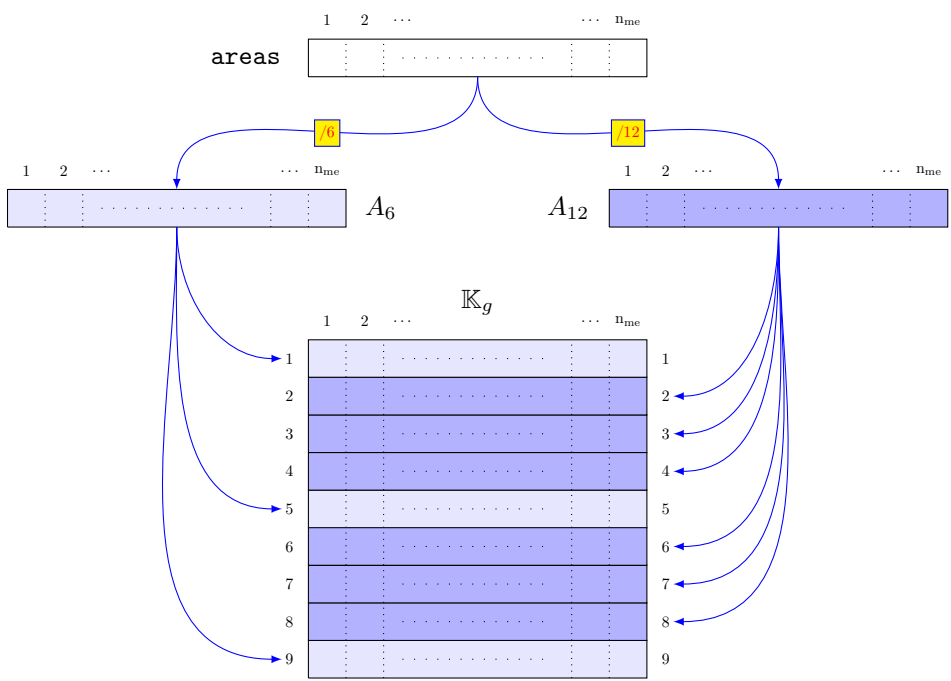

The element mass matrix $\MAT{M}^e(T_k)$ associated to the triangle
$T_k$ is given by (\ref{MassElem}).  The array $\MAT{K}_g$ is defined
by : $\forall k\in\ENS{1}{\nme},$
\begin{eqnarray*}
\MAT{K}_g(\il,k)&=&\frac{|T_k|}{6},\ \forall \il\in\{1,5,9\}, \\
\MAT{K}_g(\il,k)&=&\frac{|T_k|}{12},\  \forall \il\in\{2,3,4,6,7,8\}. 
\end{eqnarray*}
Then we build two arrays $A_6$ and $A_{12}$ of size $1\times \nme$ such that $\forall k\in\ENS{1}{\nme}$ :
\begin{eqnarray*}
A_6(k)&=&\frac{|T_k|}{6}, \quad\quad A_{12}(k)=\frac{|T_k|}{12}.
\end{eqnarray*}
The rows $\{1,5,9\}$ in the array $\MAT{K}_g$ correspond to $A_{6}$ and the rows $\{2,3,4,6,7,8\}$ to $A_{12}$,
see Figure~\ref{AssMatV2c}.
The Matlab/Octave code associated to this technique is :
\begin{center}
\begin{minipage}[t]{\ListingWidth}
\MonMatlab
\begin{lstlisting}[caption=Optimized matrix assembly code - version 2 (Mass matrix)]
function [M]=MassAssemblingP1OptV2(nq,nme,me,areas)
Ig = me([1 2 3 1 2 3 1 2 3],:);
Jg = me([1 1 1 2 2 2 3 3 3],:);
A6=areas/6;
A12=areas/12;
Kg = [A6;A12;A12;A12;A6;A12;A12;A12;A6];
M = sparse(Ig(:),Jg(:),Kg(:),nq,nq);
\end{lstlisting}
\end{minipage}
\end{center}

\subsection{Weighted mass matrix assembly}
The element weighted mass matrices $\MassWElem{\Local{w}}(T_k)$
are given by (\ref{MassWElem}).  We introduce the array $\vecb{T_w}$
of length $\nq$ defined by $\vecb{T_w}(i)=w(\q^i),$ for all $i\in
\ENS{1}{\nq}$ and the three arrays $\vecb{W}_\il, \ 1\le \il\le 3$, of
length $\nme$, defined for all $k\in\ENS{1}{\nme}$ by
$\vecb{W}_\il(k)=\frac{|T_k|}{30}\vecb{T_w}(\me(\il,k))$.

The code for computing these three arrays is given below, 
in a non-vectorized form (on the left) and in a vectorized form (on the right):
\begin{center}
\hspace{2mm}
\begin{minipage}[t]{0.48\textwidth}
\MonMatlabNonumber
\begin{lstlisting}
W1=zeros(1,nme);
W2=zeros(1,nme);
W3=zeros(1,nme);
for k=1:nme
  W1(k)=Tw(me(1,k))*areas(k)/30;
  W2(k)=Tw(me(2,k))*areas(k)/30;
  W3(k)=Tw(me(3,k))*areas(k)/30;
end
\end{lstlisting}
\end{minipage}\hspace{0.05\textwidth}
\begin{minipage}[t]{0.4\textwidth}
\MonMatlabNonumber
\begin{lstlisting}
W1=Tw(me(1,:)).*areas/30;
W2=Tw(me(2,:)).*areas/30;
W3=Tw(me(3,:)).*areas/30;
\end{lstlisting}
\end{minipage}
\end{center}
We follow the method described on Figure~\ref{AssMatV2}.
We have to vectorize the computation of $\MAT{K}_g$ (Listing~\ref{UsualAssemblycwc}).
Let $\vecb{K}_{1}$, $\vecb{K}_{2}$, $\vecb{K}_{3}$, $\vecb{K}_{5}$,
$\vecb{K}_{6}$, $\vecb{K}_{9}$ be six arrays of length $\nme$
defined by
$$
\begin{array}{lll}
\vecb{K}_{1}=3\vecb{W}_1 + \vecb{W}_2 +\vecb{W}_3 ,& 
\vecb{K}_{2}=\vecb{W}_1+\vecb{W}_2+\displaystyle\frac{\vecb{W}_3}{2},& 
\vecb{K}_{3}=\vecb{W}_1+\displaystyle\frac{\vecb{W}_2}{2}+\vecb{W}_3,\\ 
\vecb{K}_{5}=\vecb{W}_1 + 3\vecb{W}_2 +\vecb{W}_3,& 
\vecb{K}_{6}=\displaystyle\frac{\vecb{W}_1}{2}+\vecb{W}_2+\vecb{W}_3,& 
\vecb{K}_{9}=\vecb{W}_1 + \vecb{W}_2 +3\vecb{W}_3.
\end{array}
$$
The element weighted mass matrix and the $k$-th column of $\MAT{K}_g$
are respectively :
$$
\MassWElem{\Local{w}}(T_k)=\begin{pmatrix}
  \vecb{K}_1(k)  & \vecb{K}_2(k) & \vecb{K}_3(k)\\
  \vecb{K}_2(k)  & \vecb{K}_5(k) & \vecb{K}_6(k)\\
  \vecb{K}_3(k) & \vecb{K}_6(k) & \vecb{K}_9(k)
\end{pmatrix},\quad \quad
\MAT{K}_g(:,k)=
\begin{pmatrix}
\vecb{K}_1(k)  \\ \vecb{K}_2(k) \\ \vecb{K}_3(k)\\
  \vecb{K}_2(k)  \\ \vecb{K}_5(k) \\ \vecb{K}_6(k)\\
  \vecb{K}_3(k) \\ \vecb{K}_6(k) \\ \vecb{K}_9(k) 
\end{pmatrix}.
$$
Thus we obtain the following vectorized code for $\MAT{K}_g$ :
\begin{center}
\begin{minipage}[t]{\ListingWidth}
\MonMatlabNonumber
\begin{lstlisting}
K1 = 3*W1+W2+W3;
K2 = W1+W2+W3/2;
K3 = W1+W2/2+W3;
K5 = W1+3*W2+W3;
K6 = W1/2+W2+W3;
K9 = W1+W2+3*W3;
Kg = [K1;K2;K3;K2;K5;K6;K3;K6;K9];
\end{lstlisting}
\end{minipage}
\end{center}
We represent this technique on Figure~\ref{AssMatWV2}.
\begin{figure}[H]
\definecolor{ffqqqq}{rgb}{1,0,0}
\definecolor{cqcqcq}{rgb}{0.75,0.75,0.75}
\definecolor{color1}{rgb}{0.75,0.75,0.}
\colorlet{colorK1}{Orange}
\colorlet{colorK2}{Red}
\colorlet{colorK3}{Brown}
\colorlet{colorK5}{SeaGreen}
\colorlet{colorK6}{DodgerBlue}
\colorlet{colorK9}{MidnightBlue}
\newcommand{\DrawHorTabBis}[5]{%
  \draw [#4] (#1,#2)-- ({#1+#3},{#2})-- ({#1 + #3},{#2 +1}) -- ({#1}, {#2+1}) -- (#1,#2);
  \draw [dotted] ({#1+1},{#2}) -- ({#1+1},{#2+1});
  \draw [dotted] ({#1+2},{#2}) -- ({#1+2},{#2+1});
  
  \draw [loosely dotted] (#1+2.5,#2+0.5) -- ({#1+3.5},{#2+0.5});
  \draw [dotted] ({#1+4},{#2}) -- ({#1+4},{#2+1});
  \draw [color=black](#1+4.5,#2+0.5) node[anchor=center,scale=0.5] {#5};
  \draw [dotted] ({#1+5},{#2}) -- ({#1+5},{#2+1});
  \draw [loosely dotted] (#1+5.5,#2+0.5) -- ({#1+#3-2.5},{#2+0.5});
  \draw [dotted] ({#1+#3-2},{#2}) -- ({#1+#3-2},{#2+1});
  \draw [dotted] ({#1+#3-1},{#2}) -- ({#1+#3-1},{#2+1});
}

\newcommand{\DrawHorTab}[4]{%
  \draw [#4] (#1,#2)-- ({#1+#3},{#2})-- ({#1 + #3},{#2 +1}) -- ({#1}, {#2+1}) -- (#1,#2);
  \draw [dotted] ({#1+1},{#2}) -- ({#1+1},{#2+1});
  \draw [dotted] ({#1+2},{#2}) -- ({#1+2},{#2+1});
  \draw [loosely dotted] (#1+2.5,#2+0.5) -- ({#1+#3-2.5},{#2+0.5});
  \draw [dotted] ({#1+#3-2},{#2}) -- ({#1+#3-2},{#2+1});
  \draw [dotted] ({#1+#3-1},{#2}) -- ({#1+#3-1},{#2+1});
}
\newcommand{\DrawHorInd}[3]{%
\draw [color=red](#1+0.5,#2+1) node[anchor=center,scale=0.6] {$ 1 $};
 \draw [color=red](#1+1.5,#2+1) node[anchor=center,scale=0.6] {$ 2 $};
 \draw [color=red](#1+2.5,#2+1) node[anchor=center,scale=0.6] {$ \hdots $};
 \draw [color=red](#1+#3-1.5,#2+1) node[anchor=center,scale=0.6] {$ \hdots $};
 \draw [color=red](#1+#3-0.5,#2+1) node[anchor=center,scale=0.6] {$ \nme $};
}

\newcommand{\DrawHorIndBis}[3]{%
\draw [color=red](#1+0.5,#2+1) node[anchor=center,scale=0.6] {$ 1 $};
 \draw [color=red](#1+1.5,#2+1) node[anchor=center,scale=0.6] {$ 2 $};
 \draw [color=red](#1+2.5,#2+1) node[anchor=center,scale=0.6] {$ \hdots $};
 \draw [color=red](#1+4.5,#2+1) node[anchor=center,scale=0.6] {$ k $};
 \draw [color=red](#1+#3-1.5,#2+1) node[anchor=center,scale=0.6] {$ \hdots $};
 \draw [color=red](#1+#3-0.5,#2+1) node[anchor=center,scale=0.6] {$ \nme $};
}

\newcommand{\DrawHorIndLabel}[4]{%
\draw [color=red](#1+0.5,#2+1) node[anchor=center,scale=0.6] {$ 1 $};
 \draw [color=red](#1+1.5,#2+1) node[anchor=center,scale=0.6] {$ 2 $};
 \draw [color=red](#1+2.5,#2+1) node[anchor=center,scale=0.6] {$ \hdots $};
 \draw [color=red](#1+#3-1.5,#2+1) node[anchor=center,scale=0.6] {$ \hdots $};
 \draw [color=red](#1+#3-0.5,#2+1) node[anchor=center,scale=0.6] {#4};
}

\newcommand{\DrawVerInd}[2]{
 \foreach \i in {1,...,9}
   \draw [color=red](#1,#2+9-\i+0.5) node[anchor=center,scale=0.6] {$ \i $};
}

\newcommand{\DrawKgTer}[2]{%
\DrawHorTabBis{#1}{#2+8}{9}{fill=colorK1!60}{$\vecb{K}_1(k)$}; %
\DrawHorTabBis{#1}{#2+7}{9}{fill=colorK2!60}{$\vecb{K}_2(k)$};
\DrawHorTabBis{#1}{#2+6}{9}{fill=colorK3!60}{$\vecb{K}_3(k)$};
\DrawHorTabBis{#1}{#2+5}{9}{fill=colorK2!60}{$\vecb{K}_2(k)$};
\DrawHorTabBis{#1}{#2+4}{9}{fill=colorK5!60}{$\vecb{K}_5(k)$}; 
\DrawHorTabBis{#1}{#2+3}{9}{fill=colorK6!60}{$\vecb{K}_6(k)$}; 
\DrawHorTabBis{#1}{#2+2}{9}{fill=colorK3!60}{$\vecb{K}_3(k)$}; 
\DrawHorTabBis{#1}{#2+1}{9}{fill=colorK6!60}{$\vecb{K}_6(k)$}; 
\DrawHorTabBis{#1}{#2}{9}{fill=colorK9!60}{$\vecb{K}_9(k)$};
}

\begin{center}
\begin{tikzpicture}[>=latex,scale=0.4]
\draw  [fill=gray!10] (-4,0) -- (10,0) -- (10,6) -- (-4,6) -- (-4,0);
\draw (-0.5,4) node[anchor=east] { $\vecb{T_w}$ };
\DrawHorTab{0}{3.5}{9}{}
\DrawHorIndLabel{0}{4.}{9}{ $\nq $}

\draw (-0.5,1.5) node[anchor=east] { \texttt{areas} };
\DrawHorTab{0}{1.}{7}{}
\DrawHorInd{0}{1.5}{7}

\draw [->,color=black] (3.5,0) to[in=0,out=270] (0.75,-3.5);
\draw [->,color=black] (3.5,0) to[in=0,out=270] (0.75,-5.5);
\draw [->,color=black] (3.5,0) to[in=0,out=270] (0.75,-7.5);

\draw [->,color=black] (3.5,0) to[in=180,out=270] (6.25,-3.5);
\draw [->,color=black] (3.5,0) to[in=180,out=270] (6.25,-5.5);
\draw [->,color=black] (3.5,0) to[in=180,out=270] (6.25,-7.5);

\draw (-0.5,-3.5) node[anchor=west] { $\vecb{K}_1$ };
\DrawHorTab{-8}{-4}{7}{fill=colorK1!60}
\DrawHorInd{-8}{-3.5}{7}

\draw (-0.5,-5.5) node[anchor=west] { $\vecb{K}_2$ };
\DrawHorTab{-8}{-6}{7}{fill=colorK2!60}

\draw (-0.5,-7.5) node[anchor=west] { $\vecb{K}_3$ };
\DrawHorTab{-8}{-8}{7}{fill=colorK3!60}

\draw (7.5,-3.5) node[anchor=east] { $\vecb{K}_5$ };
\DrawHorTab{8}{-4}{7}{fill=colorK5!60}
\DrawHorInd{8}{-3.5}{7}

\draw (7.5,-5.5) node[anchor=east] { $\vecb{K}_6$ };
\DrawHorTab{8}{-6}{7}{fill=colorK6!60}

\draw (7.5,-7.5) node[anchor=east] { $\vecb{K}_9$ };
\DrawHorTab{8}{-8}{7}{fill=colorK9!60}

\draw [->,color=colorK1] (-8,-3.5) to[in=180,out=180] (-1.75,-11.5);
\draw [->,color=colorK2] (-8,-5.5) to[in=180,out=180] (-1.75,-12.5);
\draw [->,color=colorK2] (-8,-5.5) to[in=180,out=180] (-1.75,-14.5);
\draw [->,color=colorK3] (-8,-7.5) to[in=180,out=180] (-1.75,-13.5);
\draw [->,color=colorK3] (-8,-7.5) to[in=180,out=180] (-1.75,-17.5);

\draw [->,color=colorK5] (15,-3.5) to[in=0,out=0] (8.75,-15.5);
\draw [->,color=colorK6] (15,-5.5) to[in=0,out=0] (8.75,-16.5);
\draw [->,color=colorK6] (15,-5.5) to[in=0,out=0] (8.75,-18.5);
\draw [->,color=colorK9] (15,-7.5) to[in=0,out=0] (8.75,-19.5);

\draw (3.5,-21.5) node[anchor=center] { $\MAT{K}_g$ };
\DrawKgTer{-1}{-20}
\DrawHorIndBis{-1}{-11.5}{9}
\DrawVerInd{-1.5}{-20}
\DrawVerInd{8.5}{-20}

\end{tikzpicture}
\end{center}

\caption{Weighted mass matrix assembly\label{AssMatWV2} - Version 2}
\end{figure}
Finally, the complete vectorized code using element matrix symmetry is :
\begin{center}
\begin{minipage}[t]{\ListingWidth}
\MonMatlab
\begin{lstlisting}[caption=Optimized assembly - version 2 (Weighted mass matrix)]
function M=MassWAssemblingP1OptV2(nq,nme,me,areas,Tw)
W1=Tw(me(1,:)).*areas/30;
W2=Tw(me(2,:)).*areas/30;
W3=Tw(me(3,:)).*areas/30;
Kg=zeros(9,nme);
Kg(1,:) = 3*W1+W2+W3;
Kg(2,:) = W1+W2+W3/2;
Kg(3,:) = W1+W2/2+W3;
Kg(5,:) = W1+3*W2+W3;
Kg(6,:) = W1/2+W2+W3;
Kg(9,:) = W1+W2+3*W3;
Kg([4, 7, 8],:)=Kg([2, 3, 6],:);
clear W1 W2 W3
Ig = me([1 2 3 1 2 3 1 2 3],:);
Jg = me([1 1 1 2 2 2 3 3 3],:);
M = sparse(Ig(:),Jg(:),Kg(:),nq,nq);
\end{lstlisting}
\end{minipage}
\end{center}

\subsection{Stiffness matrix assembly}
\label{sec:stiffmatP1}
The vertices of the triangle $T_k$ are $\q^{\me(\il,k)}$, $1 \le \il \le 3$.
We define $\VEC{u}^k=\q^{\me(2,k)} - \q^{\me(3,k)},$
$\VEC{v}^k=\q^{\me(3,k)} - \q^{\me(1,k)}$ and $\VEC{w}^k=\q^{\me(1,k)} - \q^{\me(2,k)}$.
Then, the element stiffness matrix $\StiffElem(T_k)$ associated to $T_k$ is defined
by~\eqref{StiffMatElem} with $\VEC{u}=\VEC{u}^k$, $\VEC{v}=\VEC{v}^k$,
$\VEC{w}=\VEC{w}^k$ and $T=T_k$.  Let $\vecb{K}_{1}$, $\vecb{K}_{2}$,
$\vecb{K}_{3}$, $\vecb{K}_{5}$, $\vecb{K}_{6}$ and $\vecb{K}_{9}$ be
six arrays of length $\nme$ such that, for all $k\in\ENS{1}{\nme},$
$$
\begin{array}{lclclclclcl}
\vecb{K}_{1}(k)&=&\displaystyle\frac{\DOT{\VEC{u}^k}{\VEC{u}^k}}{4|T_k|},& & 
\vecb{K}_{2}(k)&=&\displaystyle\frac{\DOT{\VEC{u}^k}{\VEC{v}^k}}{4|T_k|},& & 
\vecb{K}_{3}(k)&=&\displaystyle\frac{\DOT{\VEC{u}^k}{\VEC{w}^k}}{4|T_k|},\\ \\
\vecb{K}_{5}(k)&=&\displaystyle\frac{\DOT{\VEC{v}^k}{\VEC{v}^k}}{4|T_k|},& & 
\vecb{K}_{6}(k)&=&\displaystyle\frac{\DOT{\VEC{v}^k}{\VEC{w}^k}}{4|T_k|},& & 
\vecb{K}_{9}(k)&=&\displaystyle\frac{\DOT{\VEC{w}^k}{\VEC{w}^k}}{4|T_k|}.
\end{array}
$$
With these arrays, the vectorized assembly method is similar to the one shown
in Figure~\ref{AssMatWV2} and the corresponding code is :
\begin{center}
\begin{minipage}[t]{\ListingWidth}
\MonMatlabNonumber
\begin{lstlisting}
Kg = [K1;K2;K3;K2;K5;K6;K3;K6;K9];
S = sparse(Ig(:),Jg(:),Kg(:),nq,nq);
\end{lstlisting}
\end{minipage}
\end{center}
We now describe the vectorized computation of these six arrays. We introduce
the $2$-by-$\nme$ arrays $\vecb{q}_\il, \,\il\in\ENS{1}{3},$
containing the coordinates of the three vertices  of the triangle~$T_k :$ 
$$
\vecb{q}_\il(1,k)=\q(1,\me(\il,k)),\ \quad \vecb{q}_\il(2,k)=\q(2,\me(\il,k)).
$$
We give below the code to compute these arrays, in a non-vectorized form (on the left) and
in a vectorized form (on the right) :\\
\begin{center}
\begin{minipage}[t]{0.71\textwidth}
\MonMatlabNonumber
\begin{lstlisting}
q1=zeros(2,nme);q2=zeros(2,nme);q3=zeros(2,nme);
for k=1:nme
  q1(:,k)=q(:,me(1,k));
  q2(:,k)=q(:,me(2,k));
  q3(:,k)=q(:,me(3,k));
end
\end{lstlisting}
\end{minipage}\hspace{0.05\textwidth}
\begin{minipage}[t]{0.23\textwidth}
\MonMatlabNonumber
\begin{lstlisting}
q1=q(:,me(1,:));
q2=q(:,me(2,:));
q3=q(:,me(3,:));
\end{lstlisting}
\end{minipage}
\end{center}
We trivially obtain the $2$-by-$\nme$ arrays $\VEC{u},$ $\VEC{v}$ and $\VEC{w}$ whose $k$-th column is
$\VEC{u}^k, \ \VEC{v}^k$ and $\VEC{w}^k$ respectively.

The associated code is :
\begin{center}
\begin{minipage}[t]{\ListingWidth}
\MonMatlabNonumber
\begin{lstlisting}
u=q2-q3;
v=q3-q1;
w=q1-q2;
\end{lstlisting}
\end{minipage}
\end{center}
The operators \lstinline{.*}, \lstinline{./} (element-wise arrays multiplication and division)
and the function \lstinline{sum(.,1)} (row-wise sums) allow to compute all arrays.
For example, $\vecb{K}_{2}$ is computed using the following vectorized code :
\begin{center}
\begin{minipage}[t]{\ListingWidth}
\MonMatlabNonumber
\begin{lstlisting}
K2=sum(u.*v,1)./(4*areas);
\end{lstlisting}
\end{minipage}
\end{center}
Then, the complete vectorized function using element matrix symmetry is :
\begin{center}
\begin{minipage}[t]{\ListingWidth}
\MonMatlab
\begin{lstlisting}[caption=Optimized matrix assembly code - version 2 (Stiffness matrix)]
function S=StiffAssemblingP1OptV2(nq,nme,q,me,areas)
q1 =q(:,me(1,:)); q2 =q(:,me(2,:)); q3 =q(:,me(3,:));
u = q2-q3; v=q3-q1; w=q1-q2;
areas4=4*areas;
Kg=zeros(9,nme);
Kg(1,:)=sum(u.*u,1)./areas4; % K1 
Kg(2,:)=sum(v.*u,1)./areas4; % K2 
Kg(3,:)=sum(w.*u,1)./areas4; % K3 
Kg(5,:)=sum(v.*v,1)./areas4; % K5 
Kg(6,:)=sum(w.*v,1)./areas4; % K6 
Kg(9,:)=sum(w.*w,1)./areas4; % K9 
Kg([4, 7, 8],:)=Kg([2, 3, 6],:);
clear q1 q2 q3 areas4 u v w
Ig = me([1 2 3 1 2 3 1 2 3],:);
Jg = me([1 1 1 2 2 2 3 3 3],:);
S = sparse(Ig(:),Jg(:),Kg(:),nq,nq);
\end{lstlisting}
\end{minipage}
\end{center}
\vspace{-1.1mm}

\subsection{Comparison with FreeFEM++}
\begin{figure}[H]
\dblimageps{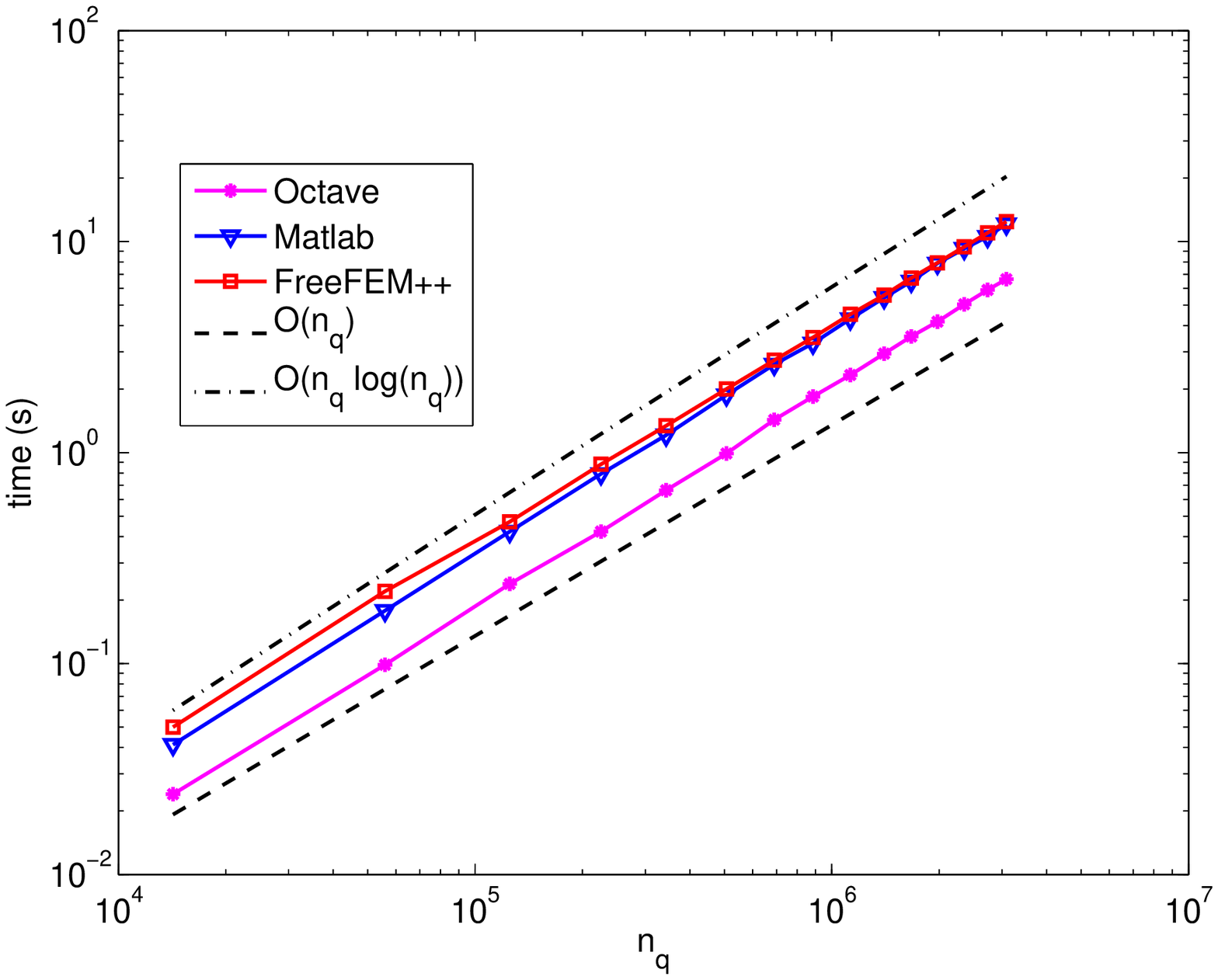}{\zoomdeux}
{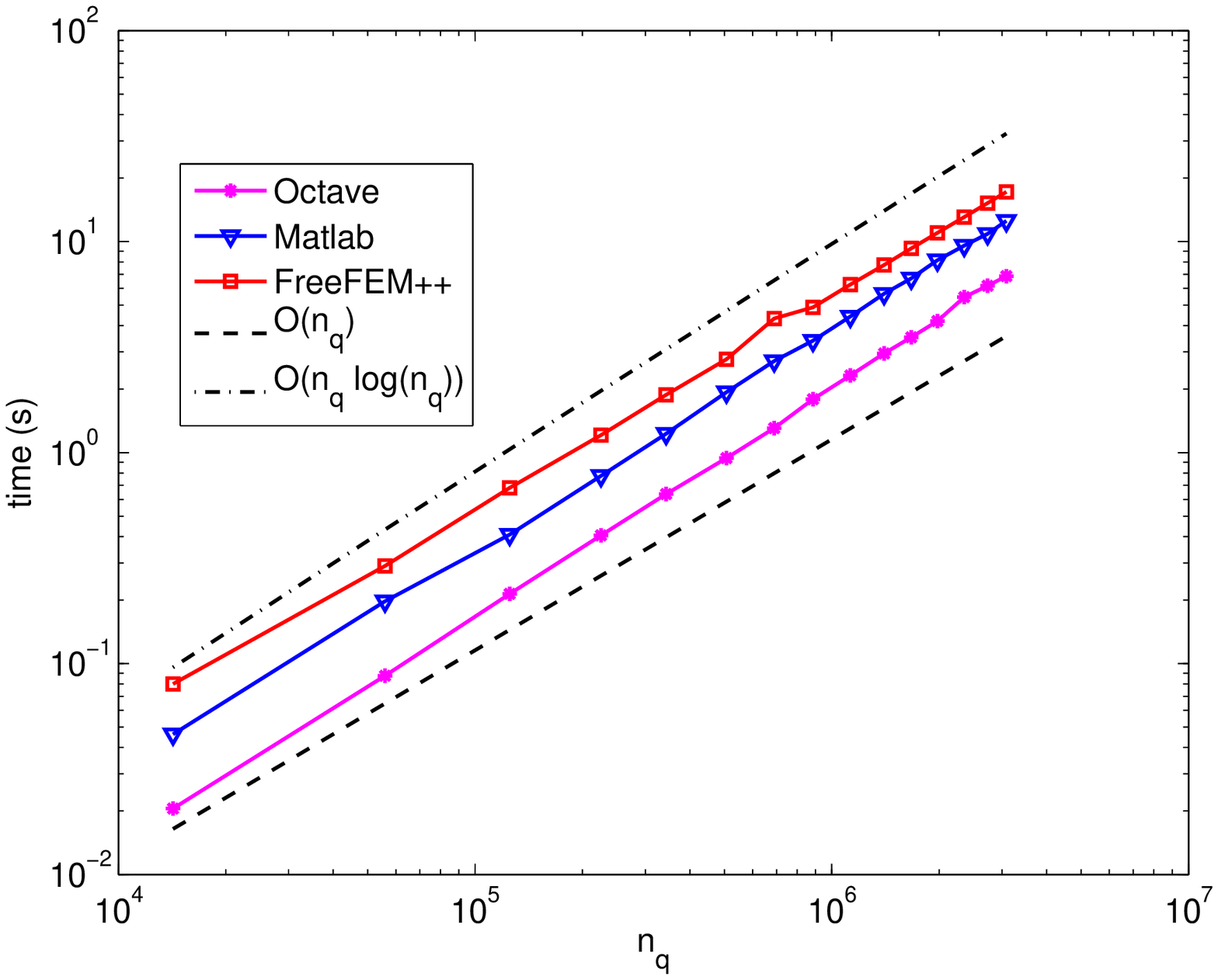}{\zoomdeux}
\imageps{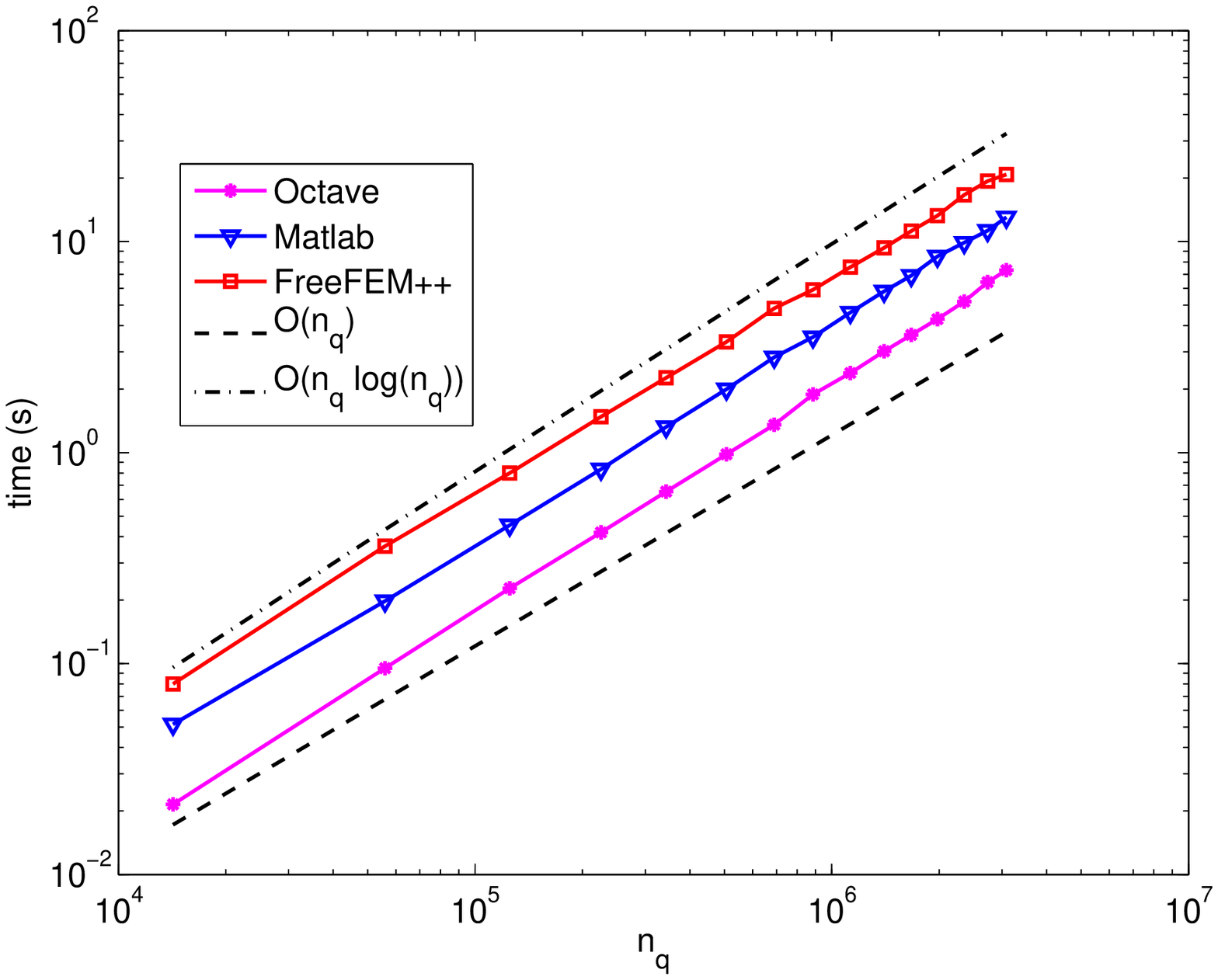}{\zoomdeux}
\caption{Comparison of the matrix assembly codes : \texttt{OptV2} in Matlab/Octave and FreeFEM++,
  for the mass (top left), weighted mass (top right) and stiffness (bottom) matrices.
  \label{OptV2VsFreeFEM}}
\end{figure}

On Figure~\ref{OptV2VsFreeFEM}, we show the computation times of the FreeFEM++ and
\texttt{OptV2} Matlab/Octave codes, versus $\nq$.
\vspace{-4mm}
\noindent
The computation times values are given in Appendix~\ref{App:OptV2vsFreeFem}.
The complexity of the Matlab/Octave codes is still linear ($\GrandO(\nq)$) and slightly better
than the one of FreeFEM++.
\vspace{1mm}
\begin{remark}\label{rem:JIT}
We observed that only with the \texttt{OptV2} codes, Octave gives better results than Matlab.
For the other versions of the codes, not fully vectorized, the JIT-Accelerator (Just-In-Time)
of Matlab allows significantly better performances than Octave 
(JIT compiler for GNU Octave is under development).

Furthermore, we can improve Matlab performances using SuiteSparse packages
from T. Davis \cite{Davis:SSP:2012}, which is originally used in Octave. 
In our codes, using \texttt{cs\_sparse} function from SuiteSparse instead of
Matlab \texttt{sparse} function is approximately 1.1 times faster for \texttt{OptV1} version
and 2.5 times for \texttt{OptV2} version (see also Section~\ref{sec:Elasticity}).
\end{remark}

\subsection{Comparison with other matrix assembly codes}
\label{subsec:CompChenetal}
\begin{figure}[H]
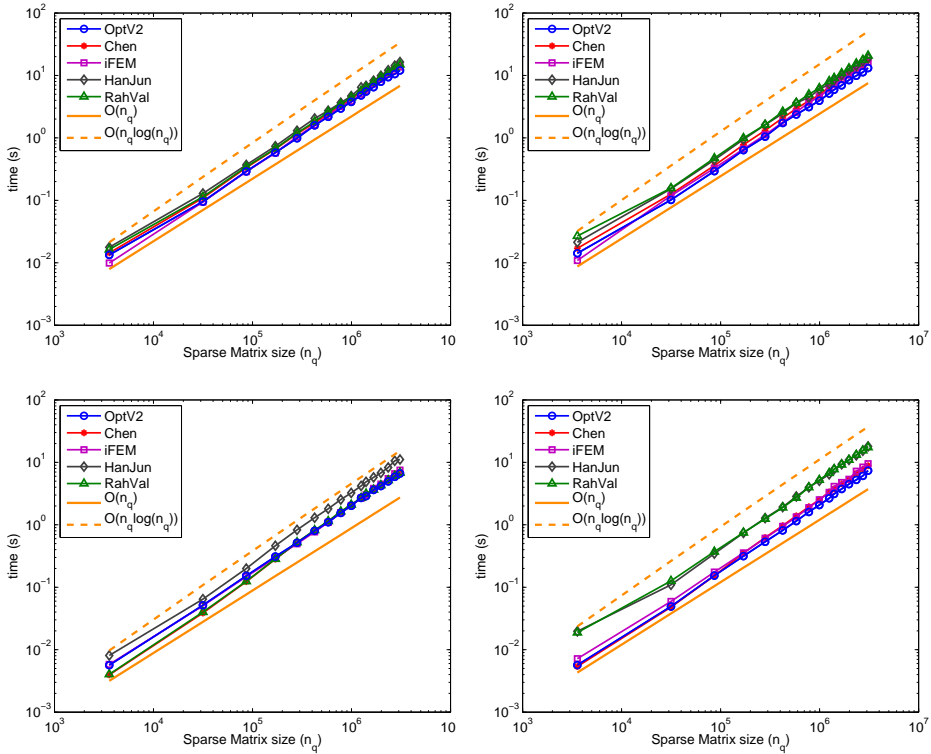

\dblimagepss{MassAssembling_disque4-1_gpuschwarz_\OptvsHanJunvsRahValDATE_computetime}{0.43}
{StiffAssembling_disque4-1_gpuschwarz_\OptvsHanJunvsRahValDATE_computetime}{0.43}
\dblimagepss{MassAssembling_disque4-1_gpuschwarz_\OptvsHanJunvsRahValDATEOctave_computetime}{0.43}
{StiffAssembling_disque4-1_gpuschwarz_\OptvsHanJunvsRahValDATEOctave_computetime}{0.43}
\caption{Comparison of the assembly codes in  Matlab R2012b (top) and Octave~3.6.3 (bottom):
  \texttt{OptV2} and~\cite{Chen:PFE:2011,Chen:iFEM:2013,Hannukainen:IFE:2012,Rahman:FMA:2011},
  for the mass (left) and stiffness (right) matrices.\label{OptV2vsHJRVMO}}
\end{figure}

We compare the matrix assembly codes proposed by
L.~Chen~\cite{Chen:PFE:2011,Chen:iFEM:2013}, A.~Hannukainen and
M.~Juntunen~\cite{Hannukainen:IFE:2012} and T.~Rahman and
J.~Valdman~\cite{Rahman:FMA:2011} to the \texttt{OptV2} version
developed in this paper, for the mass and stiffness matrices.  The
domain $\Omega$ is the unit disk.  The computations have been done on
our reference computer.
On Figure~\ref{OptV2vsHJRVMO}, with Matlab (top) and Octave (bottom),
we show the computation times versus the number of vertices of the mesh, for these different codes.
The associated values are given in Tables \ref{MassMatAll} to \ref{StiffOctAll}.
For large sparse matrices, our \texttt{OptV2} version allows gains in computational time 
of $5\%$ to $20\%$, compared to the other vectorized codes (for sufficiently large meshes).
\vspace{-1mm}
\begin{table}[H]\scriptsize
\begin{center}

\end{center}
\caption{\label{MassMatAll}Computational cost, in Matlab (R2012b), of the \texttt{Mass} matrix assembly versus $\nq,$ with the \texttt{OptV2} version (column $2$) and with the codes in \cite{Chen:PFE:2011,Chen:iFEM:2013,Hannukainen:IFE:2012,Rahman:FMA:2011} (columns $3$-$6$) : time in seconds (top value) and speedup (bottom value). The speedup reference is \texttt{OptV2} version.}
\end{table}

\vspace{-0.1cm}
\begin{table}[H]\scriptsize
\begin{center}
%
\end{center}
\caption{\label{StiffMatAll}Computational cost, in Matlab (R2012b), of the \texttt{Stiffness} matrix assembly versus $\nq,$ with the \texttt{OptV2} version (column $2$) and with the codes in \cite{Chen:PFE:2011,Chen:iFEM:2013,Hannukainen:IFE:2012,Rahman:FMA:2011} (columns $3$-$6$) : time in seconds (top value) and speedup (bottom value). The speedup reference is \texttt{OptV2} version.}
\end{table}

\newpage
\begin{table}[H]\scriptsize
\begin{center}
%
\end{center}
\caption{\label{MassOctAll}Computational cost, in Octave (3.6.3), of the \texttt{Mass} matrix assembly versus $\nq,$ with the \texttt{OptV2} version (column $2$) and with the codes in \cite{Chen:PFE:2011,Chen:iFEM:2013,Hannukainen:IFE:2012,Rahman:FMA:2011} (columns $3$-$6$) : time in seconds (top value) and speedup (bottom value). The speedup reference is \texttt{OptV2} version.}
\end{table}

\vspace{-0.1cm}
\begin{table}[H]\scriptsize
\begin{center}
%
 \\ \hline
\end{tabular}
\end{center}
\caption{\label{StiffOctAll}Computational cost, in Octave (3.6.3), of the \texttt{Stiffness} matrix assembly versus $\nq,$ with the \texttt{OptV2} version (column $2$) and with the codes in \cite{Chen:PFE:2011,Chen:iFEM:2013,Hannukainen:IFE:2012,Rahman:FMA:2011} (columns $3$-$6$) : time in seconds (top value) and speedup (bottom value). The speedup reference is \texttt{OptV2} version.}
\end{table}

\newpage
%
\section{Extension to linear elasticity}
\label{sec:Elasticity}
In this part we extend the codes of the previous sections to a linear
elasticity matrix assembly.

Let $\HUnHD{\DOMH}$ be the finite dimensional space spanned by
the $P_1$ Lagrange basis functions $\{\varphi_i\}_{i\in\ENS{1}{\nq}}$.
Then, the space $(\HUnHD{\DOMH})^2$ is spanned by
$\mathcal{B}=\{\BasisFuncTwoD_l\}_{1\le l \le 2\nq}$, with
$\BasisFuncTwoD_{2i-1}= \begin{pmatrix} \FoncBase_i \\ 0 \end{pmatrix}$,
$\BasisFuncTwoD_{2i}= \begin{pmatrix} 0 \\ \FoncBase_i  \end{pmatrix}$, $1\le i \le \nq$.

The example we consider is the elastic stiffness matrix $\StiffElas$, defined by
\begin{equation*}
\StiffElas_{m,l}=\int_{\DOMH} \Odv^t(\BasisFuncTwoD_m) 
\Ocv(\BasisFuncTwoD_l)dT, \ \forall (m,l)\in\ENS{1}{2\,\nq}^2,
\end{equation*}
where $\Ocv=(\Occ_{xx},\Occ_{yy},\Occ_{xy})^t$ and $\Odv=(\Odc_{xx},\Odc_{yy},2\Odc_{xy})^t$
are the elastic stress and strain tensors respectively. 
We consider here linearized elasticity with small strain hypothesis (see for example \cite{Dhatt:FEM:2012}).
Consequently, let $\mathcal{D}$ be the differential operator which links displacements $\VEC{u}$ to strains:
\begin{equation*}\label{OpDiffD}
\Odv(\VEC{u}) = \mathcal{D}(\VEC{u}) = \frac 12 \left(\GRAD(\VEC{u})+\GRAD^t(\VEC{u})\right).
\end{equation*}
This gives, in vectorial form and after reduction to the plane,
\begin{equation*}\label{OpDiffPlan}
\mathcal{D}=\begin{pmatrix}
\DP{}{x} & 0 \\
0 & \DP{}{y} \\
\DP{}{y} & \DP{}{x}
\end{pmatrix}.
\end{equation*}
For the constitutive equation, Hooke's law is used and the material is supposed to be isotropic.
Thus, the elasticity tensor denoted by $\mathbb{C}$ becomes a 3-by-3 matrix and can be defined
by the Lamé parameters $\lambda$ and $\mu$, which are supposed constant on $\Omega$ and satisfying $\lambda+\mu >0$. Thus, the constitutive equation writes
\begin{equation*}\label{eqCompPlan}
\Ocv=\mathbb{C}\Odv=
\begin{pmatrix}
\lambda+2\mu & \lambda&  0\\
\lambda & \lambda+2\mu &  0\\
 0 & 0 & \mu 
\end{pmatrix}\Odv.
\end{equation*}
Using the triangulation $\DOMH$ of $\DOM$, we have
\begin{equation*}
\StiffElas_{m,l}=\sum_{k=1}^{\nme}\StiffElas_{m,l}(T_k), 
 \quad \mbox{with} \quad
\StiffElas_{m,l}(T_k)=\int_{T_k} \Odv^t(\BasisFuncTwoD_m) 
\Ocv(\BasisFuncTwoD_l)dT, \ \ \forall (m,l)\in\ENS{1}{2\nq}^2.
\end{equation*}
Let $\mathcal{I}_k=[ 2\me(1,k)-1,2\me(1,k), \ 2\me(2,k)-1,2\me(2,k), \ 2\me(3,k)-1,2\me(3,k)].$
Due to the support of functions $\BasisFuncTwoD_l$, we have
$\forall (l,m)\in  \left(\ENS{1}{\nq}\backslash \mathcal{I}_k\right)^2,
\ \StiffElas_{m,l}(T_k)=0.$
\noindent
Thus, we only have to compute 
$\StiffElas_{m,l}(T_k), \ \forall (m,l) \in \mathcal{I}_k \times \mathcal{I}_k$,
the other terms being zeros. We denote by
$\StiffElasElem_{\il,\jl}(T_k)=\StiffElas_{\mathcal{I}_k(\il),\mathcal{I}_k(\jl)}(T_k),
\ \forall (\il,\jl) \in \ENS{1}{6}^2.$ Therefore, we introduce
$\Local{\mathcal{B}}(T_k)=\{\BasisFuncTwoDL_\il\}_{1\le \il\le 6}$ the
local basis associated to a triangle $T_k$ with
$\BasisFuncTwoDL_\il=\BasisFuncTwoD_{\mathcal{I}_k(\il)},$
$1 \le \il \le 6.$
We thus have
$\BasisFuncTwoDL_{2\gamma-1}= \begin{pmatrix} \FoncBaseLocale_\gamma
  \\ 0 \end{pmatrix},\ \BasisFuncTwoDL_{2\gamma}= \begin{pmatrix} 0
  \\ \FoncBaseLocale_\gamma \end{pmatrix}
\ 1 \le \gamma \le 3.$

\begin{figure}[H]
\dblimageps{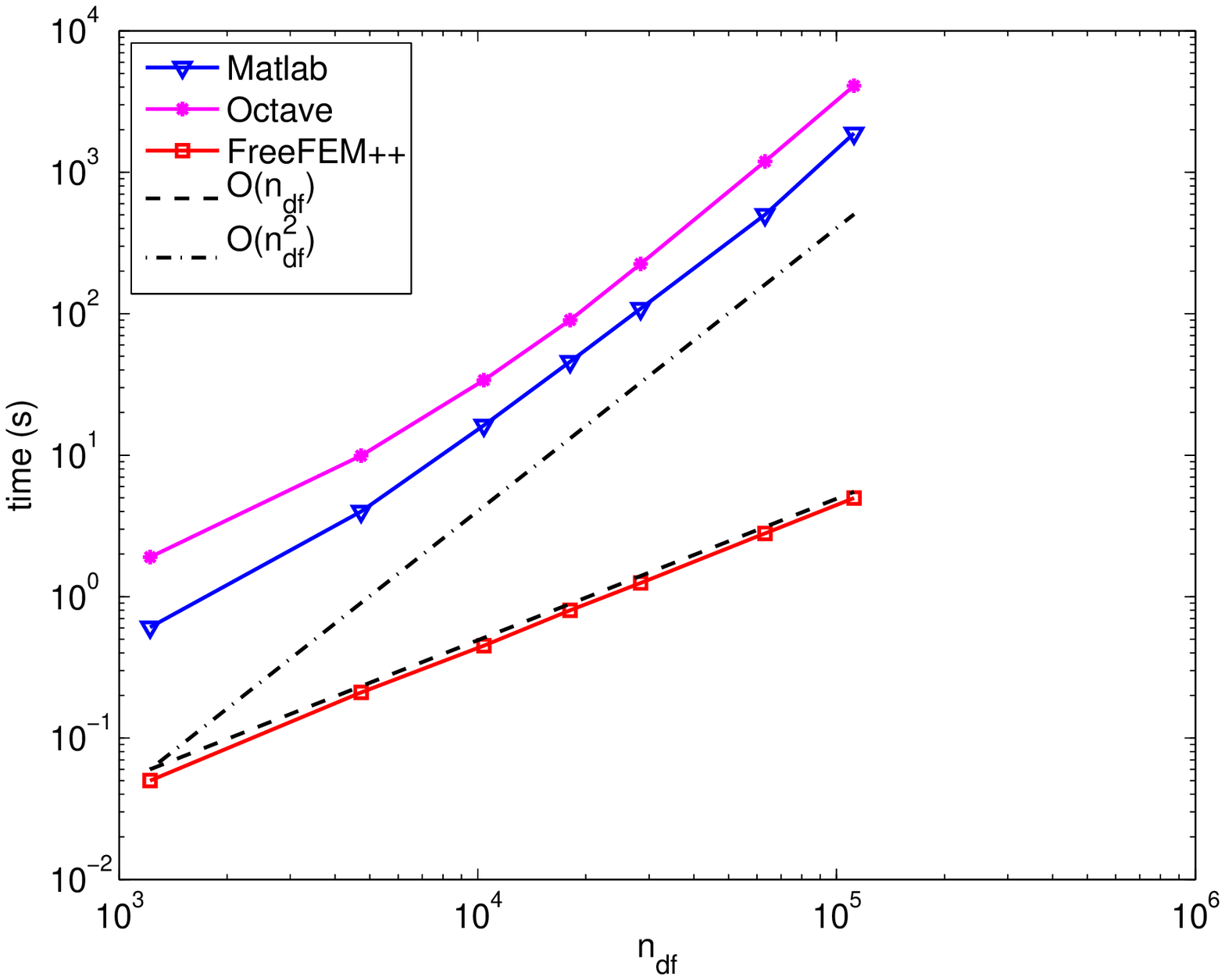}{\zoom}{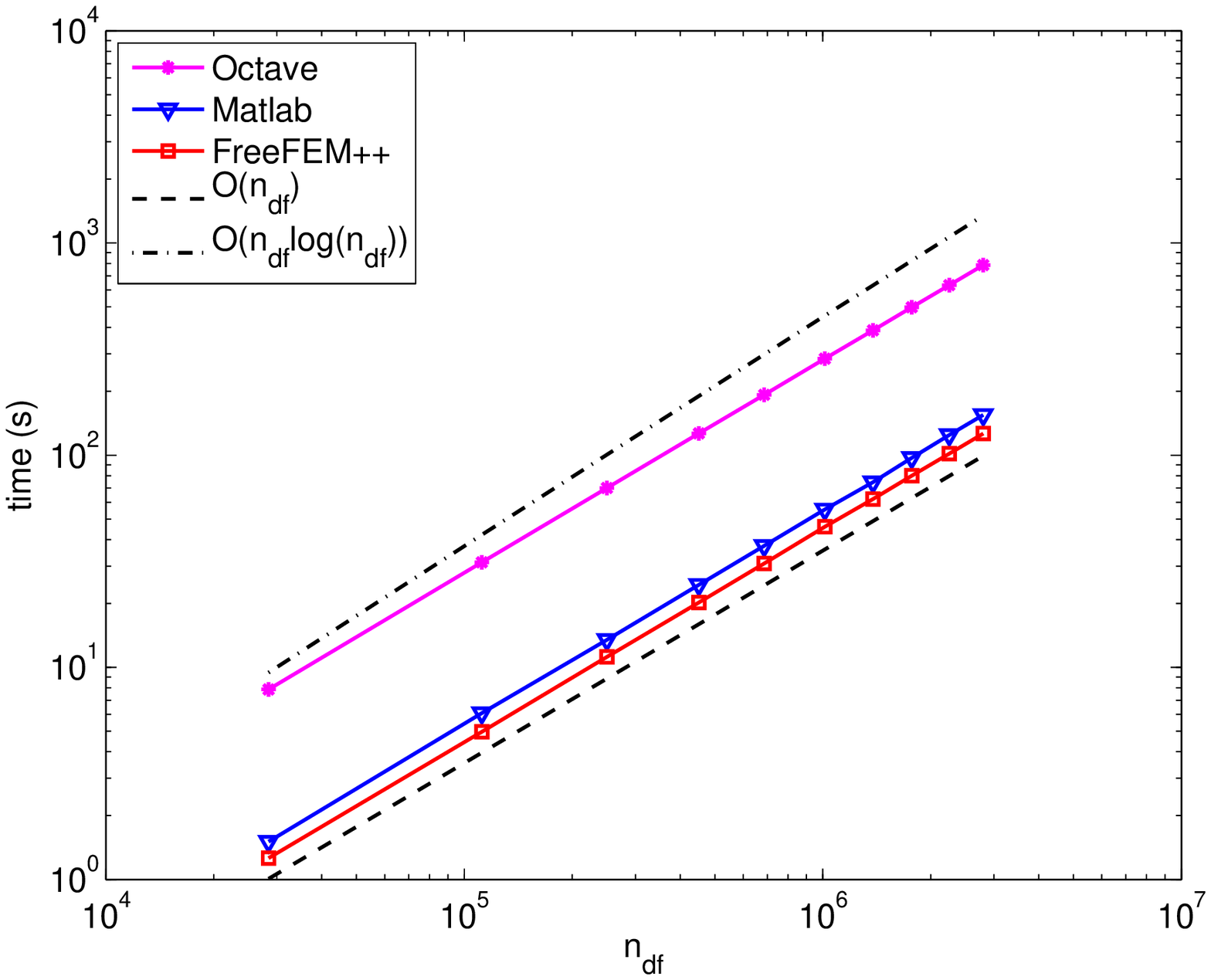}{\zoom}
\caption{Comparison of the matrix assembly codes : usual assembly (left) and \texttt{OptV1} (right)
  in Matlab/Octave and FreeFEM++, for stiffness elasticity matrix.
  \label{VFBaseAndOptV1VsFreeFEM}}
\end{figure}

The element stiffness matrix $\StiffElasElem$ is given by
\begin{equation*}
\StiffElasElem_{\il,\jl}(T_k)=\int_{T_k} \Odv^t(\FoncBaseDeuxDLocale_\il) 
\mathbb{C}\Odv(\FoncBaseDeuxDLocale_\jl)dT, \quad \forall (\il,\jl)\in\ENS{1}{6}^2.
\end{equation*}
Denoting, as in Section~\ref{sec:stiffmatP1},
\begin{equation*}\label{eq:uk-vk-wk}
\VEC{u}^k=\q^{\me(2,k)} - \q^{\me(3,k)}, \quad \VEC{v}^k=\q^{\me(3,k)} - \q^{\me(1,k)},\quad \text{and} \quad
\VEC{w}^k=\q^{\me(1,k)} - \q^{\me(2,k)},
\end{equation*}
with $\q^{\me(\il,k)}, \ 1 \le \il \le 3,$ the three vertices of $T_k$, 
then the gradients of the local functions $\FoncBaseLocale_\il^k=\FoncBase_{\me(\il,k)|_{T_k}}, \ 1 \le \il \le 3$,
associated to $T_k$, are constants and given respectively by
\begin{equation}\label{eq:GradFoncBas}
 \nabla \FoncBaseLocale_1^k=\frac{1}{2|T_k|}\begin{pmatrix} u^k_2 \\ -u^k_1\end{pmatrix},\ 
 \nabla \FoncBaseLocale_2^k=\frac{1}{2|T_k|}\begin{pmatrix} v^k_2 \\ -v^k_1\end{pmatrix},\ 
 \nabla \FoncBaseLocale_3^k=\frac{1}{2|T_k|}\begin{pmatrix} w^k_2 \\ -w^k_1\end{pmatrix}.
\end{equation}
So, we can rewrite 
the matrix $\StiffElasElem(T_k)$ in the form
\begin{equation*}\label{eq:RigiditeElasElemBa}
\StiffElasElem(T_k)=|T_k| \MATT{B}_k\MAT{C}\MAT{B}_k,
\end{equation*}
where
\begin{eqnarray*}
 \MAT{B}_k=\left(\begin{array}{ccccc}
\Odv(\FoncBaseDeuxDLocale_1) & | & \hdots & | & \Odv(\FoncBaseDeuxDLocale_6)
\end{array}\right)=
\frac{1}{2|T_k|}
\begin{pmatrix} 
u_2^k & 0 & v_2^k & 0 & w_2^k & 0\\
0 & -u^k_1 & 0 & -v^k_1 & 0 & -w^k_1\\
-u_1^k & u_2^k & -v_1^k & v_2^k & -w_1^k & w_2^k
\end{pmatrix}.
\end{eqnarray*}
We give the Matlab/Octave code for computing $\StiffElasElem(T_k)$:
\begin{center}
\begin{minipage}[t]{\ListingWidth}
\begin{lstlisting}[caption=Element matrix code (elastic stiffness matrix),label=ElemElas:OptV0]
function Ke=ElemStiffElasMatP1(qm,area,C)
% qm=[q1,q2,q3]
u=qm(:,2)-qm(:,3);   
v=qm(:,3)-qm(:,1); 
w=qm(:,1)-qm(:,2);
B=[u(2),0,v(2),0,w(2),0; ...
   0,-u(1),0,-v(1),0,-w(1); ...
   -u(1),u(2),-v(1),v(2),-w(1),w(2)];
Ke=B'*C*B/(4*area);
\end{lstlisting}
\end{minipage}
\end{center}
Then, the classical matrix assembly code using the element matrix $\StiffElasElem(T_k)$
with a loop through the triangles is
\begin{center}
\begin{minipage}[t]{0.92\textwidth}
\begin{lstlisting}[caption=Classical matrix assembly code (elastic stiffness matrix),label=AssemblageElas:basique]
function K=StiffElasAssemblingP1(nq,nme,q,me,areas,lam,mu)
K=sparse(2*nq,2*nq);
C=[lam+2*mu,lam,0;lam,lam+2*mu,0;0,0,mu];
for k=1:nme
    MatElem=ElemStiffElasMatP1(q(:,me(:,k)),areas(k),C); 
    I=[2*me(1,k)-1, 2*me(1,k), 2*me(2,k)-1, ...
       2*me(2,k), 2*me(3,k)-1, 2*me(3,k)];
    for il=1:6
        for jl=1:6
            K(I(il),I(jl))=K(I(il),I(jl))+MatElem(il,jl);           
        end
    end   
end
\end{lstlisting}
\end{minipage}
\end{center}
On Figure~\ref{VFBaseAndOptV1VsFreeFEM} on the left, we show the computation times (in seconds)
versus the matrix size $\ndf=2\nq$, for the classical matrix assembly code and the FreeFEM++ code
given in Listing \ref{AssembElas:FF}.
We observe that the complexity is $\GrandO(\ndfs)$ for the Matlab/Octave codes,
while the complexity seems to be $\GrandO(\ndf)$ for FreeFEM++.

\subsection{Optimized matrix assembly - version 1 (OptV1)}
\label{sec:CodeStiffElasOptV1}
We define the three local arrays $\vecb{K}^e_k,$ $\vecb{I}^e_k$ and   
$\vecb{J}^e_k$ of $36$ elements by
\begin{center}
\begin{tabular}{lcl}
$\vecb{K}^e_k$ &:& elements of the matrix $\StiffElasElem(T_k)$ stored column-wise,\\
$\vecb{I}^e_k$ &:& global row indices associated to the elements stored in $\vecb{K}^e_k$,\\
$\vecb{J}^e_k$ &:& global column indices associated to the elements stored in $\vecb{K}^e_k.$
\end{tabular}
\end{center}
Using the definition of $\mathcal{I}_k$ in the introduction of Section~\ref{sec:Elasticity}, we have
\begin{equation*}
\forall (\il,\jl)\in\ENS{1}{6},\ 
\left\{\begin{array}{lcl}
\vecb{K}^e_k\bigl(6(\jl-1)+\il \bigr)&=& \StiffElasElem_{\il,\jl}(T_k),\\
\vecb{I}^e_k\bigl(6(\jl-1)+\il \bigr)&=& \mathcal{I}_k(\il),\\
\vecb{J}^e_k\bigl(6(\jl-1)+\il \bigr)&=& \mathcal{I}_k(\jl).
\end{array}\right.
\end{equation*}
Thus, from the matrix $\StiffElasElem(T_k)=(K_{i,j}^k)_{1\le i,j \le 6}$, we obtain
$$
\begin{array}{cclc*{10}{c}r}
\vecb{K}^e_k&=& \bigl(  K_{1,1}^k & \hdots & K_{6,1}^k & , & K_{1,2}^k& \hdots & 
                    K_{6,2}^k & , & \hdots & , &  K_{1,6}^k& \hdots & K_{6,6}^k  \bigr)\\
\vecb{I}^e_k&=& \bigl(  \mathcal{I}_k(1) & \hdots & \mathcal{I}_k(6) & , & \mathcal{I}_k(1)& \hdots & 
                    \mathcal{I}_k(6) & , & \hdots & , &  \mathcal{I}_k(1)& \hdots & \mathcal{I}_k(6)  \bigr)\\
\vecb{J}^e_k&=& \bigl(  \mathcal{I}_k(1) & \hdots & \mathcal{I}_k(1) & , & \mathcal{I}_k(2)& \hdots & 
                    \mathcal{I}_k(2) & , & \hdots & , &  \mathcal{I}_k(6)& \hdots & \mathcal{I}_k(6) \bigr)
\end{array}
$$
We give below the associated Matlab/Octave code :
\vspace{-2mm}
\begin{center}
\begin{minipage}[t]{0.95\textwidth}
\MonMatlab
\begin{lstlisting}[caption=Optimized matrix assembly code - version 1 (elastic stiffness matrix),label=AssemblageElas:OptV1]
function K=StiffElasAssemblingP1OptV1(nq,nme,q,me,areas,lam,mu)
Ig=zeros(36*Th.nme,1);Jg=zeros(36*Th.nme,1);
Kg=zeros(36*Th.nme,1);
kk=1:36;
C=[lam+2*mu,lam,0;lam,lam+2*mu,0;0,0,mu];
for k=1:nme
  Me=ElemStiffElasMatP1(q(:,me(:,k)),areas(k),C);
  I=[2*me(1,k)-1, 2*me(1,k), 2*me(2,k)-1, ...
     2*me(2,k), 2*me(3,k)-1, 2*me(3,k)];
  je=ones(6,1)*I; ie=je';
  Ig(kk)=ie(:); Jg(kk)=je(:);
  Kg(kk)=Me(:);
  kk=kk+36;
end
K=sparse(Ig,Jg,Kg,2*nq,2*nq);
\end{lstlisting}
\end{minipage}
\end{center}
\newpage
On Figure~\ref{VFBaseAndOptV1VsFreeFEM} on the right, we show the 
computation times of the \texttt{OptV1} codes in Matlab/Octave and of the FreeFEM++ codes
versus the number of degrees of freedom on the mesh (unit disk). 
The complexity of the Matlab/Octave codes seems now linear (i.e. $\GrandO(\ndf)$) as for FreeFEM++.
Also, FreeFEM++ is slightly faster than Matlab, while much more faster than Octave (about a factor 10).
\begin{center}
\begin{minipage}[t]{0.97\textwidth}
\MonFreeFEMNonumber
\begin{lstlisting}[caption=Matrix assembly code in FreeFEM++ (elastic stiffness matrix),label=AssembElas:FF]
mesh Th(...);
fespace Wh(Th,[P1,P1]);
Wh [u1,u2],[v1,v2];
real lam=...,mu=...;
func C=[[lam+2*mu,lam,0],[lam,lam+2*mu,0],[0,0,mu]];
macro epsilon(ux,uy) [dx(ux),dy(uy),(dy(ux)+dx(uy))]
macro sigma(ux,uy) ( C*epsilon(ux,uy) )
varf vStiffElas([u1,u2],[v1,v2])=
                int2d(Th)(epsilon(u1,u2)'*sigma(v1,v2));
matrix K = vStiffElas(Wh,Wh);
\end{lstlisting}
\end{minipage}
\end{center}
To further improve the efficiency of the matrix assembly code, we introduce now the optimized version 2.

\subsection{Optimized matrix assembly - version 2 (OptV2)}
\label{sec:CodeStiffElasOptV2}
In this version, no loop is used. As in Section \ref{sssec:CodeMassOptV2},
we define three 2d-arrays that allow to store all the element matrices 
as well as their positions in the global matrix. We denote by $\MAT{K}_g,$ $\MAT{I}_g$ and $\MAT{J}_g$ these
$36$-by-$\nme$ arrays, defined $\forall k \in \ENS{1}{\nme},$ $\forall il \in\ENS{1}{36}$ by
\begin{eqnarray*}
\MAT{K}_g(il,k)=\vecb{K}^e_k(il),\quad\quad
\MAT{I}_g(il,k)=\vecb{I}^e_k(il),\quad\quad
\MAT{J}_g(il,k)=\vecb{J}^e_k(il).
\end{eqnarray*}
Thus, the local arrays $\vecb{K}^e_k,$ $\vecb{I}^e_k$ and $\vecb{J}^e_k$ are stored
in the $k$-th column of the global arrays $\MAT{K}_g,$ $\MAT{I}_g$ and $\MAT{J}_g$ respectively.
Once these arrays are determined, the assembly matrix is obtained with the Matlab/Octave command

\vspace{5mm}
\centerline{\lstinline{M = sparse(Ig(:),Jg(:),Kg(:),2*nq,2*nq);}}
\vspace{3mm}
In order to vectorize the computation of $\MAT{I}_g$ and $\MAT{J}_g$,
we generalize the technique introduced in the Remark~\ref{remq:IgJg}
and denote by $\MAT{T}$ the $6 \times \nme$ array defined by
\begin{equation*}\label{def:T}
\MAT{T}=\begin{pmatrix}
\mathcal{I}_1(1) & \hdots & \mathcal{I}_k(1) & \hdots & \mathcal{I}_\nme(1)\\ 
\mathcal{I}_1(2) & \hdots & \mathcal{I}_k(2) & \hdots &\mathcal{I}_\nme(2)\\ 
\vdots & \vdots & \vdots & \vdots & \vdots \\ 
\mathcal{I}_1(6) & \hdots & \mathcal{I}_k(6) & \hdots &\mathcal{I}_\nme(6)
\end{pmatrix}.
\end{equation*}
Then $\MAT{I}_g$ is computed by duplicating $\MAT{T}$ six times, column-wise.
The array $\MAT{J}_g$ is computed from $\MAT{T}$ by duplicating each line, six times, successively.
We give in Listing \ref{listing:IgetJgVF}, the Matlab/Octave vectorized function
which enables to compute $\MAT{I}_g$ and $\MAT{J}_g.$

\begin{figure}[H]
  \centering
\imageps{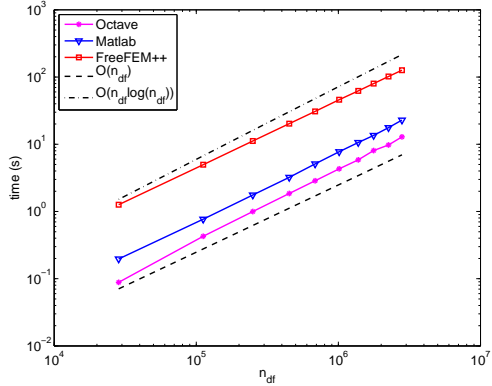}{\zoom}
\caption{Comparison of the assembly code : \texttt{OptV2} in Matlab/Octave and FreeFEM++,
  for elastic stiffness matrix.
  \label{VFOptV2VsFreeFEM}}
\end{figure}

It remains to vectorize the computation of the 2d-array $\MAT{K}_g$.
Using formulas (\ref{eq:GradFoncBas}), for $1 \le \il \le 3$, we define the $2$-by-$\nme$ array $\vecb{G}_\il$,
the $k$-th column of which contains $\nabla \FoncBaseLocale_\alpha^k$.
\begin{center}
\begin{minipage}[t]{0.55\textwidth}
\MonMatlab
\begin{lstlisting}[caption={Vectorized code for computing
      $\MAT{I}_g$ and  $\MAT{J}_g$},label=listing:IgetJgVF] 
function [Ig,Jg]=BuildIgJgP1VF(me)
T= [2*me(1,:)-1; 2*me(1,:); ...
    2*me(2,:)-1; 2*me(2,:); ...
    2*me(3,:)-1; 2*me(3,:)];

ii=[1 1 1 1 1 1; ...
    2 2 2 2 2 2; ...
    3 3 3 3 3 3; ...
    4 4 4 4 4 4; ...
    5 5 5 5 5 5; ...
    6 6 6 6 6 6];

jj=ii';

Ig=T(ii(:),:);
Jg=T(jj(:),:);
\end{lstlisting}
\end{minipage}
\end{center}
Let us focus on the first column of $\MAT{K}^e(T_k)$. It is given by
{\small
\begin{align*}
\MAT{K}^e_{1,1}(T_k)
&=|T_k|\left({\left(\lambda + 2 \, \mu\right)} \frac{\dsp\partial\FoncBaseLocale_{1}}{\dsp\partial x}^{2}
  + \mu \frac{\dsp\partial\FoncBaseLocale_{1}}{\dsp\partial y}^{2} \right), 
& \hspace{-3mm}\MAT{K}^e_{2,1}(T_k)
&=|T_k|\left(\lambda \frac{\dsp\partial\FoncBaseLocale_{1}}{\dsp\partial x} \frac{\dsp\partial\FoncBaseLocale_{1}}{\dsp\partial y}
  + \mu \frac{\dsp\partial\FoncBaseLocale_{1}}{\dsp\partial x} \frac{\dsp\partial\FoncBaseLocale_{1}}{\dsp\partial y} \right)\\
\MAT{K}^e_{3,1}(T_k)
&=|T_k|\left({\left(\lambda + 2 \, \mu\right)} \frac{\dsp\partial\FoncBaseLocale_{1}}{\dsp\partial x}
   \frac{\dsp\partial\FoncBaseLocale_{2}}{\dsp\partial x}
   + \mu \frac{\dsp\partial\FoncBaseLocale_{1}}{\dsp\partial y} \frac{\dsp\partial\FoncBaseLocale_{2}}{\dsp\partial y} \right), 
&\hspace{-3mm}\MAT{K}^e_{4,1}(T_k)
&=|T_k|\left(\lambda \frac{\dsp\partial\FoncBaseLocale_{1}}{\dsp\partial x} \frac{\dsp\partial\FoncBaseLocale_{2}}{\dsp\partial y}
   + \mu \frac{\dsp\partial\FoncBaseLocale_{1}}{\dsp\partial y} \frac{\dsp\partial\FoncBaseLocale_{2}}{\dsp\partial x} \right)\\
\MAT{K}^e_{5,1}(T_k)
&=|T_k|\left({\left(\lambda + 2 \, \mu\right)} \frac{\dsp\partial\FoncBaseLocale_{1}}{\dsp\partial x}
\frac{\dsp\partial\FoncBaseLocale_{3}}{\dsp\partial x}
+ \mu \frac{\dsp\partial\FoncBaseLocale_{1}}{\dsp\partial y} \frac{\dsp\partial\FoncBaseLocale_{3}}{\dsp\partial y} \right), 
&\hspace{-3mm}\MAT{K}^e_{6,1}(T_k)
&=|T_k|\left(\lambda \frac{\dsp\partial\FoncBaseLocale_{1}}{\dsp\partial x} \frac{\dsp\partial\FoncBaseLocale_{3}}{\dsp\partial y}
+ \mu \frac{\dsp\partial\FoncBaseLocale_{1}}{\dsp\partial y} \frac{\dsp\partial\FoncBaseLocale_{3}}{\dsp\partial x} \right)
\end{align*}
}
This gives, on the triangle $T_k,$
{\small
\begin{eqnarray*}
\MAT{K}^e_{1,1}(T_k)&=&|T_k|\left({\left(\lambda + 2 \, \mu\right)} G_1(1,k)^{2} + \mu G_1(2,k)^{2} \right) \\
\MAT{K}^e_{2,1}(T_k)&=&|T_k|\left(\lambda G_1(1,k) G_1(2,k) + \mu G_1(1,k) G_1(2,k) \right)\\ 
\MAT{K}^e_{3,1}(T_k)&=&|T_k|\left({\left(\lambda + 2 \, \mu\right)} G_1(1,k) G_2(1,k) + \mu G_1(2,k) G_2(2,k) \right)\\
\MAT{K}^e_{4,1}(T_k)&=&|T_k|\left(\lambda G_1(1,k) G_2(2,k) + \mu G_1(2,k) G_2(1,k) \right)\\
\MAT{K}^e_{5,1}(T_k)&=&|T_k|\left({\left(\lambda + 2 \, \mu\right)} G_1(1,k) G_3(1,k) + \mu G_1(2,k) G_3(2,k) \right)\\
\MAT{K}^e_{6,1}(T_k)&=&|T_k|\left(\lambda G_1(1,k) G_3(2,k) + \mu G_1(2,k) G_3(1,k) \right)
\end{eqnarray*}
}
Thus, the computation of the first six lines of $\MAT{K}_g$ may be vectorized under the form:
\begin{center}
\begin{minipage}[t]{1\textwidth}
\MonMatlabNonumber
\begin{lstlisting}
Kg(1,:)=((lam+2*mu)*G1(1,:).^2 + mu*G1(2,:).^2).*area;
Kg(2,:)=(lam*G1(1,:).*G1(2,:) + mu*G1(1,:).*G1(2,:)).*area;
Kg(3,:)=((lam+2*mu)*G1(1,:).*G2(1,:) + mu*G1(2,:).*G2(2,:)).*area;
Kg(4,:)=(lam*G1(1,:).*G2(2,:) + mu*G1(2,:).*G2(1,:)).*area;
Kg(5,:)=((lam+2*mu)*G1(1,:).*G3(1,:) + mu*G1(2,:).*G3(2,:)).*area;
Kg(6,:)=(lam*G1(1,:).*G3(2,:) + mu*G1(2,:).*G3(1,:)).*area;
\end{lstlisting}
\end{minipage}
\end{center}
The other columns of $\MAT{K}_g$ are computed on the same principle, using the symmetry of the matrix.
We give in the Listings \ref{listing:ElemStiffElasMatBaVecP1} and \ref{listing:StiffElasAssemblingP1OptV2}
the complete vectorized Matlab/Octave functions for computing $\MAT{K}_g$ and the elastic stiffness matrix assembly
respectively.
\begin{center}
\begin{minipage}[t]{0.95\textwidth}
\MonMatlab
\begin{lstlisting}[caption=Vectorized code for computing $\MAT{K}_g$,label=listing:ElemStiffElasMatBaVecP1] 
function [Kg]=ElemStiffElasMatVecP1(q,me,areas,lam,mu)
u=q(:,me(2,:))-q(:,me(3,:));  % q2-q3
G1=[u(2,:);-u(1,:)];
u=q(:,me(3,:))-q(:,me(1,:));  % q3-q1
G2=[u(2,:);-u(1,:)];
u=q(:,me(1,:))-q(:,me(2,:));  % q1-q2
G3=[u(2,:);-u(1,:)];
clear u
coef=ones(2,1)*(0.5./sqrt(areas));
G1=G1.*coef;
G2=G2.*coef;
G3=G3.*coef;
clear coef
Kg=zeros(36,size(me,2));
Kg(1,:)=(lam + 2*mu)*G1(1,:).^2 + mu*G1(2,:).^2;
Kg(2,:)=lam.*G1(1,:).*G1(2,:) + mu*G1(1,:).*G1(2,:);
Kg(3,:)=(lam + 2*mu)*G1(1,:).*G2(1,:) + mu*G1(2,:).*G2(2,:);
Kg(4,:)=lam.*G1(1,:).*G2(2,:) + mu*G1(2,:).*G2(1,:);
Kg(5,:)=(lam + 2*mu)*G1(1,:).*G3(1,:) + mu*G1(2,:).*G3(2,:);
Kg(6,:)=lam.*G1(1,:).*G3(2,:) + mu*G1(2,:).*G3(1,:);
Kg(8,:)=(lam + 2*mu)*G1(2,:).^2 + mu*G1(1,:).^2;
Kg(9,:)=lam.*G1(2,:).*G2(1,:) + mu*G1(1,:).*G2(2,:);
Kg(10,:)=(lam + 2*mu)*G1(2,:).*G2(2,:) + mu*G1(1,:).*G2(1,:);
Kg(11,:)=lam.*G1(2,:).*G3(1,:) + mu*G1(1,:).*G3(2,:);
Kg(12,:)=(lam + 2*mu)*G1(2,:).*G3(2,:) + mu*G1(1,:).*G3(1,:);
Kg(15,:)=(lam + 2*mu)*G2(1,:).^2 + mu*G2(2,:).^2;
Kg(16,:)=lam.*G2(1,:).*G2(2,:) + mu*G2(1,:).*G2(2,:);
Kg(17,:)=(lam + 2*mu)*G2(1,:).*G3(1,:) + mu*G2(2,:).*G3(2,:);
Kg(18,:)=lam.*G2(1,:).*G3(2,:) + mu*G2(2,:).*G3(1,:);
Kg(22,:)=(lam + 2*mu)*G2(2,:).^2 + mu*G2(1,:).^2;
Kg(23,:)=lam.*G2(2,:).*G3(1,:) + mu*G2(1,:).*G3(2,:);
Kg(24,:)=(lam + 2*mu)*G2(2,:).*G3(2,:) + mu*G2(1,:).*G3(1,:);
Kg(29,:)=(lam + 2*mu)*G3(1,:).^2 + mu*G3(2,:).^2;
Kg(30,:)=lam.*G3(1,:).*G3(2,:) + mu*G3(1,:).*G3(2,:);
Kg(36,:)=(lam + 2*mu)*G3(2,:).^2 + mu*G3(1,:).^2;
Kg([7,13,14,19,20,21,25,26,27,28,31,32,33,34,35],:)= ...
  Kg([2,3,9,4,10,16,5,11,17,23,6,12,18,24,30],:);
\end{lstlisting}
\end{minipage}
\end{center}
\begin{center}
\begin{minipage}[t]{0.98\textwidth}
\MonMatlab
\begin{lstlisting}[caption=Optimized matrix assembly code - version 2 (elastic stiffness matrix),label=listing:StiffElasAssemblingP1OptV2] 
function [K]=StiffElasAssemblingP1OptV2(nq,nme,q,me,areas,lam,mu)
[Ig,Jg]=BuildIgJgP1VF(me);
Kg=ElemStiffElasMatVecP1(q,me,areas,lam,mu);
K = sparse(Ig(:),Jg(:),Kg(:),2*nq,2*nq);
\end{lstlisting}
\end{minipage}
\end{center}
On Figure~\ref{VFOptV2VsFreeFEM}, we show the computation times of the \texttt{OptV2} 
(in Matlab/Octave) and FreeFEM++ codes, versus the number of degrees of freedom on the mesh.
The computation times values are given in Table~\ref{StiffVFMat}. 
The complexity of the Matlab/Octave codes is still linear ($\GrandO(\ndf)$). Moreover, the computation times
are $10$ (resp. $5$) times faster with Octave (resp. Matlab) than those obtained with FreeFEM++. 
\vspace{-2mm}
\begin{table}[H]\scriptsize
\begin{center}
\begin{tabular}{|r|r||*{3}{c|}}
  \hline 
  $n_q$ &  $n_{df}$ &  \begin{tabular}{c} Octave \\ (3.6.3) \end{tabular} &  \begin{tabular}{c} Matlab \\ (R2012b) \end{tabular} &  \begin{tabular}{c} FreeFEM++ \\ (3.20) \end{tabular}  \\ \hline \hline
$14222$ & $28444$ & \begin{tabular}{c} 0.088 (s)\\ \texttt{x 1.00} \end{tabular} & \begin{tabular}{c} 0.197 (s)\\ \texttt{x 0.45} \end{tabular} & \begin{tabular}{c} 1.260 (s)\\ \texttt{x 0.07} \end{tabular}\\ \hline
$55919$ & $111838$ & \begin{tabular}{c} 0.428 (s)\\ \texttt{x 1.00} \end{tabular} & \begin{tabular}{c} 0.769 (s)\\ \texttt{x 0.56} \end{tabular} & \begin{tabular}{c} 4.970 (s)\\ \texttt{x 0.09} \end{tabular}\\ \hline
$125010$ & $250020$ & \begin{tabular}{c} 0.997 (s)\\ \texttt{x 1.00} \end{tabular} & \begin{tabular}{c} 1.757 (s)\\ \texttt{x 0.57} \end{tabular} & \begin{tabular}{c} 11.190 (s)\\ \texttt{x 0.09} \end{tabular}\\ \hline
$225547$ & $451094$ & \begin{tabular}{c} 1.849 (s)\\ \texttt{x 1.00} \end{tabular} & \begin{tabular}{c} 3.221 (s)\\ \texttt{x 0.57} \end{tabular} & \begin{tabular}{c} 20.230 (s)\\ \texttt{x 0.09} \end{tabular}\\ \hline
$343082$ & $686164$ & \begin{tabular}{c} 2.862 (s)\\ \texttt{x 1.00} \end{tabular} & \begin{tabular}{c} 5.102 (s)\\ \texttt{x 0.56} \end{tabular} & \begin{tabular}{c} 30.840 (s)\\ \texttt{x 0.09} \end{tabular}\\ \hline
$506706$ & $1013412$ & \begin{tabular}{c} 4.304 (s)\\ \texttt{x 1.00} \end{tabular} & \begin{tabular}{c} 7.728 (s)\\ \texttt{x 0.56} \end{tabular} & \begin{tabular}{c} 45.930 (s)\\ \texttt{x 0.09} \end{tabular}\\ \hline
$689716$ & $1379432$ & \begin{tabular}{c} 5.865 (s)\\ \texttt{x 1.00} \end{tabular} & \begin{tabular}{c} 10.619 (s)\\ \texttt{x 0.55} \end{tabular} & \begin{tabular}{c} 62.170 (s)\\ \texttt{x 0.09} \end{tabular}\\ \hline
$885521$ & $1771042$ & \begin{tabular}{c} 8.059 (s)\\ \texttt{x 1.00} \end{tabular} & \begin{tabular}{c} 13.541 (s)\\ \texttt{x 0.60} \end{tabular} & \begin{tabular}{c} 79.910 (s)\\ \texttt{x 0.10} \end{tabular}\\ \hline
$1127090$ & $2254180$ & \begin{tabular}{c} 9.764 (s)\\ \texttt{x 1.00} \end{tabular} & \begin{tabular}{c} 17.656 (s)\\ \texttt{x 0.55} \end{tabular} & \begin{tabular}{c} 101.730 (s)\\ \texttt{x 0.10} \end{tabular}\\ \hline
$1401129$ & $2802258$ & \begin{tabular}{c} 12.893 (s)\\ \texttt{x 1.00} \end{tabular} & \begin{tabular}{c} 22.862 (s)\\ \texttt{x 0.56} \end{tabular} & \begin{tabular}{c} 126.470 (s)\\ \texttt{x 0.10} \end{tabular}\\ \hline
\end{tabular}
\end{center}
\caption{\label{StiffVFMat}Computational cost of the \texttt{StiffElas} matrix assembly versus $\nq/n_{df},$ with the \texttt{OptV2} Matlab/Octave codes (columns 3,4) and with FreeFEM++ (column 5) : time in seconds (top value) and speedup (bottom value). The speedup reference is \texttt{OptV2} Octave version.}
\end{table}

\begin{figure}[htp]
\dblimageps{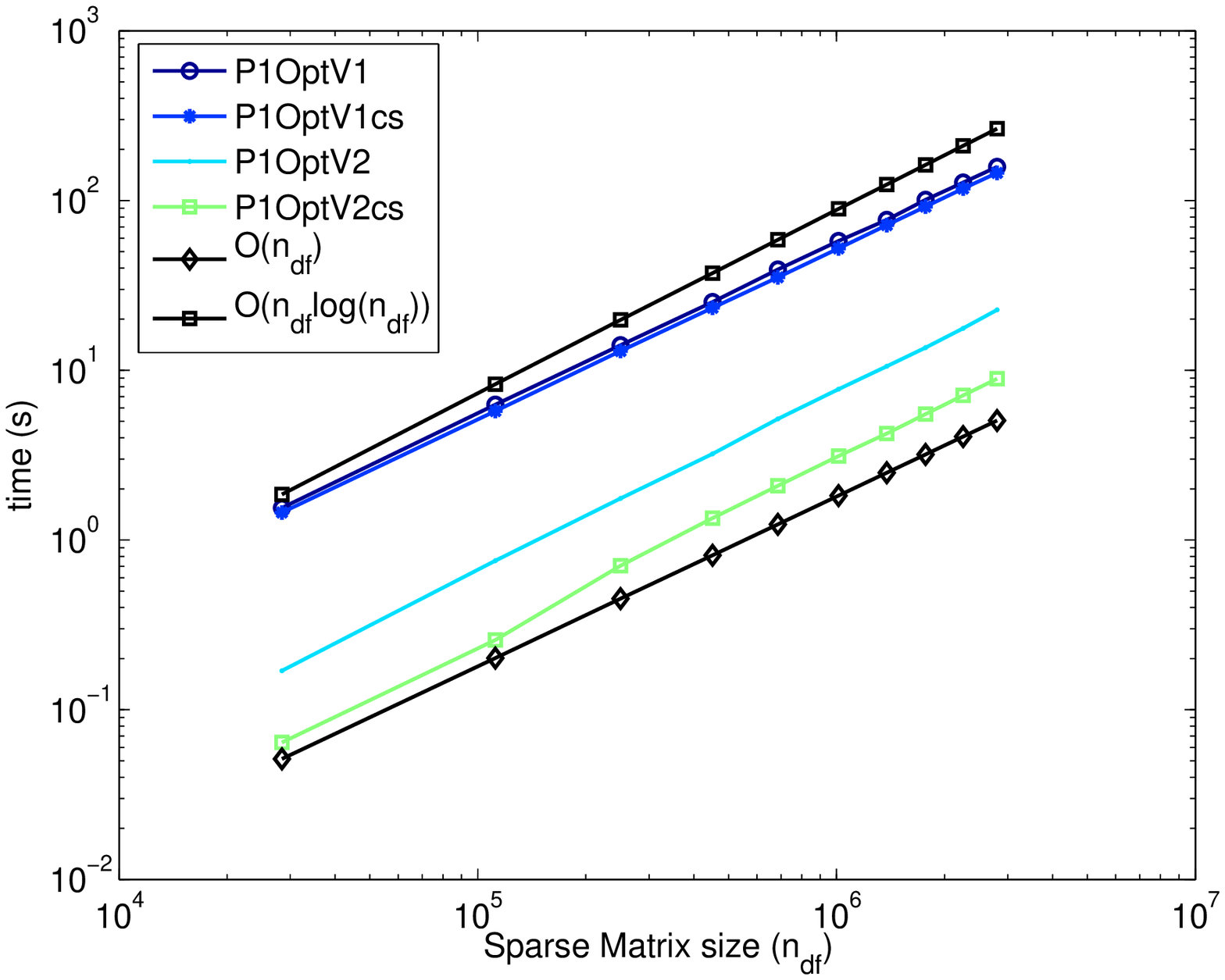}{\zoom}
{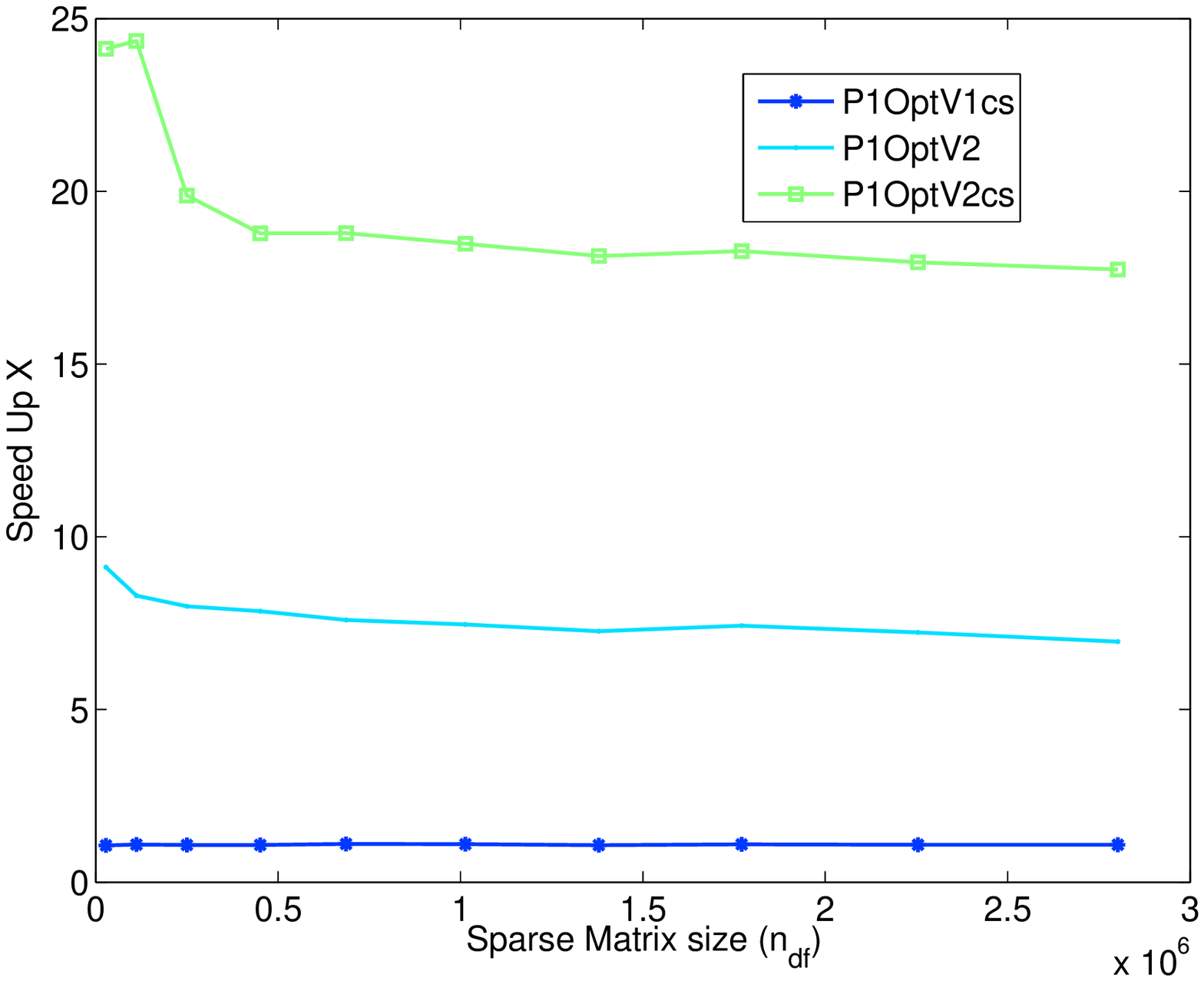}{\zoom}
\caption{Computational cost of the StiffElasAssembling functions versus $\ndf$, with Matlab (R2012b) : time
in seconds (left) and speedup (right). The speedup reference is \texttt{OptV1} version.
 \label{VFOptV2CS}}
\end{figure}

\vspace{-3mm}
\noindent
As observed in Remark \ref{rem:JIT}, Octave gives better results than Matlab only
for the \texttt{OptV2} codes. Using \texttt{cs\_sparse} function instead of
Matlab \texttt{sparse} function is approximately 1.1 (resp. 2.5) times faster for \texttt{OptV1}
(resp.  \texttt{OptV2}) version as shown on Figure~\ref{VFOptV2CS}.

%

\section{Conclusion}
For several examples of matrices, from the classical code we have built step by step
the codes to perform the assembly of these matrices to obtain a fully vectorized form.
For each version, we have described the algorithm and estimated its  numerical complexity.
The assembly of the mass, weighted mass and stiffness matrices of size~$10^6,$ on our reference computer,
is obtained in less than $4$ seconds (resp. about $2$ seconds) with Matlab (resp. with Octave).
The assembly of the elastic stiffness matrix of size $10^6,$ is computed in less than $8$ seconds
(resp. about $4$ seconds) with Matlab (resp. with Octave).

These optimization techniques in Matlab/Octave may be extended to other types of matrices,
for higher order or others finite elements ($P_k,$ $Q_k,$ ...) and in $3$D.

In Matlab, it is possible to further improve
the performances of the \texttt{OptV2} codes by using  Nvidia GPU cards. Preliminary Matlab results  
give a computation time divided by a factor $5$ on a Nvidia \texttt{GTX 590} GPU card (compared to the \texttt{OptV2} without GPU).


\appendix

\section{Comparison of the performances with FreeFEM++}

\subsection{Classical matrix assembly code vs FreeFEM++}
\label{App:basevsFreeFem}
\noindent

\begin{table}[H]\scriptsize
\begin{center}

\end{center}
\caption{Computational cost of the \texttt{Mass} matrix assembly versus $\nq,$ with the \texttt{basic} Matlab/Octave version (columns $2,3$) and with FreeFEM++ (column $4$) : 
time in seconds (top value) and speedup (bottom value). The speedup reference is \texttt{basic} Matlab version.\label{tab:Mass:base}}
\end{table}

\begin{table}[H]\scriptsize
\begin{center}
%
\end{center}
\caption{Computational cost of the \texttt{MassW} matrix assembly versus $\nq,$ with the \texttt{basic} Matlab/Octave version (columns $2,3$) and with FreeFEM++ (column $4$) : 
time in seconds (top value) and speedup (bottom value). The speedup reference is \texttt{basic} Matlab version.\label{tab:MassW:base}}
\end{table}

\begin{table}[H]\scriptsize
\begin{center}
%
\end{center}
\caption{Computational cost of the \texttt{Stiff} matrix assembly versus $\nq,$ with the \texttt{basic} Matlab/Octave version (columns $2,3$) and with FreeFEM++ (column $4$) : 
time in seconds (top value) and speedup (bottom value). The speedup reference is \texttt{basic} Matlab version.\label{tab:Stiff:base}}
\end{table}

\newpage
\subsection{OptV0 matrix assembly code vs FreeFEM++}\label{App:OptV0vsFreeFem}
\noindent

\begin{table}[H]\scriptsize
\begin{center}
%
\end{center}
\caption{Computational cost of the \texttt{Mass} matrix assembly versus $\nq,$ with the \texttt{OptV0} Matlab/Octave version (columns $2,3$) and with FreeFEM++ (column $4$) : 
time in seconds (top value) and speedup (bottom value). The speedup reference is \texttt{OptV0} Matlab version.\label{tab:Mass:OptV0}}
\end{table}

\begin{table}[H]\scriptsize
\begin{center}
%
\end{center}
\caption{Computational cost of the \texttt{MassW} matrix assembly versus $\nq,$ with the \texttt{OptV0} Matlab/Octave version (columns $2,3$) and with FreeFEM++ (column $4$) : 
time in seconds (top value) and speedup (bottom value). The speedup reference is \texttt{OptV0} Matlab version.\label{tab:MassW:OptV0}}
\end{table}

\begin{table}[H]\scriptsize
\begin{center}
%
\end{center}
\caption{Computational cost of the \texttt{Stiff} matrix assembly versus $\nq,$ with the \texttt{OptV0} Matlab/Octave version (columns $2,3$) and with FreeFEM++ (column $4$) : 
time in seconds (top value) and speedup (bottom value). The speedup reference is OptV0 Matlab version.\label{tab:Stiff:OptV0}}
\end{table}

\vspace{-0.95cm}
\subsection{OptV1 matrix assembly code vs FreeFEM++}\label{App:OptV1vsFreeFem}
\noindent

\begin{table}[H]\scriptsize
\begin{center}
%
\end{center}
\caption{Computational cost of the \texttt{Mass} matrix assembly versus $\nq,$ with the \texttt{OptV1} Matlab/Octave version (columns $2,3$) and with FreeFEM++ (column $4$) : 
time in seconds (top value) and speedup (bottom value). The speedup reference is \texttt{OptV1} Matlab version.}
\end{table}

\begin{table}[H]\scriptsize
\begin{center}
%
\end{center}
\caption{Computational cost of the \texttt{MassW} matrix assembly versus $\nq,$ with the \texttt{OptV1} Matlab/Octave version (columns $2,3$) and with FreeFEM++ (column $4$) : 
time in seconds (top value) and speedup (bottom value). The speedup reference is \texttt{OptV1} Matlab version.}
\end{table}

\begin{table}[H]\scriptsize
\begin{center}
%
\end{center}
\caption{Computational cost of the \texttt{Stiff} matrix assembly versus $\nq,$ with the \texttt{OptV1} Matlab/Octave version (columns $2,3$) and with FreeFEM++ (column $4$) : 
time in seconds (top value) and speedup (bottom value). The speedup reference is \texttt{OptV1} Matlab version.}
\end{table}

\newpage
\subsection{OptV2 matrix assembly code vs FreeFEM++}\label{App:OptV2vsFreeFem}
\noindent

\begin{table}[H]\scriptsize
\begin{center}
%
\end{center}
\caption{Computational cost of the \texttt{Mass} matrix assembly versus $\nq,$ with the \texttt{OptV2} Matlab/Octave codes (columns $2,3$) and with FreeFEM++ (column $4$) : time in seconds (top value) and speedup (bottom value). The speedup reference is \texttt{OptV2} Octave version.}
\end{table}

\vspace{-5mm}
\begin{table}[H]\scriptsize
\begin{center}
%
\end{center}
\caption{Computational cost of the \texttt{MassW} matrix assembly versus $\nq,$ with the \texttt{OptV2} Matlab/Octave codes (columns $2,3$) and with FreeFEM++ (column $4$) : time in seconds (top value) and speedup (bottom value). The speedup reference is \texttt{OptV2} Octave version.}
\end{table}

\begin{table}[H]\scriptsize
\begin{center}
%
\\ \hline
\end{tabular}
\end{center}
\caption{Computational cost of the \texttt{Stiff} matrix assembly versus $\nq,$ with the \texttt{OptV2} Matlab/Octave codes (columns $2,3$) and with FreeFEM++ (column $4$) : time in seconds (top value) and speedup (bottom value). The speedup reference is \texttt{OptV2} Octave version.}
\end{table}

\section{Matrix assembly codes}\label{App:codes}
\subsection{Element matrices}\label{App:ElemMat}
\noindent
\begin{lstlisting}[caption=ElemMassMatP1.m]
function AElem=ElemMassMatP1(area)
AElem=(area/12)*[2 1 1; 1 2 1; 1 1 2];
\end{lstlisting}
\begin{lstlisting}[caption=ElemMassWMatP1.m]
function AElem=ElemMassWMatP1(area,w)
AElem=(area/30)* ...
          [3*w(1)+w(2)+w(3), w(1)+w(2)+w(3)/2, w(1)+w(2)/2+w(3); ...
           w(1)+w(2)+w(3)/2, w(1)+3*w(2)+w(3), w(1)/2+w(2)+w(3); ...
           w(1)+w(2)/2+w(3), w(1)/2+w(2)+w(3), w(1)+w(2)+3*w(3)];
\end{lstlisting}
\begin{lstlisting}[caption=ElemStiffMatP1.m]
function AElem=ElemStiffMatP1(q1,q2,q3,area)
M=[q2-q3, q3-q1, q1-q2];
AElem=(1/(4*area))*M'*M;
\end{lstlisting}
\newpage
\subsection{Classical matrix assembly code}\label{App:codebase}
\noindent
\vspace{0.3cm}
\begin{lstlisting}[caption=MassAssemblingP1base.m]
function M=MassAssemblingP1base(nq,nme,me,areas)
M=sparse(nq,nq);
for k=1:nme
    E=ElemMassMatP1(areas(k));
    for il=1:3
        i=me(il,k);
        for jl=1:3
            j=me(jl,k);
            M(i,j)=M(i,j)+E(il,jl);
        end
    end
end
\end{lstlisting}
\begin{lstlisting}[caption=MassWAssemblingP1base.m]
function M=MassWAssemblingP1base(nq,nme,me,areas,Tw)
M=sparse(nq,nq);
for k=1:nme
    for il=1:3
        i=me(il,k);
	Twloc(il)=Tw(i);
    end
    E=ElemMassWMatP1(areas(k),Twloc);
    for il=1:3
        i=me(il,k);
        for jl=1:3
            j=me(jl,k);
            M(i,j)=M(i,j)+E(il,jl);
        end
    end
end
\end{lstlisting}
\begin{lstlisting}[caption=StiffAssemblingP1base.m]
function R=StiffAssemblingP1base(nq,nme,q,me,areas)
R=sparse(nq,nq);
for k=1:nme
    E=ElemStiffMatP1(q(:,me(1,k)),q(:,me(2,k)),q(:,me(3,k)),areas(k));
    for il=1:3
        i=me(il,k);
        for jl=1:3
            j=me(jl,k);
            R(i,j)=R(i,j)+E(il,jl);
        end
    end
end
\end{lstlisting}
\subsection{Optimized matrix assembly codes - Version 0}\label{App:codeV0}
\noindent
\vspace{0.3cm}
\begin{lstlisting}[caption=MassAssemblingP1OptV0.m]
function M=MassAssemblingP1OptV0(nq,nme,me,areas)
M=sparse(nq,nq);
for k=1:nme
  I=me(:,k);
  M(I,I)=M(I,I)+ElemMassMatP1(areas(k));
end
\end{lstlisting}
\begin{lstlisting}[caption=MassWAssemblingP1OptV0.m]
function M=MassWAssemblingP1OptV0(nq,nme,me,areas,Tw)
M=sparse(nq,nq);
for k=1:nme
  I=me(:,k);
  M(I,I)=M(I,I)+ElemMassWMatP1(areas(k),Tw(me(:,k)));
end
\end{lstlisting}
\begin{lstlisting}[caption=StiffAssemblingP1OptV0.m]
function R=StiffAssemblingP1OptV0(nq,nme,q,me,areas)
R=sparse(nq,nq);
for k=1:nme
  I=me(:,k);
  Me=ElemStiffMatP1(q(:,me(1,k)),q(:,me(2,k)),q(:,me(3,k)),areas(k));
  R(I,I)=R(I,I)+Me;
end
\end{lstlisting}
\vspace{.5cm}
\subsection{Optimized matrix assembly codes - Version 1}\label{App:codeV1}
\noindent
\vspace{0.3cm}
\begin{lstlisting}[caption=MassAssemblingP1OptV1.m]
function M=MassAssemblingP1OptV1(nq,nme,me,areas)
Ig=zeros(9*nme,1);Jg=zeros(9*nme,1);Kg=zeros(9*nme,1);

ii=[1 2 3 1 2 3 1 2 3];
jj=[1 1 1 2 2 2 3 3 3];
kk=1:9;
for k=1:nme
  E=ElemMassMatP1(areas(k));
  Ig(kk)=me(ii,k); 
  Jg(kk)=me(jj,k);
  Kg(kk)=E(:);
  kk=kk+9;
end
M=sparse(Ig,Jg,Kg,nq,nq);
\end{lstlisting}
\newpage
\begin{lstlisting}[caption=MassWAssemblingP1OptV1.m]
function M=MassWAssemblingP1OptV1(nq,nme,me,areas,Tw)
Ig=zeros(9*nme,1);Jg=zeros(9*nme,1);Kg=zeros(9*nme,1);

ii=[1 2 3 1 2 3 1 2 3];
jj=[1 1 1 2 2 2 3 3 3];
kk=1:9;
for k=1:nme
  E=ElemMassWMat(areas(k),Tw(me(:,k)));
  Ig(kk)=me(ii,k); 
  Jg(kk)=me(jj,k);
  Kg(kk)=E(:);
  kk=kk+9;
end
M=sparse(Ig,Jg,Kg,nq,nq);
\end{lstlisting}
\begin{lstlisting}[caption=StiffAssemblingP1OptV1.m]
function R=StiffAssemblingP1OptV1(nq,nme,q,me,areas)
Ig=zeros(nme*9,1);Jg=zeros(nme*9,1);
Kg=zeros(nme*9,1);

ii=[1 2 3 1 2 3 1 2 3];
jj=[1 1 1 2 2 2 3 3 3];
kk=1:9;
for k=1:nme
  Me=ElemStiffMatP1(q(:,me(1,k)),q(:,me(2,k)),q(:,me(3,k)),areas(k));
  Ig(kk)=me(ii,k); 
  Jg(kk)=me(jj,k);
  Kg(kk)=Me(:);
  kk=kk+9;
end
R=sparse(Ig,Jg,Kg,nq,nq);
\end{lstlisting}

\section{Matlab sparse trouble\label{MatlabTrouble}}\label{App:MatlabTrouble}
In this part, we illustrate a problem that we encountered in the
development of our codes : decrease of the performances of the
assembly codes, for the classical and \texttt{OptV0} versions, when
migrating from release R2011b to release R2012a or R2012b
independently of the operating system used.  In fact, this comes
from the use of the command \lstinline{M = sparse(nq,nq)}.  We
illustrate this for the mass matrix assembly, by giving in
Table~\ref{BenchMatlab} the computation time of the function
\lstinline{MassAssemblingP1OptV0} for different Matlab releases.

This problem has been reported to the MathWorks's development team :
\begin{verbatim}
As you have correctly pointed out, MATLAB 8.0 (R2012b) seems to perform
slower than the previous releases for this specific case of 
reallocation in sparse matrices. I will convey this information to the
development team for further investigation and a possible fix in the 
future releases of MATLAB. I apologize for the inconvenience.
\end{verbatim}
\vspace{3mm}
To fix this issue in the releases R2012a and R2012b, it is recommended by the Matlab's technical support
to use the function \lstinline{spalloc} instead of the function \lstinline{sparse} :
\begin{verbatim}
The matrix, 'M' in the function 'MassAssemblingP1OptV0' was initialized 
using SPALLOC instead of the SPARSE command since the maximum number of
non-zeros in M are already known. 

Previously existing line of code: 
M = sparse(nq,nq);

Modified line of code: 
M = spalloc(nq, nq, 9*nme);
\end{verbatim}

\vspace{-1mm}
\begin{table}[htbp]\scriptsize
\begin{center}
\begin{tabular}{|r||*{4}{c|}}
  \hline 
  Sparse dim &  R2012b &  R2012a &  R2011b &  R2011a  \\ \hline \hline
$1600$ & \begin{tabular}{c} 0.167 (s) \\ $(\times 1.00)$\end{tabular} & \begin{tabular}{c} 0.155 (s) \\ $(\times 1.07)$\end{tabular} & \begin{tabular}{c} 0.139 (s) \\ $(\times 1.20)$\end{tabular} & \begin{tabular}{c} 0.116 (s) \\ $(\times 1.44)$\end{tabular}\\ \hline
$3600$ & \begin{tabular}{c} 0.557 (s) \\ $(\times 1.00)$\end{tabular} & \begin{tabular}{c} 0.510 (s) \\ $(\times 1.09)$\end{tabular} & \begin{tabular}{c} 0.461 (s) \\ $(\times 1.21)$\end{tabular} & \begin{tabular}{c} 0.355 (s) \\ $(\times 1.57)$\end{tabular}\\ \hline
$6400$ & \begin{tabular}{c} 1.406 (s) \\ $(\times 1.00)$\end{tabular} & \begin{tabular}{c} 1.278 (s) \\ $(\times 1.10)$\end{tabular} & \begin{tabular}{c} 1.150 (s) \\ $(\times 1.22)$\end{tabular} & \begin{tabular}{c} 0.843 (s) \\ $(\times 1.67)$\end{tabular}\\ \hline
$10000$ & \begin{tabular}{c} 4.034 (s) \\ $(\times 1.00)$\end{tabular} & \begin{tabular}{c} 2.761 (s) \\ $(\times 1.46)$\end{tabular} & \begin{tabular}{c} 1.995 (s) \\ $(\times 2.02)$\end{tabular} & \begin{tabular}{c} 1.767 (s) \\ $(\times 2.28)$\end{tabular}\\ \hline
$14400$ & \begin{tabular}{c} 8.545 (s) \\ $(\times 1.00)$\end{tabular} & \begin{tabular}{c} 6.625 (s) \\ $(\times 1.29)$\end{tabular} & \begin{tabular}{c} 3.734 (s) \\ $(\times 2.29)$\end{tabular} & \begin{tabular}{c} 3.295 (s) \\ $(\times 2.59)$\end{tabular}\\ \hline
$19600$ & \begin{tabular}{c} 16.643 (s) \\ $(\times 1.00)$\end{tabular} & \begin{tabular}{c} 13.586 (s) \\ $(\times 1.22)$\end{tabular} & \begin{tabular}{c} 6.908 (s) \\ $(\times 2.41)$\end{tabular} & \begin{tabular}{c} 6.935 (s) \\ $(\times 2.40)$\end{tabular}\\ \hline
$25600$ & \begin{tabular}{c} 29.489 (s) \\ $(\times 1.00)$\end{tabular} & \begin{tabular}{c} 27.815 (s) \\ $(\times 1.06)$\end{tabular} & \begin{tabular}{c} 12.367 (s) \\ $(\times 2.38)$\end{tabular} & \begin{tabular}{c} 11.175 (s) \\ $(\times 2.64)$\end{tabular}\\ \hline
$32400$ & \begin{tabular}{c} 47.478 (s) \\ $(\times 1.00)$\end{tabular} & \begin{tabular}{c} 47.037 (s) \\ $(\times 1.01)$\end{tabular} & \begin{tabular}{c} 18.457 (s) \\ $(\times 2.57)$\end{tabular} & \begin{tabular}{c} 16.825 (s) \\ $(\times 2.82)$\end{tabular}\\ \hline
$40000$ & \begin{tabular}{c} 73.662 (s) \\ $(\times 1.00)$\end{tabular} & \begin{tabular}{c} 74.188 (s) \\ $(\times 0.99)$\end{tabular} & \begin{tabular}{c} 27.753 (s) \\ $(\times 2.65)$\end{tabular} & \begin{tabular}{c} 25.012 (s) \\ $(\times 2.95)$\end{tabular}\\ \hline
\end{tabular}
\end{center}
\caption{MassAssemblingP1OptV0 for different Matlab releases : computation times and speedup\label{BenchMatlab}}
\end{table}


\vspace{-14mm}
\begin{thebibliography}{10} 
\bibitem{Bernardi:NNA:1989} {\sc C. Bernardi, Y. Maday, and A. Patera}, {\em A new nonconforming approach
  to domain decomposition: the mortar element method}, in Nonlinear Partial Differential Equations and their Applications,
  H. Brezis and J.L. Lions, eds., 1989.

\bibitem{Chen:PFE:2011}
{\sc L. Chen}, {\em Programming of Finite Element Methods in Matlab}, 
  Preprint, University of California Irvine, \url{http://math.uci.edu/~chenlong/226/Ch3FEMCode.pdf}, 2011.

\bibitem{Chen:iFEM:2013}
{\sc L. Chen}, {\em iFEM, a Matlab software package}, 
  University of California Irvine, \url{http://math.uci.edu/~chenlong/programming.html}, 2013. 

\bibitem{Ciarlet:TFE:2002}
{\sc P.~G. Ciarlet}, {\em The finite element method for elliptic problems}, 
  SIAM, Philadelphia, 2002.
 
\bibitem{Cuvelier:MEF:2008}
{\sc F. Cuvelier}, {\em Méthodes des éléments finis. {D}e la théorie à la programmation}, 
  Lecture Notes, Université Paris 13, \url{http://www.math.univ-paris13.fr/~cuvelier/docs/poly/polyFEM2D.pdf}, 2008.

\bibitem{Cuvelier:OptFEM2DP1:2012}
{\sc F. Cuvelier, C. Japhet, and G. Scarella}, {\em OptFEM2DP1, a MATLAB/Octave software package codes}, 
  Université Paris 13, \url{http://www.math.univ-paris13.fr/~cuvelier}, 2012.

\bibitem{Cuvelier:AEW:2013}
{\sc F. Cuvelier, C. Japhet, and G. Scarella}, {\em An efficient way to perform the assembly of finite element matrices
  in Matlab and Octave: the $P_k$ finite element case}, in preparation.

\bibitem{Davis:DMS:2006}
{\sc T.~A.~Davis}, {\em Direct Methods for Sparse Linear Systems}, SIAM, 2006.

\bibitem{Davis:SSP:2012}
{\sc T.~A. Davis}, {\em SuiteSparse packages, release 4.0.12}, 
  University of Florida, \url{http://www.cise.ufl.edu/research/sparse/SuiteSparse}, 2012.

\bibitem{Dhatt:FEM:2012}
{\sc G. Dhatt, E. Lefrançois, and G. Touzot}, {\em Finite Element Method}, Wiley, 2012.

\bibitem{Octave:2012}
{\em GNU Octave},  \url{http://www.gnu.org/software/octave}, 2012.

\bibitem{Hannukainen:IFE:2012}
{\sc A. Hannukainen, and M. Juntunen}, {\em Implementing the Finite Element Assembly in Interpreted Languages}, 
  Preprint, Aalto University, \url{http://users.tkk.fi/~mojuntun/preprints/matvecSISC.pdf}, 2012.

\bibitem{Hecht:FreeFEM:2012}
{\sc F. Hecht}, {\em FreeFEM++}, \url{http://www.freefem.org/ff++/index.htm}, 2012.

\bibitem{Johnson:NSP:2009}
{\sc C. Johnson}, {\em Numerical  Solution of Partial Differential Equations by the Finite Element Method}, 
  Dover Publications, Inc, 2009.

\bibitem{Lucquin:ISC:1998}
{\sc B.~Lucquin, and O.~Pironneau}, {\em Introduction to Scientific Computing}, 
  John Wiley \& Sons Ltd, 1998.

\bibitem{Matlab:2012}
{\em MATLAB}, \url{http://www.mathworks.com}, 2012.

\bibitem{Rahman:FMA:2011} {\sc T. Rahman, and J. Valdman}, {\em Fast MATLAB assembly of FEM matrices in 2D and 3D:
  Nodal elements}, Appl. Math. Comput., 219(13) (2013), pp.~7151--7158.

\end{thebibliography}
\end{document}